\documentclass[12pt]{article}
\usepackage{makeidx}
\usepackage{amssymb,amsthm,amsmath,amsfonts,enumerate}
\usepackage{float}
\usepackage{setspace,graphicx,rotating,caption,subcaption,relsize}
\usepackage{url}
\usepackage{hyperref}
\usepackage[authoryear,longnamesfirst]{natbib}
\usepackage{array}
\usepackage{indentfirst}

\hypersetup{
  colorlinks = true, urlcolor = cyan, linkcolor = blue, citecolor = red }
\renewenvironment{abstract}
 {\small
  \begin{center}
  \bfseries \abstractname\vspace{-.0em}\vspace{0pt}
  \end{center}
  \list{}{    \setlength{\leftmargin}{0mm}
    \setlength{\rightmargin}{\leftmargin}  }  \item\relax}
 {\endlist}
\makeatletter
\def\maketag@@@#1{\hbox{\m@th\normalfont\normalsize#1}}
\makeatother

\allowdisplaybreaks
\newtheorem {theorem}{Theorem}[section]
\newtheorem {assumption}{Assumption}

\newtheorem{example}{EXAMPLE}

\newtheorem{lemma}[theorem]{Lemma}

\newtheorem{remark}{Remark}[section]

\allowdisplaybreaks
\counterwithin*{equation}{section}
\numberwithin{equation}{section}
\begin{document}

\title{Max-laws of large numbers for weakly dependent\\
high dimensional arrays with applications}
\author{Jonathan B. Hill\thanks{%
Department of Economics, University of North Carolina, Chapel Hill, North
Carolina, E-mail:\texttt{jbhill@email.unc.edu}; \texttt{%
https://tarheels.live/jbhill}.}\medskip \\
Dept. of Economics, University of North Carolina, Chapel Hill, NC}
\date{{\large This draft:} \today
}
\maketitle

\begin{abstract}
We derive weak \textit{max-laws of large numbers} for $\max_{1\leq i\leq
k_{n}}|1/n\sum_{t=1}^{n}x_{i,n,t}|$ for zero mean stochastic triangular
arrays $\{x_{i,n,t}$ $:$ $1$ $\leq $ $t$ $\leq n\}_{n\geq 1}$, with
dimension counter $i$ $=$ $1,...,k_{n}$ and dimension $k_{n}$ $\rightarrow $ 
$\infty $. Rates of convergence are also analyzed based on feasible
sequences $\{k_{n}\}$. We work in three dependence settings: independence,
Dedecker and Prieur's (\citeyear{Dedecker_Prieur_2004}) $\tau $-mixing and
Wu's (\citeyear{Wu2005}) physical dependence. We initially ignore
cross-coordinate $i$ dependence as a benchmark. We then work with
martingale, nearly martingale, and mixing coordinates to deliver improved
bounds on $k_{n}$. Finally, we use the results in three applications, each
representing a key novelty: we ($i$) bound $k_{n}$\ for a max-correlation
statistic for regression residuals under $\alpha $-mixing or physical
dependence; ($ii$) extend correlation screening, or marginal regressions, to
physical dependent data with diverging dimension $k_{n}$ $\rightarrow $ $%
\infty $; and ($iii$) test a high dimensional parameter after partialling
out a fixed dimensional nuisance parameter in a linear time series
regression model under $\tau $-mixing.\medskip \newline
\end{abstract}

\setstretch{1.2}

\section{Introduction\label{sec:intro}}

In this article we derive and compare laws of large numbers for the maximum
sample mean of a triangular array $\{x_{n,t}$ $:$ $1$ $\leq $ $t$ $\leq $ $%
n\}_{n\geq 1}$ where $x_{n,t}$ $=$ $[x_{i,n,t}]_{i=1}^{k}$ $\in $ $\mathbb{R}%
^{k}$ with dimension $k$ $\in $ $\mathbb{N}$ and sample size $n$. When $k$ $%
= $ $k_{n}>>n$ we have a high dimensional [HD] setting that may be
potentially huge relative to the sample size (e.g. $\ln (k_{n})$ $\sim $ $%
an^{b}$ for some $a$, $b$ $>$ $0$, or $k_{n}$ $\rightarrow $ $\infty $
arbitrarily fast, depending on available information). We are particularly
interested in disparate settings of weak dependence and their impact on
feasible sequences $\{k_{n}\}$. High dimensionality is common due to the
enormous amount of available data, survey techniques, and technology for
data collection. Examples span social, communication, bio-genetic,
electrical, and engineering sciences to name a few. See, for instance, \cite%
{FanLi2006}, \cite{BuhlmannVanDeGeer2011}, \cite{FanLvQi2011}, and \cite%
{BelloniChernozhukovHansen2014} for examples and surveys. Our main results
are then applied to three settings in econometrics and statistics detailed
below.

Assuming $\mathbb{E}x_{n,t}$ $=$ $0$ for all $n$ and $t$, we derive what we
call a \textit{max-Weak Law of Large Numbers} (max-WLLN) for certain integer
sequences $\{k_{n}\}$ by case,%
\begin{equation}
\mathcal{M}_{n}:=\max_{1\leq i\leq k_{n}}\left\vert \frac{1}{n}%
\sum_{t=1}^{n}x_{i,n,t}\right\vert \overset{p}{\rightarrow }0.
\label{maxLLN}
\end{equation}%
Typically we obtain $\mathcal{M}_{n}$ $\overset{p}{\rightarrow }$ $0$ by
proving $\mathbb{E}|\mathcal{M}_{n}|^{p}$ $\rightarrow $ $0$ for $p$ $\geq $ 
$1$, and we establish $\{k_{n}\}$ such that $\sqrt{n}\mathcal{M}_{n}$ $=$ $%
O_{p}(g(k_{n}))$ for case-specific monotonic mappings $g$. We will call the
property $\mathbb{E}|\mathcal{M}_{n}|^{p}$ $\rightarrow $ $0$\ a max-WLLN
throughout as a convenience. In the technical appendix \cite{sm_max_LLN} we
derive a max-Strong Law of Large Numbers (max-SLLN) for independent
sequences, and reveal how the classic Cauchy sequence and Kronecker lemma
argument fails to deliver a max-SLLN under dependence.

Although max-laws are implicitly used in many papers too numerous to cite,
often under sub-exponential or sub-Gaussian tails and independence, we
believe this is the first attempt to derive and compare possible laws and
their resulting bounds on $k_{n}$ under various serial (over $t$) or
cross-coordinate (over $i$) dependence and heterogeneity settings. A very
few examples where max-WLLN's appear include HD model inference under
independence \citep{Dezeure_etal_2017,Hill_2025_maxtest} or weak dependence %
\citep[e.g.][]{Babii_etal2019,Adamek_et_al2023,MiesSteland2023}, and
wavelet-like HD covariance stationary tests under linearity %
\citep{JinWangWang2015,HillLi2025}. \cite{Hill_2025_maxtest} explores
max-LLN's for standard least squares components in an iid linear regression
setting. \cite{JinWangWang2015} exploit HD theory for autocovariances dating
to \citet[Chapt.
7]{HannanDiestler1988} and \cite{Keenan1997}. They require linearity with
iid innovations, and only work with high dimensionality across
autocovariance lags and so-called systematic samples (sub-sample counters). 
\cite{HillLi2025} work in the same setting under a broader dependence
concept. Thus neither systematically presents max-LLN's for heterogeneous
high dimensional arrays.

\cite{Adamek_et_al2023} develop inference methods for debiased Lasso in a
linear time series setting. Their Lemma A.4 presents an implicit max-WLLN by
using a union bound and mixingale maximal inequality (for sub-samples). That
result is quite close to what we present here. They require uniform $%
\mathcal{L}_{p}$-boundedness for some $p$ $>$ $2$, and near epoch dependence
[NED]. We allow for trending higher moments and $p$ $>$ $1$ under physical
dependence yielding a max-WLLN, while NED implies mixingale and adapted
mixingales are physical dependent \citep{Davidson1994,Hill2025_mixg}. We
also use cross-coordinate dependence to improve $k_{n}$. Thus our results
are more general and broad in scope. See Remark \ref{rm:Ademeck} for
details. The reader should also consult \cite{Babii_etal2019} for earlier
related work.

\cite{MiesSteland2023} exploit martingale theory in \cite{Pinelis1994} to
yield an $\mathcal{L}_{q}$-maximal inequality under $\mathcal{L}_{p}$%
-physical dependence, $2$ $\leq $ $p$ $\leq $ $q$. Their upper bound appears
sharper than the one we present in Lemma \ref{lm:phys_dep_Lp} and Theorem %
\ref{thm:max_LLN_phys_dep}, also based on a martingale approximation. The
improvement, however, does not yield a faster rate $k_{n}$ $\rightarrow $ $%
\infty $, while the latter can only be deduced once $p$ $=$ $q$. Moreover,
we allow for sub-exponential tails or $\mathcal{L}_{p}$-boundedness, $p$ $>$ 
$1$, and exploit cross-coordinate dependence, each new and ignored in \cite%
{MiesSteland2023}.

Apparently only max-WLLN's exist: a max-strong LLN (max-SLLN) has not been
explored. Moreover, max-LLN's are not explicitly available for $\tau $%
-mixing and physical dependent arrays under broad tail conditions, and to
the best our of knowledge inter-coordinate dependence is universally ignored
where union bounds, Lyapunov' inequality, and log-exp bounds under
sub-exponentiality are the standard for getting around $\max_{1\leq i\leq
k_{n}}|\cdot |$ and for bounding $k_{n}$.

We work under three broad dependence and heterogeneity settings: ($i$) $\tau 
$-mixing \citep{Dedecker_Prieur_2004}, ($ii$) $\mathcal{L}_{p}$-physical
dependence ($p>1$) \citep{Wu2005,WuMin2005}, and ($iii$) independence, with
varying degrees of cross-coordinate dependence. See Tables \ref{tlb:sum} and %
\ref{tlb:rank} and the end of the paper for an exhaustive summary of
allowed dependence and tail decay properties, with corresponding
cross-coordinate restrictions, resulting rates $k_{n}$ $\rightarrow $ $%
\infty $ and theorem references.

Under $(i)$, $(ii.a)$ and $(iii)$ in Table \ref{tlb:sum}\ we do not restrict
dependence coordinate-wise. This is the seemingly universal setting in the
high dimensional literatures. A variety of mixing and related properties
promote a Bernstein-type inequality that yield (\ref{maxLLN}) and bounds on $%
k_{n}$ qualitatively similar to the independence case. We treat a recent
representative sub-exponential $\tau $-mixing %
\citep{Dedecker_Prieur_2004,DedeckerPrieur2005}. The latter construction
along with other recent mixing concepts, like mixingale and related
moment-based constructions \citep{Gordin1969,McLeish1975}, were proposed to
handle stochastic processes that are not, e.g., uniform $\sigma $-field
based $\alpha $-, $\beta $-, or $\phi $-mixing. This includes possibly
infinite order functions of mixing processes, and Markovian dynamical
systems and related expanding maps, covering simple autoregressions with
Bernoulli shocks, and various attractors in mathematical physics with
applications in atmospheric mapping, electrical components and artificial
intelligence %
\citep[e.g.][]{Chernick1981,Andrews1984,Rio1996,Collet_etal2002,DedeckerPrieur2005,ChazottesGouezel2012}%
. Thus they fill certain key gaps in the field of processes that yield
deviation or concentration bounds and central limits.

We include ($ii.b$)-($ii.d$) to show that bounds on $k_{n}$ can be improved
when cross-coordinate dependence is available. We work under serial physical
dependence to focus ideas, but the result appears to apply generally. Strong
coordinate dependence ($ii.b$), where $x_{i,n,t}$ is a martingale over $i$,
yields unbounded $k_{n}$ (the result is truly \textit{dimension-agnostic}).
Under ($ii.c$) the condition is weakened such that $x_{i,n,t} $ \textit{%
becomes} a martingale as $n$ $\rightarrow $ $\infty $: $\mathbb{P}(\mathbb{E}%
[x_{i+1,n,t}|\mathfrak{F}_{i,n}]$ $=$ $x_{i,n,t})$ $\rightarrow $ $1$ for
some filtration $\{\mathfrak{F}_{i,n}\}$.\footnote{\textit{Nearly martingale}
in this paper distinctly differs from \textit{near-martingale} ($\mathbb{E}_{%
\mathfrak{F}_{i-1,t}}x_{i,n,t}$ $=$ $\mathbb{E}_{\mathfrak{F}%
_{i-1,t}}x_{i-1,n,t}$), \textit{weak-martingale}, or \textit{local-martingale%
} \citep[cf.][]{Kallenberg2021}.} We show that even in a Gaussian setting $%
k_{n}$ must be restricted, but a better bound is yielded by using
cross-coordinate information. We obtain the same result under
cross-coordinate mixing ($ii.d$) where improvements are gained in Gaussian,
sub-exponential and heavy-tailed cases.

As a third dependence setting $(iii)$ we deliver max-LLN's under serial
independence in the supplemental material \citet[Appendix B]{sm_max_LLN}. We
prove a max-SLLN under $\mathcal{L}_{1}$-boundedness and show that $k_{n}$
is unrestricted when a cross-coordinate probability decay property holds.
The proof exploits a new necessary and sufficient HD three-series theorem.

The cases are naturally nested: mixing includes independence, and physical
dependence covers mixing and non-mixing cases. Moreover, $\tau $-mixing and
adapted mixingale properties are closely related 
\citep[Appendix
C]{sm_max_LLN}, while adapted mixingale and physical dependence properties
are asymmetrically related \citep{Hill2025_mixg}. Mixingale-like constructs
date at least to \citet{Gordin1969}, \citet[eq.
(4)]{Hannan1973}, and \cite{McLeish1975}, with expansions to $\mathcal{L}%
_{p} $-arrays in, e.g., \cite{Andrews1988} and \cite{Hansen1991}. In the $%
\mathcal{L}_{p}$-physical dependence case if the coefficients grow in $p$ at
a polynomial rate then a Bernstein inequality promotes an exponential bound
on $k_{n}$.

Key technical tools, depending on the dependence property, are: \textit{%
log-exp} (or \textquotedblleft \textit{log-sum-exp}") bound on the maximum
of a sequence when a moment generating function exists; Bernstein,
Fuk-Naegev, and \cite{Nemirovski2000} inequalities; and maximal
inequalities, e.g. for physical dependent arrays. The log-exp transform
yields a \textquotedblleft smooth-max" approximation that has been broadly
exploited when cross-coordinate dependence is not modeled 
\citep[see,
e.g.,][]{Talagrand2003,BuhlmannVanDeGeer2011,Chernozhukov_etal2013}.

Bernstein-type inequalities exist for iid and various mixing and related
sequences, covering $\alpha $-, $\beta $-, $\phi $-, $\tilde{\phi}$-, $%
\varphi $-, $\tau $- and $\mathcal{C}$- mixing random variables in array,
random field and lattice forms %
\citep[e.g.][]{Rio1995,Samson2000,Merlevede_et_al_2011,HangSteinwart2017},
and physical dependent processes \citep{Wu2005}.\footnote{%
Consult, e.g., \cite{Dedecker_etal2007} for many mixing definitions, cf. 
\cite{Rio1996}, \citet{Dedecker_Prieur_2004,DedeckerPrieur2005}, and \cite%
{Maume-Deschamps2006}.} In most cases the random variables are assumed
bounded or sub-exponential, and in many cases only $1$-Lipschitz functions
are treated. We generalize the $\tau $-mixing $\mathcal{L}_{1}$ metric to an 
$\mathcal{L}_{p}$ metric, $p$ $\geq $ $1$, and derive a Bernstein inequality
under so-called $\tau ^{(p)}$-mixing by closely following \cite%
{Merlevede_et_al_2011}.

We do not attempt to use the sharpest available bounds within the
Bernstein-Hoeffding class, or under physical dependence. This is both for
clarity and ease of presenting proofs, and generally because sharp bounds
will only lead to modest, or no, improvements for $k_{n}$. See %
\citet{Talagrand1995a,Talagrand1995b}, \citet{Bentkus2008} and \cite%
{Dumbgen_etal_2010} for many results and suggested readings.

Bernstein and\ Fuk-Nagaev inequalities that can be used for max-LLN's have
been expanded beyond classic settings, covering bounded or sub-exponential $%
\alpha $- and $\beta $-mixing random variables %
\citep{Viennet1997,Bosq1993,Krebs2018b} with exponential memory decay %
\citep[e.g.][]{Merlevede_et_al_2011}, or geometric or even hyperbolic decay 
\citep[see][for bounded $\varphi $-mixing 1-Lipschitz
functions]{Wintenberger2010}. Results allowing for strong (or similar)
mixing have gone much farther to include spatial lattices %
\citep{Valenzuela-Dominguez_et_al_2017}, random fields \citep{Krebs2018a},
and less conventional mixing properties \citep{HangSteinwart2017}. Seminal
generic results are due to \cite{Talagrand1995a,Talagrand1995b}, leading to
inequalities for bounded stochastic objects 
\citep[see,
e.g.,][who work with bounded envelopes of $f$-mixing processes]{Samson2000}.

As a secondary contribution that will be of independent interest, we apply
the max-LLN's to three settings in order to yield new results. In each case
a bootstrap theory would complement the application but is ignored here for
brevity. We first consider a serial max-correlation statistic derived from a
model residual. \cite{HillMotegi2020} exploit Ramsey theory in order to
yield a complete bootstrap theory under a broad Near Epoch Dependence
property, yet without being able to characterize an upper bound on the
number of lags $k_{n}$. We provide new bounds on $k_{n}$ under $\alpha $%
-mixing and physical dependence.

The second application extends the marginal screening method to allow for an
increasing number of covariates under weak dependence. Marginal regressions
with \textquotedblleft optimal" covariate selection is also called sure
screening and correlation learning; see \cite{Genovese_etal_2012} for
references and historical details. In a recent contribution \cite%
{McKeague_Qian_2015} regress some $y_{t}$ on each covariate $\{x_{i,t}$ $:$ $%
1$ $\leq $ $i$ $\leq $ $k\}$ one at a time for \textit{fixed} $k$ that is
allowed to be larger than $n$ (note $t$ $=$ $1,...,n$). This yields marginal
coefficients $\hat{\theta}_{n,i}$ $=$ $\widehat{cov}(y,x_{i})/\widehat{var}%
(x_{i})$, max-index $\hat{l}_{n}$ $=$ $\arg \max_{1\leq l\leq k}|\hat{\theta}%
_{n,l}|$ ideally representing the most informative regressor, and therefore $%
\hat{\theta}_{n,\hat{l}_{n}}$. Let $\theta _{0,i}$ $=$ $%
cov(y,x_{i})/var(x_{i})$. An implicit iid assumption is imposed in order to
study $\hat{\theta}_{n,\hat{l}_{n}}$ as a vehicle for testing that no
regressor is correlated with $y_{t}$, $H_{0}$ : $\theta _{0,i^{\ast }}$ $=$ $%
0$ where $i^{\ast }$ $=$ $\arg \max_{1\leq l\leq k}|cov(y,x_{l})/var(x_{l})|$%
. See \cite{McKeague_Qian_2015} for discussion, and resulting non-standard
asymptotics for $\sqrt{n}(\hat{\theta}_{n,\hat{l}_{n}}$ $-$ $\theta
_{0,i^{\ast }})$.

We instead study $\max_{1\leq i\leq k_{n}}|\sqrt{n}\hat{\theta}_{n,i}|$ to
test $H_{0}$ : $\theta _{0,i}$ $=$ $0$ $\forall i$ $\Leftrightarrow $ $H_{0}$
: $\theta _{0,i^{\ast }}$ $=$ $0$, under weak dependence, allowing for
non-stationarity \textit{and} high dimensionality $k_{n}$ $>>$ $n$, where $%
k_{n}$ $\rightarrow $ $\infty $ and $k_{n}/n$ $\rightarrow $ $\infty $ are
allowed. We do not explore, nor do we need, an endogenously selected optimal
covariate index $\hat{l}_{n}$ under weak dependence. This narrowly relates
to work in \cite{Hill_2025_maxtest} where low dimensional models with a
fixed dimension nuisance covariate are used to test a HD parameter in an iid
regression setting.

The third application rests in the settings of \cite%
{CattaneoJanssonNewey2018} and \cite{Hill_2025_maxtest}. \cite%
{CattaneoJanssonNewey2018} study post-estimation inference when there are
many \textquotedblleft nuisance" parameters $\delta _{n,0}$ in a linear
regression model $y_{n,t}=\delta _{n,0}^{\prime }w_{n,t}+\theta
_{n,0}^{\prime }x_{n,t}+u_{n,t}$. Allowing for arbitrary in-group dependence
of finite group size, they deliver a heteroscedasticity-robust limit theory
for an estimator of the low dimensional $\theta _{n,0}$ by partialling out $%
\delta _{n,0}$. We extend their idea to weakly dependent and heterogeneous
data, but focus instead on testing the HD parameter $\delta _{n,0}$.

Finally, we focus on pointwise convergence throughout, ignoring uniform
convergence for high dimensional measurable mappings $x_{i,n,t}(\theta )$
with finite or infinite dimensional $\theta $. Generic results are well
known in low dimensional settings: see, e.g., \cite{Andrews1987} and \cite%
{Newey1991} for weak laws, \cite{PotscherPrucha1991} for a strong law, and 
\cite{vanderVaartWellner1996} for classic results for low dimensional $%
x_{i,n,t}(\cdot )$ with infinite dimensional $\theta $. Sufficient
conditions generally reduce to pointwise convergence, plus stochastic
equicontinuity (or related) conditions. The same generality likely extends
to a high dimensional setting, but this is left for future work.

The remainder of the paper is as follows. In Section \ref{sec:maxLLNs} we
present max-LLN's for mixing and physical dependent arrays. Sections \ref%
{sec:applic:mc}-\ref{sec:applic:param_test} contain applications, with
concluding remarks in Section \ref{sec:conclud}. Technical proofs of the
main results are presented in Appendix \ref{app:proofs}, and omitted content
is relegated to \cite{sm_max_LLN}.\medskip

We assume all random variables exist on the same complete measure space $%
(\Omega ,\mathcal{F},\mathbb{P})$ in order to side-step any measurability
issues concerning suprema \citep[e.g.][Appendix C]{Pollard1984}. $|x|$ $=$ $%
\sum_{i,j}|x_{i,j}|$ is the $l_{1}$-norm, $|x|_{2}$ $=$ $%
(\sum_{i,j}x_{i,j}^{2})^{1/2}$ is the Euclidean, Frobenius or $l_{2}$ norm; $%
||\cdot ||$ is the spectral norm; $||\cdot ||_{p}$ denotes the $\mathcal{L}%
_{p}$-norm ($||x||_{p}$ $:=$ $(\sum_{i=1}^{k}\mathbb{E}|x_{i}|^{p})^{1/p}$). 
$\mathbb{E}$ is the expectations operator; $\mathbb{E}_{\mathcal{A}}$ is
expectations conditional on $\mathcal{F}$-measurable $\mathcal{A}$. $\overset%
{p}{\rightarrow }$ and $\overset{\mathcal{L}_{p}}{\rightarrow }$\ denote
convergence in probability and in $\mathcal{L}_{p}$ norm. $a.s.$ means 
\textit{almost surely}. $o_{p}(1)$ depicts little \textquotedblleft $o$%
\textquotedblright\ convergence\textit{\ in probability}. \textit{awp1} =
\textquotedblleft asymptotically with probability approaching one". $\kappa $%
-\textit{Lipschitz} functions $f$ $:$ $\mathbb{R}^{r}$ $\rightarrow $ $%
\mathbb{R}$ satisfy $|f(x)$ $-$ $f(y)|$ $\leq $ $\kappa |x$ $-$ $y|$ for $%
\kappa $ $>$ $0$. $\{k_{n}\}_{n\in \mathbb{N}}$ is monotonically increasing. 
$\mathbb{N}$ $:=$ $\{1,2,...\}$ and $\mathbb{Z}$ $:=$ $\{...-1,0,1,2,...\}$. 
$K$ $>$ $0$ and tiny $\iota $ $>$ $0$ are constants that may change from
line to line. $O(m^{-\lambda })$ for integers $m$ $\geq $ $0$ and $\lambda $ 
$>$ $0$ implies $O((m$ $\vee $ $1)^{-\lambda })$. All stated distributions
are non-degenerate (with respect to Lebesgue measure).

\section{High dimensional maximal inequalities and LLN's\label{sec:maxLLNs}}

Let $\mathbb{E}x_{n,t}$ $=$ $0$ throughout. We first work with mixing arrays.

\subsection{Sub-Exponential $\protect\tau $-Mixing arrays\label%
{sec:bounded_mix}}

We discuss $\tau $-mixing in this section. Any other mixing condition cited
in the introduction with a sub-exponential condition will generally yield (%
\ref{maxLLN}) under an exponential bound for $k_{n}$. Define the filtration $%
\mathcal{F}_{i,n,s}^{t}$ $=$ $\sigma (x_{i,n,\tau }$ $:$ $-\infty $ $\leq $ $%
s$ $\leq $ $\tau $ $\leq $ $t$ $\leq $ $n)$.

The first result exploits a Fuk-Nagaev type Bernstein-inequality in %
\citet[Theorem 1]{Merlevede_et_al_2011} for \ geometric $\tau $-mixing
processes. $\tau $-mixing nests some non-$\alpha $-mixing processes %
\citep[cf.][]{Dedecker_Prieur_2004}, and is implied by $\alpha $-mixing and
nests $\mathcal{L}_{p}$-mixingales \citep{Hill2025_mixg,sm_max_LLN}. Let $%
\Lambda _{1}(\mathbb{R}^{r})$ denote the class of $1$-Lipschitz functions $f$
$:$ $\mathbb{R}^{r}$ $\rightarrow $ $\mathbb{R}$, let $\mathcal{A}$ be a $%
\sigma $-subfield of $\mathcal{F}$, and define for $\mathbb{R}^{r}$-valued
random variable $X$ as in \cite{Dedecker_Prieur_2004}: $\tau ^{(1)}(\mathcal{%
A},X)$ $:=$ $\left\Vert \sup_{f\in \Lambda _{1}(\mathbb{R}^{r})}\left\vert 
\mathbb{E}_{\mathcal{A}}f(X)-\mathbb{E}f(X)\right\vert \right\Vert _{1}$. If
we write for any $l$-tuple $J_{l}$ $:=$ $(j_{1},...,j_{l})$ $\in $ $\mathbb{N%
}^{l}$%
\begin{equation*}
\boldsymbol{X}_{i,n}(J_{l}):=\{x_{i,n,j_{1}},...,x_{i,n,j_{l}}\},
\end{equation*}%
then we have a generalization of the $\tau $-mixing coefficient defined only
for $p$ $=$ $1$ (see \citet[Defn. 2]{Dedecker_Prieur_2004} and 
\citet[eq.'s
(2.2)-(2.3)]{Merlevede_et_al_2011}) 
\begin{eqnarray*}
\tau _{i,n}^{(1)}(m) &=&\sup_{r\geq 1}\max_{1\leq l\leq r}\max_{1\leq t\leq
n}\max_{t\leq j_{1}<\cdots <j_{l}}\frac{1}{l}\tau ^{(1)}\left( \mathcal{F}%
_{i,n,-\infty }^{t-m},\boldsymbol{X}_{i,n}(J_{l})\right) \\
&=&\sup_{l\geq 1}\max_{1\leq t\leq n}\max_{t\leq j_{1}<\cdots <j_{l}}\frac{1%
}{l}\tau ^{(1)}\left( \mathcal{F}_{i,n,-\infty }^{t-m},\boldsymbol{X}%
_{i,n}(J_{l})\right) .
\end{eqnarray*}%
We use a trivial shift $\sup_{t\leq j_{1}<\cdots <j_{l}}\tau (\mathcal{F}%
_{i,n,-\infty }^{t-m},\mathcal{\cdot })$ rather than $\sup_{t+m\leq
j_{1}<\cdots <j_{l}}\tau (\mathcal{F}_{i,n,-\infty }^{t},\mathcal{\cdot })$
in \cite{Dedecker_Prieur_2004} in order to draw a direct comparison with
mixingales in \cite{sm_max_LLN}: the two versions are identical as $n$ $%
\rightarrow $ $\infty $, or under stationarity $\forall n$. Notice we do not
restrict dependence across coordinates $(x_{i,n,s},x_{j,n,t})$ for $i$ $\neq 
$ $j$. See \citet{Dedecker_Prieur_2004,DedeckerPrieur2005} and \cite%
{Dedecker_etal2007} for examples and further theory related to $\tau $%
-mixing and its relationship to other mixing properties.

Now use the $\mathcal{L}_{p}$ metric to define in general for $p$ $\geq $ $1$%
\begin{eqnarray*}
&&\tau ^{(p)}(\mathcal{A},X):=\mathbb{E}\left\vert \sup_{f\in \Lambda _{1}(%
\mathbb{R}^{r})}\left\vert \mathbb{E}_{\mathcal{A}}f(X)-\mathbb{E}%
f(X)\right\vert \right\vert ^{p} \\
&&\tau _{i,n}^{(p)}(m):=\sup_{l\geq 1}\max_{1\leq t\leq n}\max_{t\leq
j_{1}<\cdots <j_{l}}\frac{1}{l}\tau ^{(p)}\left( \mathcal{F}_{i,n,-\infty
}^{t-m},\boldsymbol{X}_{i,n}(J_{l})\right) .
\end{eqnarray*}%
Add and subtract $f(0)$ in $\mathbb{E}_{\mathcal{A}}f(X)$ $-$ $\mathbb{E}%
f(X) $, and use the $1$-Lipschitz property coupled with Minkowski and Jensen
inequalities to deduce an upper bound $\tau ^{(p)}(\mathcal{A},X)$ $\leq $ $2%
\mathbb{E}|X|^{p}$, hence $\lim \sup_{n\rightarrow \infty }\tau
_{i,n}^{(p)}(m)$ $\leq $ $2\lim \sup_{n\rightarrow \infty }\max_{1\leq t\leq
n}\mathbb{E}|x_{i,n,t}|^{p}$. We say $x_{i,n,t}$ is $\tau ^{(p)}$-mixing
when $\lim_{m\rightarrow \infty }\tau _{i,n}^{(p)}(m)$ $\rightarrow $ $0$.
Clearly $\tau _{i,n}^{(p)}(m)$ $\leq $ $\tau _{i,n}^{(q)}(m)^{p/q}$, $p$ $%
\leq $ $q$. The $\mathcal{L}_{p}$-variants ($\tau ^{(p)},\tau _{i,n}^{(p)}$)
share the same properties as ($\tau ^{(1)},\tau _{i,n}^{(1)}$), typically by
simple adjustments to existing proofs in \cite{Dedecker_Prieur_2004}, cf. 
\cite{Peligrad2002}. In particular, it retains a useful coupling property:
for any $\sigma $-field $\mathcal{A}$ of $\mathcal{F}$, there exists a
random variable $X^{\ast }$ distributed as $X$ and independent of $\mathcal{A%
}$ such that $\tau ^{(p)}(\mathcal{A},X)$ $=$ $\mathbb{E}|X$ $-$ $X^{\ast
}|^{p}$.

Assume geometric mixing decay and a sub-exponential tail condition in order
to focus ideas,%
\begin{eqnarray}
&&\max_{1\leq i\leq k_{n}}\tau _{i,n}^{(p)}(m)\leq ae^{-bm^{\gamma _{1}}}%
\text{ for some }p\geq 1\text{, }\forall n\in \mathbb{N}  \label{tao} \\
&&\max_{1\leq i\leq k_{n},1\leq t\leq n}\mathbb{P}\left( \left\vert
x_{i,n,t}\right\vert >\epsilon \right) \leq d\exp \left\{ -c\epsilon
^{\gamma _{2}}\right\} \text{ }\forall \epsilon >0\text{, }\forall n\in 
\mathbb{N}  \label{sube}
\end{eqnarray}%
where $(a,b,c,d,\gamma _{1},\gamma _{2})$ $>$ $0$ are constants. Define $%
\gamma $ by%
\begin{equation}
1/\gamma :=1/\gamma _{1}+1/\gamma _{2}\leq 1.  \label{ggg}
\end{equation}%
The latter imposes $\gamma _{i}$ $>$ $1$, forcing a trade-off between tail
decay and memory decay. When $\gamma _{2}$ $=$ $2$ we have the sub-Gaussian
class \citep[e.g.][Chapt. 2.5]{Vershynin2018}.

We now have the following Fuk-Nagaev type inequality under $\tau ^{(p)}$.

\begin{lemma}
\label{lm:tp_Bern}Under (\ref{tao})-(\ref{ggg}), for some constants $(%
\mathcal{K}_{j}\}_{j=1}^{5}$ $>$ $0$ depending only on $\{a,b,c,d,\gamma
,\gamma _{1},p\}$, any $n$ $\geq $ $4$, and any $\epsilon $ $>$ $1/4$ 
\begin{eqnarray}
&&\max_{1\leq i\leq k_{n}}\mathbb{P}\left( \max_{1\leq l\leq n}\left\vert 
\frac{1}{n}\sum_{t=1}^{l}x_{i,n,t}\right\vert \geq \epsilon \right)
\label{BE_t} \\
&&\text{ \ \ \ \ \ \ \ \ \ }\leq n\exp \left\{ -\mathcal{K}_{1}\epsilon
^{\gamma }n^{\gamma }\right\} +\exp \left\{ -\mathcal{K}_{2}\frac{\epsilon
^{2}n^{2}}{1+\mathcal{K}_{3}n}\right\} +\exp \left\{ -\mathcal{K}%
_{4}\epsilon ^{2}ne^{\frac{\mathcal{K}_{5}\left( \epsilon n\right) ^{\gamma
(1-\gamma )}}{[\ln (\epsilon n)]^{\gamma }}}\right\} .  \notag
\end{eqnarray}
\end{lemma}

The subsequent max-WLLN is proved using Lemma \ref{lm:tp_Bern} and a log-exp
bound.

\begin{theorem}[max-WLLN: $\protect\tau ^{(p)}$-mixing]
\label{thm:max_LLN_t_mix}Let $\{x_{i,n,t}$ $:$ $1$ $\leq $ $i$ $\leq $ $%
k_{n}\}_{t=1}^{n}$ satisfy (\ref{tao})-(\ref{ggg}). Then $\mathcal{M}_{n}$ $%
\overset{\mathcal{L}_{1}}{\rightarrow }$ $0$ provided $\ln (k_{n})$ $=$ $%
o(n) $. Moreover, $\sqrt{n}\mathcal{M}_{n}$ $=$ $O_{p}(\sqrt{\ln (k_{n})})$
if $\ln (k_{n})$ $=$ $O(n)$.
\end{theorem}

\begin{remark}
\normalfont The first result, $\mathcal{M}_{n}$ $\overset{\mathcal{L}_{1}}{%
\rightarrow }$ $0$ if $\ln (k_{n})$ $=$ $o(n)$, relies heavily on
sub-exponential tail decay and therefore Lemma \ref{lm:tp_Bern}. The second
result\ somewhat remarkably does not rely on the degree of dependence or
tail decay, $(\gamma _{1},\gamma _{2})$: the iid case $\sqrt{n}\mathcal{M}%
_{n}$ $=$ $O_{p}(\sqrt{\ln (k_{n})})$ is achieved generally. That said, it
arguably becomes less remarkable in view of the coupling between $\tau
^{(p)} $-mixing and independence, cf. \citet[Lemma
5]{Dedecker_Prieur_2004} and \citet[Lemma C.2]{sm_max_LLN}.
\end{remark}

\begin{example}[\textbf{Linear Processes}]
\normalfont Let $x_{i,t}$ $=$ $\sum_{j=0}^{\infty }\psi _{i,j}\epsilon
_{i,t-j}$, assume $\{\epsilon _{i,t}\}$ are iid for each $i$, $\mathbb{P}%
(|\epsilon _{i,t}|$ $>$ $u)$ $\leq $ $d\exp \{-cu^{\gamma _{2}}\}$ $\forall
i,t$ and constants $d,\gamma _{2}$ $>$ $0$, $\psi _{i,0}$ $=$ $1$ and $%
\sum_{j=0}^{\infty }|\psi _{i,j}|$ $<$ $\infty $ for each $i$. By exploiting
coupling results in \cite{MerlevedePeligrad2002}, and in view of $\tau
^{(p)} $-coupling \citep[Lemma C.2]{sm_max_LLN}, arguments in 
\citet[p.
214]{DedeckerPrieur2005} yield $\tau _{i,n}^{(p)}(m)$ $\leq $ $%
K(\sum_{j=m}^{\infty }|\psi _{i,j}|)^{p}$. Thus if $\max_{i\in \mathbb{N}%
}|\psi _{i,m}|$ $=$ $O(e^{-bm^{\gamma _{1}}})$ then by Theorem \ref%
{thm:max_LLN_t_mix} $\sqrt{n}\mathcal{M}_{n}$ $=$ $O_{p}(\sqrt{\ln (k_{n})})$%
\ whenever $\ln (k_{n})$ $=$ $o(n)$.
\end{example}

\begin{example}[$\boldsymbol{\protect\rho }$\textbf{-Lipschitz Markov Chains}%
]
\label{ex:p-Lip Mark}\normalfont Let $x_{i,t}$ $=$ $f_{i}(x_{i,t-1})$ $+$ $%
\epsilon _{i,t}$ with $\epsilon _{i,t}$ as above, where $f_{i}$ is $\rho
_{i} $--Lipschitz $\rho _{i}$ $\in $ $[0,1)$. If $x_{i,t}$ is $\mathcal{L}%
_{p}$-bounded\ then $\tau _{i,n}^{(p)}(m)$ $\leq $ $K\rho _{i}^{m}$ (see 
\citet[p.
215]{DedeckerPrieur2005}). Theorem \ref{thm:max_LLN_t_mix} applies when $%
\rho _{i}$ $\in $ $(0,e^{-b}]$, $b$ $>$ $0$.
\end{example}

\subsection{Physical dependence\label{sec:phys_dep}}

Next we augment the physical dependence measure in \citet[Defn.
1]{Wu2005} to cover non-stationary arrays, similar to 
\citet[Section
2.1.3]{ChangChenWu2024}. We initially ignore dependence across coordinates
as is standard in the literature, and then control cross-coordinate
dependence to improve the bound on $k_{n}$.

Suppose for measurable functions $g_{i,n,t}(\cdot )$ that may depend on $%
(i,n,t)$, 
\begin{equation}
x_{i,n,t}=g_{i,n,t}\left( \epsilon _{i,t},\epsilon _{i,t-1},\ldots \right)
\label{xe}
\end{equation}%
where $\{\epsilon _{i,t}\}$ are for each $i$ iid sequences. Examples include
linear and nonlinear time series like switching, random coefficient and
(non)linear Markov processes. Let $\{\epsilon _{i,t}^{\prime }\}$ be an
independent copy of $\{\epsilon _{i,t}\}$, and define the coupled process%
\begin{equation*}
x_{i,n,t}^{\prime }(m):=g_{i,n,t}\left( \epsilon _{i,t},\ldots ,\epsilon
_{i,t-m+1},\epsilon _{i,t-m}^{\prime },\epsilon _{i,t-m-1},\ldots \right) 
\text{, }m=0,1,2,...
\end{equation*}%
The (serial) $\mathcal{L}_{p}$-\textit{physical dependence} measure $\theta
_{i,n,t}^{(p)}(m)$ and its accumulation are defined as 
\begin{equation*}
\theta _{i,n,t}^{(p)}(m):=\left\Vert x_{i,n,t}-x_{i,n,t}^{\prime
}(m)\right\Vert _{p}\text{ and }\Theta _{i,n,t}^{(p)}:=\sum_{m=0}^{\infty
}\theta _{i,n,t}^{(p)}(m).
\end{equation*}%
We say $x_{i,n,t}$ is $\mathcal{L}_{p}$-physical dependent (over $t$, for
each $i$) when $\Theta _{i,n,t}^{(p)}$ $<$ $\infty $. This covers $\alpha $%
-mixing, $\tau ^{(p)}$-mixing, non-mixing and mixingale arrays %
\citep{Wu2005,Hill2025_mixg,sm_max_LLN}

\subsubsection{Unrestricted Coordinates\label{sec:unrestr_units}}

The following generalizes $\mathcal{L}_{p}$-moment (i.e. Rosenthal-type) and
Bernstein inequalities in \citet[Theorem 2]{Wu2005} to possibly
non-stationary arrays, complementing the HD central limit theory in 
\citet[Scetion
2.1.3]{ChangChenWu2024}. See \citet[Theorem 1]{LiuXiaoWu2013} for a modest
improvement in the stationary case. We need 
\begin{equation*}
\mathcal{Z}_{i,l}:=\frac{1}{\sqrt{n}}\sum_{t=1}^{l}x_{i,n,t}
\end{equation*}%
and%
\begin{equation}
\gamma _{i}(\alpha ):=\limsup_{p\rightarrow \infty }p^{1/2-1/\alpha }\Theta
_{i}^{(p)}\text{ with }\alpha >0\text{ and }\Theta
_{i}^{(p)}:=\limsup_{n\rightarrow \infty }\max_{1\leq t\leq n}\Theta
_{i,n,t}^{(p)}.  \label{gam_th}
\end{equation}

\begin{lemma}
\label{lm:phys_dep_Lp}Assume $\Theta _{i,n,t}^{(p)}$ $\in $ $(0,\infty )$
for $p$ $>$ $1$, and each $1$ $\leq $ $i$ $\leq $ $k_{n}$ and $1$ $\leq $ $t$
$\leq $ $n$. Write $p^{\prime }$ $:=$ $p$ $\wedge $ $2$.$\medskip $\newline
$a.$ $\left\Vert \max_{1\leq l\leq n}\left\vert \mathcal{Z}_{i,l}\right\vert
\right\Vert _{p}$ $\leq $ $\mathcal{B}_{p}n^{1/p^{\prime }-1/2}\max_{1\leq
t\leq n}\Theta _{i,n,t}^{(p)}$, where $\mathcal{B}_{p}$ $=$ $18[p^{5/2}/(p$ $%
-$ $1)^{3/2}]$ if $p$ $\in $ $(1,2)$, else $\mathcal{B}_{p}$ $=$ $\sqrt{2}%
[p^{3/2}/(p$ $-$ $1)]$.$\medskip $\newline
$b.$ If $\max_{i\in \mathbb{N}}\gamma _{i}(\alpha )$ $<$ $\infty $ for some $%
1$ $<$ $\alpha $ $\leq $ $2$, then $\max_{i\in \mathbb{N}}\mathbb{P}%
(\max_{1\leq l\leq n}|\mathcal{Z}_{i,l}|$ $>$ $u)$ $\leq $ $\mathcal{C}\exp
\{-\mathcal{K}u^{\alpha }\}$ for some $\mathcal{C},\mathcal{K}$ $\in $ $%
(0,\infty )$ that depend on $\gamma _{i}(\alpha )$\ (uniformly) and $\alpha $%
.
\end{lemma}

\begin{remark}
\normalfont($b$) exploits ($a$). ($a$) uses a martingale difference
approximation \citep[e.g.][]{Wu2005,Wu2011}, with Doob's inequality, and
either Burkholder's inequality when $p$ $\in $ $(1,2)$, or when $p$ $\geq $ $%
2$ a moment bound due to \cite{DedeckerDoukhan2003}, cf. 
\citet[Chapt.
2.5]{Rio2017}. See also \citet[Lemma 21]{JirakKostenberger2024}. We can
evidently also set $p$ $\in $ $(0,1]$ by using related general Doob-type
bounds \citep[e.g.][Theorem
4.4]{KuhnSchilling2023}.
\end{remark}

\begin{remark}
\label{rm:sube}\normalfont By Theorem 2 in \cite{Wu2005} the condition $%
\max_{i\in \mathbb{N}}\gamma _{i}(\alpha )$ $<$ $\infty $ yields the
existence of a pointwise moment generating function $\mathbb{E}\exp
\{\lambda \max_{1\leq l\leq n}|\mathcal{Z}_{i,l}|^{\alpha }\}$ $<$ $\infty $ 
$\forall \lambda $ $\in $ $[0,\lambda _{0})$ and some $\lambda _{0}$ $>$ $0$%
, and therefore sub-exponential tails $\mathbb{P}(\max_{1\leq l\leq n}|%
\mathcal{Z}_{i,l}|$ $>$ $u)$ $\leq $ $\mathcal{C}\exp \{-\mathcal{K}%
u^{\alpha }\}$ for some finite $\mathcal{C}$,$\mathcal{K}$ $>$ $0$ that may
take different values below.

It holds, for example, if $\theta _{i,n,t}^{(p)}(m)$ $\leq $ $%
d_{i,n,t}^{(p)}\psi _{i,m}$ where $\max_{i\in \mathbb{N}}\psi _{i,m}$ $=$ $%
O(m^{-\lambda -\iota })$, $\lambda $ $\geq $ $1$ and tiny $\iota $ $>$ $0$.
By construction $\theta _{i,n,t}^{(p)}(m)$ $\leq $ $2\left\Vert
x_{i,n,t}\right\Vert _{p}$ hence $d_{i,n,t}^{(p)}$ $\leq $ $2\left\Vert
x_{i,n,t}\right\Vert _{p}$. Then by a change of variable $p^{-b}\left\Vert
x_{i,n,t}\right\Vert _{p}$ $<$ $\infty $ uniformly in $(i,n,p,t)$ for some $%
b $ $\in $ $[0,\infty )$ yields $\max_{i\in \mathbb{N}}\gamma _{i}(\alpha )$ 
$< $ $\infty $ with $\alpha $ $=$ $2/(1$ $+$ $2b)$. When $b$ $\leq $ $1$ we
have classically defined sub-exponential tails 
\citep[cf.][Proposition
2.7.1(b)]{Vershynin2018}. The latter holds, for example, when $\mathbb{P}%
(|x_{i,n,t}|$ $>$ $u)$ $\leq $ $\mathcal{C}\exp \{-\mathcal{K}u^{\alpha }\}$
uniformly in $(i,t)$. Thus $\alpha $ $=$ $1$ (i.e. $b$ $=$ $1/2$) implies
sub-Gaussian tails.
\end{remark}

We now have a max-WLLN under physical dependence. The result allows for
trending dependence coefficients $\Theta _{n}^{(p)}$ $:=$ $\max_{1\leq i\leq
k_{n},1\leq t\leq n}\Theta _{i,n,t}^{(p)}$ $\rightarrow $ $\infty $ as $n$ $%
\rightarrow $ $\infty $. We work under $\mathcal{L}_{p}$-boundedness or
sub-exponential tails. In the former case, without cross-coordinate
dependence information we use Lyapunov's inequality and $\max_{1\leq i\leq
k}|x_{i}|$ $\leq $ $\sum_{i=1}^{k}|x_{i}|$ to obtain%
\begin{equation}
\mathbb{E}\mathcal{M}_{n}\leq \left\Vert \mathcal{M}_{n}\right\Vert _{p}\leq
\left( \sum_{i=1}^{k_{n}}\mathbb{E}\left\vert \frac{1}{n}%
\sum_{t=1}^{n}x_{i,n,t}\right\vert ^{p}\right) ^{1/p}\leq
k_{n}^{1/p}\max_{1\leq i\leq k_{n}}\left\Vert \frac{1}{n}%
\sum_{t=1}^{n}x_{i,n,t}\right\Vert _{p}.  \label{Mn2}
\end{equation}%
A max-WLLN thus rests on a maximal inequality to bound $\left\Vert
1/n\sum_{t=1}^{n}x_{i,n,t}\right\Vert _{p}$.

\begin{theorem}[max-WLLN: physical dependence]
\label{thm:max_LLN_phys_dep}Let $x_{i,n,t}$ be $\mathcal{L}_{p}$-physical
dependent, $p$ $>$ $1$, with $\Theta _{i,n,t}^{(p)}$ $\in $ $(0,\infty )$
for each $1$ $\leq $ $i$ $\leq $ $k_{n}$ and $1$ $\leq $ $t$ $\leq $ $n$.
Write $p^{\prime }$ $:=$ $p$ $\wedge $ $2$.$\medskip $\newline
$a.$ $\mathcal{M}_{n}$ $\overset{\mathcal{L}_{p}}{\rightarrow }$ $0$ if $%
k_{n}$ $=$ $o(n^{p(1-1/p^{\prime })}/\Theta _{n}^{(p)})$, and $\sqrt{n}%
\mathcal{M}_{n}$ $=$ $O_{p}(k_{n}^{1/p}n^{1/p^{\prime }-1/2}\Theta
_{n}^{(p)})$.$\medskip $\newline
$b.$ If additionally $\max_{i\in \mathbb{N}}\gamma _{i}(\alpha )$ $<$ $%
\infty $ for some $1$ $<$ $\alpha $ $\leq $ $2$, then $\mathcal{M}_{n}$ $%
\overset{\mathcal{L}_{1}}{\rightarrow }$ $0$ for any $\ln (k_{n})$ $=$ $o(%
\sqrt{n})$, and $\sqrt{n}\mathcal{M}_{n}$ $=$ $O_{p}(\ln (k_{n}))$.
\end{theorem}

\begin{remark}
\normalfont Under ($b$) we need $1$ $<$ $\alpha $ $\leq $ $2$ in order to
bound the moment generating function $\mathbb{E}\exp \{\lambda |1/\sqrt{n}%
\sum_{t=1}^{n}x_{i,n,t}|\}$ following a \emph{log-exp} bound and Bernstein
inequality Lemma \ref{lm:phys_dep_Lp}.b. This is not inconsequential
considering $\Theta _{i}^{(p)}$ from (\ref{gam_th}) and in $\gamma
_{i}(\alpha )$ is non-decreasing in $p$ by Lyapunov's inequality. This rules
out $\lim \sup_{p\rightarrow \infty }\Theta _{i}^{(p)}/p^{1/\alpha -1/2}$ $<$
$\infty $ for small $\alpha $ $<$ $1$, thus excluding \textquotedblleft
heavier tailed\textquotedblright\ cases where $\Theta _{i}^{(p)}$ $%
\rightarrow $ $\infty $ rapidly in $p$.
\end{remark}

\begin{remark}
\label{rm:Ademeck}\normalfont\citet[Lemma A.4]{Adamek_et_al2023} derive a
concentration inequality under an NED property with uniform $\mathcal{L}_{p}$%
-boundedness and $p$ $>$ $2$. Their result yields $\mathcal{M}_{n}$ $\overset%
{p}{\rightarrow }$ $0$ if $k_{n}$ $=$ $o(n^{p/2})$. We allow for trending
higher moments and $p$ $\in $ $(1,2]$, where physical dependence is implied
by the adapted mixingale property \citep[Theorem 2.1]{Hill2025_mixg}, and
mixingales nest NED \citep[Chap. 17]{Davidson1994}. Thus our max-WLLN is
broader in scope.
\end{remark}

\begin{remark}
\normalfont Using our notation and expanding terms, 
\citet[Theorem
3.2]{MiesSteland2023} prove under $\mathcal{L}_{q}$-bounded $\mathcal{L}_{p}$%
-physical dependence with coefficients $\theta _{i,n,t}^{(p)}(m)\leq
d_{i,n,t}^{(p)}$ $\times $ $(m\vee 1)^{-\beta -\iota }$, $\beta $ $\geq $ $1$%
, where $d_{i,n,t}^{(p)}$ $\leq $ $2\left\Vert x_{i,n,t}\right\Vert _{p}$,
and $2$ $\leq $ $p$ $\leq $ $q$, 
\begin{equation*}
\left\{ \mathbb{E}\max_{1\leq l\leq n}\left( \sum_{i=1}^{k_{n}}\left\vert 
\frac{1}{n}\sum_{t=1}^{l}x_{i,n,t}\right\vert ^{p}\right) ^{q/p}\right\}
^{1/q}\leq K\frac{1}{n^{1/2}}\mathcal{D}_{n}\sum_{m=1}^{\infty }\frac{1}{%
m^{\beta +\iota }},
\end{equation*}%
and $\mathcal{D}_{n}$ $:=$ $2[\max_{1\leq t\leq n}\mathbb{E}%
(\sum_{i=1}^{k_{n}}|x_{i,n,t}|^{p})^{q/p}]^{1/q}$. Cf. 
\citet[Theorem
6.6]{MiesSteland2023} and \citet[Theorem 4.1]{Pinelis1994}. Thus, they
deliver an $\mathcal{L}_{q}$-maximal inequality for the $l_{p}$-norm $%
(\sum_{i=1}^{k_{n}}|\sum_{t=1}^{l}x_{i,n,t}|^{p})^{1/p}$. The bound depends
implicitly on $k_{n}$ through $\mathcal{D}_{n}$. Set $q$ $=$ $p$ to be able
to yield%
\begin{equation*}
\left\{ \mathbb{E}\max_{1\leq l\leq n}\sum_{i=1}^{k_{n}}\left\vert \frac{1}{n%
}\sum_{t=1}^{l}x_{i,n,t}\right\vert ^{p}\right\} ^{1/p}\leq K\left( \frac{%
k_{n}}{n^{p/2}}\max_{1\leq i\leq k_{n},1\leq t\leq n}\left\Vert
x_{i,n,t}\right\Vert _{p}^{p}\right) ^{1/p}\text{ for }p\geq 2.
\end{equation*}%
Compare that with the implication of Lemma \ref{lm:phys_dep_Lp}.a and (\ref%
{Mn2}), 
\begin{eqnarray*}
\left\{ \mathbb{E}\max_{1\leq i\leq k_{n},1\leq l\leq n}\left\vert \frac{1}{n%
}\sum_{t=1}^{l}x_{i,n,t}\right\vert ^{p}\right\} ^{1/p} &\leq &\left\{ 
\mathbb{E}\max_{1\leq l\leq n}\sum_{i=1}^{k_{n}}\left\vert \frac{1}{n}%
\sum_{t=1}^{l}x_{i,n,t}\right\vert ^{p}\right\} ^{1/p} \\
&\leq &\mathcal{B}_{p}\left( \frac{k_{n}}{n^{p(1-1/p^{\prime })}}\max_{1\leq
i\leq k_{n},1\leq t\leq n}\left\{ \Theta _{i,n,t}^{(p)}\right\} ^{p}\right)
^{1/p}\text{ for }p>1.
\end{eqnarray*}%
If $p$ $\geq $ $2$ then $n^{p(1-1/p^{\prime })}$ $=$ $n^{p/2}$ and the upper
bounds are virtually identical since cosmetically $\Theta _{i,n,t}^{(p)}$ $%
\leq $ $2\left\Vert x_{i,n,t}\right\Vert _{p}^{p}$. The major differences
are \cite{MiesSteland2023} ($i$) operate on the \emph{larger} $\max_{1\leq
l\leq n}(\sum_{i=1}^{k_{n}}|1/n\sum_{t=1}^{l}x_{i,n,t}|^{p})^{q/p}$ with $q$ 
$\geq $ $p$; ($ii$) require $p$ $\geq $ $2$; ($iii$) only study $\mathcal{L}%
_{p}$-bounds; ($iv$) use telescoping sums of approximating martingales based
on arguments in \cite{Pinelis1994}. We also use martingale approximation
theory, based on classic arguments, to prove Lemma \ref{lm:phys_dep_Lp}.a
and therefore Theorem \ref{thm:max_LLN_phys_dep}.a.
\end{remark}

\begin{example}[\textbf{Iterated Random Functions}]
\normalfont Consider $x_{i,n,t}$ $=$ $F_{i,t}(x_{i,n,t-1},\epsilon _{i,t})$ $%
:=$ $F_{i,t}^{(\epsilon )}(x_{i,n,t-1})$ where $F_{i,t}$ are measurable
bivariate functions, and $\epsilon _{i,t}$ are iid. Assume as in %
\citet[Example 1]{LiuXiaoWu2013} the following fixed point and Lipschitz
properties: there exist points $z_{i,0}$, $\max_{i\in \mathbb{N}}|z_{i,0}|$ $%
<$ $\infty $, and $p$ $>$ $2$ such that 
\begin{eqnarray*}
&&\kappa _{i}(p):=\limsup_{n\rightarrow \infty }\left\{ \max_{1\leq t\leq
n}\left\Vert z_{i,0}-F_{i,t}^{(\epsilon )}(z_{i,0})\right\Vert _{p}\right\}
<\infty \\
&&\lambda _{i}(p):=\limsup_{n\rightarrow \infty }\max_{1\leq t\leq n}\left\{
\sup_{x\neq x^{\prime }}\frac{\left\Vert F_{i,t}^{(\epsilon
)}(x)-F_{i,t}^{(\epsilon )}(x^{\prime })\right\Vert _{p}}{\left\vert
x-x^{\prime }\right\vert }\right\} <1\text{ uniformly in }i.
\end{eqnarray*}%
Replicating arguments in \citet[Theorem 2]{WuShao2004} and 
\citet[Example
1]{LiuXiaoWu2013} yields a uniform functional dependence bound 
\begin{equation*}
\theta _{i}^{(p)}(m):=\limsup_{n\rightarrow \infty }\max_{1\leq t\leq
n}\left\Vert x_{i,n,t}-x_{i,n,t}^{\prime }(m)\right\Vert _{p}\leq \mathcal{K}%
_{p}\lambda _{i}^{m}(p)
\end{equation*}%
for some finite universal constant $\mathcal{K}_{p}$ $>$ $0$ depending only
on $\max_{i\in \mathbb{N}}|z_{i,0}|$, $\sup_{i\in \mathbb{N}}\kappa _{i}(p)$%
, and $\max_{i\in \mathbb{N}}\lambda _{i}(p)$. Therefore $\Theta _{i}^{(p)}$ 
$\leq $ $\mathcal{K}_{p}/(1-\lambda _{i}(p))$. Thus by Theorem \ref%
{thm:max_LLN_phys_dep}.a $\mathcal{M}_{n}$ $\overset{p}{\rightarrow }$ $0$
if $k_{n}$ $=$ $o(n^{p(1-1/p^{\prime })}/\Theta _{n}^{(p)})$, and $\sqrt{n}%
\mathcal{M}_{n}$ $=$ $O_{p}(k_{n}^{1/p}n^{1/p^{\prime }-1/2}\{\Theta
_{n}^{(p)}\}^{p})$.
\end{example}

We continue to work under \textit{serial} $\mathcal{L}_{p}$-physical
dependence, but now impose restrictions across coordinates $i$ to improve
bounds on $k_{n}$.

\subsubsection{Martingale Coordinates\label{sec:mart_units}}

Write 
\begin{equation*}
\mathcal{S}_{i,n}:=\frac{1}{n}\sum_{t=1}^{n}x_{i,n,t},
\end{equation*}%
and let the filtrations $\{\mathfrak{F}_{i,n}\}_{i\in \mathbb{N}}$ be such
that $\sigma (\{x_{i,n,t}\}_{t=1}^{n})$ $\subseteq $ $\mathfrak{F}_{i,n}$
and $\mathbb{E}_{\mathfrak{F}_{i,n}}x_{i+1,n,t}$ $=$ $x_{i,n,t}$ $\forall
(i,n,t).$ Then 
\begin{equation*}
\mathbb{E}_{\mathfrak{F}_{i,n}}\mathcal{S}_{i+1,n}=\mathcal{S}_{i,n},
\end{equation*}%
hence the collection $\{\mathcal{S}_{i,n},\mathfrak{F}_{i,n}\}_{i\geq 1}$
forms a martingale. Doob's inequality applies for any $p>1$ 
\citep[Theorem
2.2]{HallHeyde1980}, 
\begin{equation*}
\mathbb{E}\max_{1\leq i\leq k_{n}}\left\vert \mathcal{S}_{i,n}\right\vert
^{p}\leq \frac{p}{p-1}\mathbb{E}\left\vert \mathcal{S}_{k_{n},n}\right\vert
^{p}.
\end{equation*}%
Now apply Lemma \ref{lm:phys_dep_Lp}.a under physical dependence to $\mathbb{%
E}|\mathcal{S}_{k_{n},n}|^{p}$ to deduce%
\begin{equation*}
\mathbb{E}\max_{1\leq i\leq k_{n}}\left\vert \mathcal{S}_{i,n}\right\vert
^{p}\leq \mathcal{B}_{p}^{p}\frac{p}{p-1}n^{p/p^{\prime }-p}\left\{
\max_{1\leq t\leq n}\Theta _{k_{n},n,t}^{(p)}\right\} ^{p}=\mathcal{C}%
_{p}n^{p/p^{\prime }-p}\left\{ \max_{1\leq t\leq n}\Theta
_{k_{n},n,t}^{(p)}\right\} ^{p},
\end{equation*}%
with $\mathcal{C}_{p}$ implicit. Since $p/p^{\prime }$ $<$ $p$ we have $%
\mathbb{E}\max_{1\leq i\leq k_{n}}|\mathcal{S}_{i,n}|^{p}$ $\rightarrow $ $0$
\textit{for any} $\{k_{n}\}_{n\in \mathbb{N}}$ when $\max_{1\leq t\leq
n}\Theta _{k_{n},n,t}^{(p)}$ $=$ $o(n^{1-1/p^{\prime }})$. The latter holds,
for example, when $\theta _{i,n,t}^{(p)}(m)$ $\leq $ $d_{i,n,t}^{(p)}\psi
_{i,m}$ with $\psi _{i,m}$ $=$ $O(m^{-1-\iota })$ and bounded $\mathcal{L}%
_{p}$-trend $d_{i,n,t}^{(p)}$ $\leq $ $2\left\Vert x_{i,n,t}\right\Vert _{p}$
$\leq $ $Kt^{1-1/p^{\prime }-\iota }$ for tiny $\iota $ $>$ $0$. The latter
holds under uniform $\mathcal{L}_{p}$-boundedness $\lim \sup_{n\rightarrow
\infty }\max_{1\leq i\leq k_{n},1\leq t\leq n}\left\Vert
x_{i,n,t}\right\Vert _{p}$ $\leq $ $K$. A trivial special case is perfect
dependence $\mathbb{P}(x_{i,n,t}$ $=$ $x_{j,n,t})$ $=$ $1$ $\forall (i,j)$,
and a similar result extends to submartingales 
\citep[e.g.][Theorem
2.1]{HallHeyde1980}.

The preceding discussion with Markov's inequality proves the next result.

\begin{theorem}[max-WLLN: physical dependence over $t$, martingale over $i$]

\label{thm:max_WLLN_mart}Let $x_{i,n,t}$ be $\mathcal{L}_{p}$-physical
dependent over $t$, $p$ $>$ $1$, with $\Theta _{i,n,t}^{(p)}$ $\in $ $%
(0,\infty )$ for each $1$ $\leq $ $i$ $\leq $ $k_{n}$ and $1$ $\leq $ $t$ $%
\leq $ $n$. If there exist filtrations $\{\mathfrak{F}_{i,n}\}_{i\in \mathbb{%
N}}$ satisfying $\sigma (\{x_{i,n,t}\}_{t=1}^{n})$ $\subseteq $ $\mathfrak{F}%
_{i,n}$ and $\mathbb{E}_{\mathfrak{F}_{i,n}}x_{i+1,n,t}$ $=$ $x_{i,n,t}$ $%
\forall i,n,t$, then $\mathcal{M}_{n}$ $\overset{\mathcal{L}_{p}}{%
\rightarrow }$ $0$ for\emph{\ any} $\{k_{n}\}$ provided $\max_{1\leq t\leq
n}\Theta _{k_{n},n,t}^{(p)}$ $=$ $o(n^{1-1/p^{\prime }})$.
\end{theorem}

\subsubsection{Nearly Martingale Gaussian Coordinates\label%
{sec:nearmart_units}}

We relax the martingale assumption to hold only as $n$ $\rightarrow $ $%
\infty $ at a sufficiently slow rate. In the following we use maximum domain
of attraction theory for a Gaussian array to explore how classic theory
yields a better bound on $k_{n}$. Write 
\begin{equation*}
\mathcal{\tilde{Z}}_{i,n}:=\frac{1}{\mathcal{V}_{i,n}}\frac{1}{\sqrt{n}}%
\sum_{t=1}^{n}x_{i,n,t}=\frac{\sqrt{n}}{\mathcal{V}_{i,n}}\mathcal{S}_{i,n}%
\text{ where }\mathcal{V}_{i,n}^{2}:=\mathbb{E}\left( \frac{1}{\sqrt{n}}%
\sum_{t=1}^{n}x_{i,n,t}\right) ^{2},
\end{equation*}%
and assume $\lim \inf_{n\rightarrow \infty }\inf_{1\leq i\leq k_{n}}\mathcal{%
V}_{i,n}^{2}$ $>$ $0$. Notice $\max_{1\leq i\leq k_{n}}\mathcal{V}_{i,n}^{2}$
$=$ $O(1)$ by Lemma \ref{lm:phys_dep_Lp}.a when $\max_{1\leq i\leq
k_{n},1\leq t\leq n}\Theta _{i,n,t}^{(q)}$ $=$ $O(1)$. Assume $x_{i,n,t}$ $%
\sim $ $N(0,1)$ is strictly stationary over ($i,t$) to ease discussion, thus 
$\{\mathcal{\tilde{Z}}_{i,n}$ $:$ $1$ $\leq $ $i$ $\leq $ $k_{n}\}_{n\in 
\mathbb{N}}$ is a stationary standard normal array.

Define cross-coordinate correlations $\rho _{n,j}$ $:=$ $\mathbb{E}\mathcal{%
\tilde{Z}}_{i,n}\mathcal{\tilde{Z}}_{i+j,n}$. Then for some $\vartheta $ $%
\in $ $[0,1]$, and sequences $a_{n}$ $=$ $[2\ln (n)]^{1/2}$ and $b_{n}$ $%
\sim $ $[2\ln (n)]^{1/2}$, under regularity conditions on $\rho _{n,j}$ that
include $(1$ $-$ $\rho _{n,j})\ln (k_{n})$ $\rightarrow $ $\delta _{j}$ $\in 
$ $(0,\infty ]$, we have (see \citet[eq.'s (2.1)-(2.3)]{HsingHuslerReiss1996}%
, cf. \cite{Berman1964}) 
\begin{equation}
\mathbb{P}\left( \max_{1\leq i\leq k_{n}}\left\vert \mathcal{\tilde{Z}}%
_{i,n}\right\vert \leq u/a_{k_{n}}+b_{k_{n}}\right) \rightarrow \exp
\{-\vartheta \exp (-u)\}\text{ }\forall u\in (-\infty ,\infty ).
\label{MDA0}
\end{equation}%
The latter permits strong dependence with $|\rho _{n,j}|$ $<$ $1$ and $|\rho
_{n,j}|$ $\rightarrow $ $1$ as $n$ $\rightarrow $ $\infty $ at a
sufficiently slow rate. For example $\rho _{n,j}$ $=$ $(1-\zeta /\ln
(k_{n}))^{j}$ hence $\delta _{j}$ $=$ $j\zeta $ (see Example \ref{ex:AR}
below). We assume convergence in (\ref{MDA0}) holds uniformly:%
\begin{equation}
\sup_{u\in \mathbb{R}}\left\vert \mathbb{P}\left( \max_{1\leq i\leq
k_{n}}\left\vert \mathcal{\tilde{Z}}_{i,n}\right\vert \leq
u/a_{k_{n}}+b_{k_{n}}\right) -\exp \{-\vartheta \exp (-u)\}\right\vert
\rightarrow 0\text{ }\forall u\in (-\infty ,\infty ).  \label{MDA}
\end{equation}%
See, e.g., \citet[Theorem 3]{Cohen1982} and 
\citet[Sect.
5]{Herrmann_etal_2014}, cf. \citet[main Theorem]{Hall1979} and 
\citet[Lemma
3.1]{Berman1964}.

Now use uniform convergence (\ref{MDA}) with $b_{k_{n}}/\sqrt{n}$ $%
\rightarrow $ $0$, and therefore $\ln (k_{n})$ $=$ $o(n)$, and $\mathcal{V}%
_{i,n}^{2}$ $\in $ $(0,\infty )$ uniformly in $i,n$ to yield for each $u$ $>$
$0$ 
\begin{eqnarray*}
\mathbb{P}\left( \mathcal{M}_{n}>u\right) &=&\mathbb{P}\left( \max_{1\leq
i\leq k_{n}}\left\vert \frac{1}{n}\sum_{t=1}^{n}x_{i,n,t}\right\vert
>u\right) \\
&\leq &\mathbb{P}\left( \max_{1\leq i\leq k_{n}}\left\vert \mathcal{\tilde{Z}%
}_{i,n}\right\vert >\frac{\sqrt{n}u}{\max_{1\leq i\leq k_{n}}\mathcal{V}%
_{i,n}}\right) \\
&=&\mathbb{P}\left( a_{k_{n}}\left( \max_{1\leq i\leq k_{n}}\left\vert 
\mathcal{\tilde{Z}}_{i,n}\right\vert -b_{k_{n}}\right) >a_{k_{n}}\left( 
\frac{\sqrt{n}u}{\max_{1\leq i\leq k_{n}}\mathcal{V}_{i,n}}-b_{k_{n}}\right)
\right) \\
&=&\mathbb{P}\left( a_{k_{n}}\left( \max_{1\leq i\leq k_{n}}\left\vert 
\mathcal{\tilde{Z}}_{i,n}\right\vert -b_{k_{n}}\right) >\sqrt{n}%
a_{k_{n}}\left( \frac{u}{\max_{1\leq i\leq k_{n}}\mathcal{V}_{i,n}}-\frac{%
b_{k_{n}}}{\sqrt{n}}\right) \right) \rightarrow 0.
\end{eqnarray*}%
Compare this to $\ln (k_{n})$ $=$ $O(\sqrt{n})$ for sub-Gaussian arrays in
Theorem \ref{thm:max_LLN_phys_dep}.b. This proves the next max-WLLN result.

\begin{theorem}[max-WLLN: physical dependence over $t$, nearly martingale
over $i$]
\label{thm:max_WLLN_nearmart}Let $x_{i,n,t}$ be stationary $\mathcal{L}_{p}$%
-physical dependent, $p$ $>$ $1$, over $t$. Assume $x_{i,n,t}$ $\sim $ $%
N(0,1)$, $\rho _{n,j}$ satisfies (2.1)-(2.3) in \cite{HsingHuslerReiss1996},
and $\lim \inf_{n\rightarrow \infty }\inf_{1\leq i\leq k_{n}}\mathcal{V}%
_{i,n}^{2}$ $>$ $0$. Then $\mathcal{M}_{n}$ $\overset{p}{\rightarrow }$ $0$
if $\ln (k_{n})$ $=$ $o(n)$.
\end{theorem}

\begin{example}[\textbf{Gaussian AR(1) Array}]
\label{ex:AR}\normalfont Assume $x_{i,n,t}$ $\sim $ $N(0,1)$ is stationary,
serially $\mathcal{L}_{p}$-physical dependent, $p$ $>$ $1$, with $\Theta
_{i,t}^{(q)}$ $\in $ $(0,\infty )$. Suppose coordinate-wise $x_{i+1,n,t}$ $=$
$d_{n}x_{i,n,t}$ $+$ $\sqrt{1-d_{n}^{2}}\varepsilon _{i+1,t}$ $\forall (i,t)$
with mutually and serially iid $\varepsilon _{i,t}$ $\sim $ $N(0,1)$, and $%
d_{n}$ $:=$ $1$ $-$ $\zeta /\ln (k_{n})$ for some $\zeta $ $>$ $0$. Hence $%
(1 $ $-$ $\rho _{n,j})\ln (k_{n})$ $\rightarrow $ $j\zeta $ 
\citep[see][Section
3]{HsingHuslerReiss1996}. Moreover, $\mathbb{E}x_{i+1,n,t}^{2}$ $=$ $1$
hence $x_{i,t}$ $\sim $ $N(0,1)$ and $\mathcal{V}_{i,n}^{2}$ $=$ $%
1/n\sum_{t=1}^{n}\sum_{j=0}^{\infty }d_{n}^{2j}(1$ $-$ $d_{n}^{2})$ $=$ $1$ $%
\forall i$. Thus we have $\mathcal{\tilde{Z}}_{i+1,n}$ $=$ $d_{n}\mathcal{%
\tilde{Z}}_{i,n}$ $+$ $\mathcal{E}_{i+1,n}$ for iid $\mathcal{E}_{i,n}$ $%
\sim $ $N(0,1)$. Then $\mathcal{M}_{n}$ $\overset{p}{\rightarrow }$ $0$ if $%
\ln (k_{n})$ $=$ $o(n)$.
\end{example}

\subsubsection{Mixing Coordinates\label{sec:mix_units}}

Finally, define serial and cross-coordinate $\alpha $-mixing coefficients
under stationarity (to ease notation below):%
\begin{eqnarray*}
&&\alpha _{n}(m):=\max_{1\leq i\leq k_{n}}\sup_{\mathcal{A}\in \mathcal{F}%
_{i,n,-\infty }^{t},\mathcal{B}\in \mathcal{F}_{i,n,t+m}^{\infty
}}\left\vert \mathbb{P}\left( \mathcal{A}\cap \mathcal{B}\right) -\mathbb{P}(%
\mathcal{A})\mathbb{P}(\mathcal{B})\right\vert \\
&&\tilde{\alpha}_{n}(m):=\max_{1\leq i\leq k_{n}}\sup_{\mathcal{A}\in 
\mathcal{G}_{n,-\infty }^{i},\mathcal{B}\in \mathcal{G}_{n,i+m}^{\infty
}}\left\vert \mathbb{P}\left( \mathcal{A}\cap \mathcal{B}\right) -\mathbb{P}(%
\mathcal{A})\mathbb{P}(\mathcal{B})\right\vert \text{,}
\end{eqnarray*}%
with filtrations $\mathcal{F}_{i,n,s}^{t}$ $:=$ $\sigma (x_{i,n,\tau }$ $:$ $%
-\infty $ $\leq $ $s$ $\leq $ $\tau $ $\leq $ $t)$ and $\mathcal{G}%
_{n,i}^{j} $ $:=$ $\sigma (\{x_{l,n,t}\}_{t=1}^{n}$ $:$ $1$ $\leq $ $i$ $%
\leq $ $l$ $\leq $ $j)$.

Notice $\lim_{m\rightarrow \infty }\tilde{\alpha}_{n}(m)$ $=$ $0$ dictates
(asymptotic) cross-coordinate independence between samples $%
\{x_{i,n,t}\}_{t=1}^{n}$ and $\{x_{i+m,n,t}\}_{t=1}^{n}$ as $m$ $\rightarrow 
$ $\infty $. We need uniform serial mixing $\alpha _{n}(m)$ $\rightarrow $ $%
0 $ fast enough to ensure the serial physical dependence property holds, and
thus supports $\max_{1\leq i\leq k_{n}}\mathcal{V}_{i,n}^{2}$ $=$ $O(1)$. We
need $\tilde{\alpha}_{n}(m)$ $\rightarrow $ $0$ fast enough to ensure
cross-coordinate tail-based $\mathcal{D}$- and $\mathcal{D}^{\prime }$%
-mixing properties hold \citep{Leadbetter1974,Leadbetter1983}, promoting a
maximum domain of attraction condition akin to (\ref{MDA}). Recall $\mathcal{%
S}_{i,n}$ $:=$ $1/n\sum_{t=1}^{n}x_{i,n,t}$ and $\mathcal{Z}_{i,l}$ $:=$ $1/%
\sqrt{n}\sum_{t=1}^{l}x_{i,n,t}$.

\begin{theorem}[max-WLLN: \textit{mixing over (}$i,t$\textit{)}]
\label{thm:max_WLLN_mix_units}Let $\{x_{i,n,t}$ $:$ $1$ $\leq $ $t$ $\leq $ $%
n\}_{n\in \mathbb{N}}$ be stationary $\mathcal{L}_{p}$-bounded, $p$ $>$ $1$,
with $\lim \sup_{n\rightarrow \infty }\alpha _{n}(m)$ $=$ $O(m^{-\lambda
-\iota })$ for some $\lambda $ $>$ $qp/(q$ $-$ $p)$ and $q$ $>$ $p$, and $%
\lim \sup_{n\rightarrow \infty }\tilde{\alpha}_{n}(m)$ $=$ $O(m^{-2-\iota })$%
. Let for any $u$ $\in $ $\mathbb{R}$ and some sequences $\{a_{n},b_{n}\}$, $%
a_{n},b_{n}$ $>$ $0$, 
\begin{equation}
\max_{1\leq i\leq k_{n}}k_{n}\mathbb{P}\left( \left\vert \sqrt{n}\mathcal{S}%
_{i,n}\right\vert >u_{k_{n}}\right) \rightarrow \tau \in \lbrack 0,\infty )%
\text{ where }u_{n}:=u/a_{n}+b_{n}.  \label{knP}
\end{equation}%
Then $\mathcal{M}_{n}$ $\overset{p}{\rightarrow }$ $0$ sufficiently if 
\begin{equation}
\sqrt{n}a_{k_{n}}\rightarrow \infty \text{ and }b_{k_{n}}/\sqrt{n}%
\rightarrow 0.  \label{ab}
\end{equation}
\end{theorem}

It remains to ensure (\ref{knP}) for $a_{n},b_{n}$ such that (\ref{ab})
holds. We look at Gaussian, sub-exponential and stable domain cases, each
yielding a distinct bound on $k_{n}$ that is larger than the bound in
Theorem \ref{thm:max_LLN_phys_dep}. Moreover, in each case we can take $%
a_{k_{n}}$ $\sim $ $b_{k_{n}}$ $>$ $0$. Throughout serial and coordinate
mixing hold: $\lim \sup_{n\rightarrow \infty }\alpha _{n}(m)$ $=$ $%
O(m^{-\lambda -\iota })$ and $\lim \sup_{n\rightarrow \infty }\tilde{\alpha}%
_{n}(m)$ $=$ $O(m^{-2-\iota })$.

\begin{example}[\textbf{Gaussian}]
\label{ex:G}\normalfont Let $x_{i,n,t}$ $\sim $ $N(0,1)$ $\forall i,n,t$,
hence $\sqrt{n}\mathcal{S}_{i,n}$ $\sim $ $N(0,\mathcal{V}_{i,n}^{2})$ where 
$\mathcal{V}_{i,n}^{2}$ $=$ $n\mathbb{E}(\mathcal{S}_{i,n}^{2})$. By the
proof of Theorem \ref{thm:max_WLLN_mix_units} we may invoke Lemma \ref%
{lm:phys_dep_Lp}.a to deduce $\mathcal{V}_{i,n}^{2}$ $=$ $O(1)$. Assume
non-degeneracy $\lim \inf_{n\rightarrow \infty }\min_{1\leq i\leq k_{n}}%
\mathcal{V}_{i,n}^{2}$ $>$ $0$. We want $\{u_{n}\}_{n\in \mathbb{N}}$ such
that $\max_{1\leq i\leq k_{n}}\mathbb{P}(|\sqrt{n}\mathcal{S}_{i,n}|$ $>$ $%
u_{k_{n}})$ $\sim $ $\tau /k_{n}$. Well known Gaussian tail properties yield 
$u_{k_{n}}$ $\simeq $ $\sqrt{2\ln (k_{n})-2\ln (\tau \sqrt{2\pi })}$ $\simeq 
$ $\sqrt{2\ln (k_{n})}$. Thus we need $\ln (k_{n})$ $=$ $o(n)$ to yield (\ref%
{ab}). Compare that to $\ln (k_{n})$ $=$ $O(\sqrt{n})$ from Theorem \ref%
{thm:max_LLN_phys_dep}.a.
\end{example}

\begin{example}[\textbf{Sub-Exponential}]
\normalfont Let $\max_{1\leq i\leq k_{n},1\leq t\leq n}\mathbb{P}%
(|x_{i,n,t}| $ $>$ $u)$ $\leq $ $\mathcal{C}\exp \{-\mathcal{K}u^{\alpha }\}$
for some $1$ $<$ $\alpha $ $\leq $ $2$. Then by Lemma \ref{lm:phys_dep_Lp}.b
and Remark \ref{rm:sube}, $\max_{i\in \mathbb{N}}\mathbb{P}(\max_{1\leq
l\leq n}|\mathcal{Z}_{i,l}|$ $>$ $u)$ $\leq $ $\mathcal{C}\exp \{-\mathcal{K}%
u^{\alpha }\}$. Thus $\tau /k_{n}$ $\sim $ $\max_{i\in \mathbb{N}}\mathbb{P}%
(|\sqrt{n}\mathcal{S}_{i,n}|$ $>$ $u_{k_{n}})\leq \mathcal{C}\exp \{-%
\mathcal{K}u_{k_{n}}^{\alpha }\}$ if $u_{k_{n}}$ $\leq $ $(\mathcal{K}%
^{-1}\ln (k_{n}))^{1/\alpha }$. We therefore need $\ln (k_{n})$ $=$ $%
o(n^{\alpha /2})$. The sub-Gaussian case holds when $\alpha $ $=$ $2$,
reducing to Example \ref{ex:G}. This is a mild improvement over Theorem \ref%
{thm:max_LLN_phys_dep}.b where $\ln (k_{n})$ $=$ $O(n^{1/2})$ \emph{for any }%
$1$ $<$ $\alpha $ $\leq $ $2$.
\end{example}

\begin{example}[\textbf{Stable Domain}]
\normalfont Suppose $\max_{1\leq i\leq k_{n}}\mathbb{P}%
(\sum_{t=1}^{n}x_{i,n,t}/[n^{1/\varphi }h(n)]$ $>$ $u)$ $\rightarrow $ $%
\mathfrak{S}_{\varphi }(u)$ $\forall u$ $\in $ $\mathbb{R}$, some $\varphi $ 
$\in $ $(1,2)$, slowly varying $h(n)$ that may be different in different
places, and some zero mean distribution $\mathfrak{S}_{\varphi }(u)$. Then
under the stated mixing conditions $\mathfrak{S}_{\varphi }(u)$ is a stable
law with infinite variance (\citet[Theorem 1.1]{Ibragimov1962}, cf. 
\citet[Theorem
2.1]{Nagave1957}: symmetry and scale parameters are not shown here). Hence 
\begin{eqnarray*}
\max_{1\leq i\leq k_{n}}k_{n}\mathbb{P}\left( \left\vert \sqrt{n}\mathcal{S}%
_{i,n}\right\vert >u\right) &=&\max_{1\leq i\leq k_{n}}k_{n}\mathbb{P}\left(
\left\vert \frac{1}{n^{1/\varphi }h(n)}\sum_{t=1}^{n}x_{i,n,t}\right\vert >%
\frac{u_{k_{n}}}{n^{1/\varphi -1/2}h(n)}\right) \\
&\sim &k_{n}h(n)\left( \frac{u_{k_{n}}}{n^{1/\varphi -1/2}h(n)}\right)
^{-\varphi },
\end{eqnarray*}%
yielding $u_{k_{n}}$ $\sim $ $n^{1/\varphi -1/2}h(n)(k_{n}/\tau )^{1/\varphi
}$. We therefore need $k_{n}$ $=$ $o\left( n^{\varphi -1}/h(n)\right) $ for
some slowly varying $h(n)$. Compare this to $k_{n}$ $=$ $o(n^{q(1-1/q^{%
\prime })})$ under Theorem \ref{thm:max_LLN_phys_dep}.a with stationarity
for $q$ $>$ $1$. The index $\varphi $ is identically the moment supremum $%
\arg \sup \{r$ $:$ $\mathbb{E}|\mathcal{Z}_{i,n}|^{r}$ $<$ $\infty \}$, thus 
$q$ $<$ $\varphi $ $<$ $2$. This implies $q(1$ $-$ $1/q^{\prime })$ $=$ $q$ $%
-$ $1$ $<$ $\varphi $ $-$ $1$, and we again yield a modest improvement.
\end{example}

\section{Application \#1: max-correlation \label{sec:applic:mc}}

We now present three applications of the main results, pointing out max-LLN
usage by case. The first is a max-correlation test under mixing and physical
dependence settings. We do not develop a bootstrap theory in any application
to focus ideas.

\subsection{Residual max-correlation}

Consider a linear regression model $y_{t}$ $=$ $\phi _{0}^{\prime }x_{t-1}$ $%
+$ $\epsilon _{t}$, where $\phi _{0}$ $\in $ $\mathbb{R}^{k_{x}}$, $k_{x}$ $%
\geq $ $0$, $\mathbb{E}\epsilon _{t}$ $=$ $0$, with zero mean covariates $%
x_{t}$ $=$ $[x_{t,j}]$ $\in $ $\mathbb{R}^{k_{x}}$. We do not require $%
\mathbb{E}_{x_{t-1}}\epsilon _{t}$ $=$ $0$ $a.s.$ thus mis-specification is
allowed. We can easily allow for a non-linear model and conditional
volatility, but at the expense of additional notation and assumptions %
\citep[see][]{HillMotegi2020}. We want to test whether the model error is
white noise,%
\begin{equation*}
H_{0}:\mathbb{E}\epsilon _{t}\epsilon _{t-h}=0\text{ }\forall h\in \mathbb{N}%
\text{ against }H_{1}:\mathbb{E}\epsilon _{t}\epsilon _{t-h}\neq 0\text{ for
some }h\in \mathbb{N}.
\end{equation*}

Assume least squares estimation when $k_{x}$ $>$ $0$, $\hat{\phi}_{n}$ $=$ $%
(\sum_{t=2}^{n}x_{t-1}x_{t-1}^{\prime })^{-1}$ $\times $ $%
\sum_{t=2}^{n}x_{t-1}y_{t}$. Define the residual and its sample serial
covariance and correlation at lag $h$ $\geq $ $1$,%
\begin{equation*}
\epsilon _{t}(\hat{\phi}_{n}):=y_{t}-\hat{\phi}_{n}^{\prime }x_{t-1}\text{ \
and \ }\hat{\gamma}_{n}(h):=\frac{1}{n}\sum_{t=1+h}^{n}\epsilon _{t}(\hat{%
\phi}_{n})\epsilon _{t-h}(\hat{\phi}_{n})\text{ \ and }\hat{\rho}_{n}(h):=%
\frac{\hat{\gamma}_{n}(h)}{\hat{\gamma}_{n}(0)}.
\end{equation*}%
The test statistic is $\sqrt{n}\max_{1\leq h\leq k_{n}}|\hat{\rho}_{n}(h)|$
for some sequence of positive lags $\{k_{n}\}$. A weighted version of $\hat{%
\rho}_{n}(h)$ is possible allowing for standardization, or weighting to
account for lagging \citep[cf.][]{LjungBox1978}. Similarly, other estimators
can be entertained, e.g. GMM, LAD, QML, and so on, although $\sqrt{n}$%
-asymptotics is assumed.

\cite{HillMotegi2020} use Ramsey theory to sidestep conventional HD
approximations, ultimately using standard theory. They therefore cannot
bound $\{k_{n}\}$ although $k_{n}$ $=$ $o(n)$ must hold for consistency of $%
\hat{\rho}_{n}(h)$. We now use HD max-LLN's in part to prove \textit{any} $%
k_{n}$ $=$ $o(n)$ is valid.

\subsection{Strong mixing}

\subsubsection{Assumptions}

Let $\{\upsilon _{t}\}$ be an $\alpha $-mixing process with $\sigma $-fields 
$\mathfrak{V}_{s}^{t}$ $:=$ $\sigma (\upsilon _{\tau }$ $:$ $s$ $\leq $ $%
\tau $ $\leq $ $t)$ and $\mathfrak{V}_{t}:=\mathfrak{V}_{-\infty }^{t}$, and
coefficients $\alpha (m)$ $=$ $\sup_{t\in \mathbb{N}}\sup_{\mathcal{A}%
\subset \mathfrak{V}_{t}^{\infty },\mathcal{B}\subset \mathfrak{V}_{-\infty
}^{t-m}}|\mathbb{P}\left( \mathcal{A}\cap \mathcal{B}\right) $ $-$ $\mathbb{P%
}\left( \mathcal{A}\right) \mathbb{P}\left( \mathcal{B}\right) |$ $%
\rightarrow $ $0$ as $m$ $\rightarrow $ $\infty $. We impose second order
stationarity to reduce notation, but otherwise allow for global
nonstationarity. Define $\widehat{\mathcal{H}}_{n}$ $:=$ $%
1/n\sum_{t=2}^{n}x_{t-1}x_{t-1}^{\prime }$, $\mathcal{H}$ $:=$ $\mathbb{E}%
x_{t}x_{t}^{\prime }$, $\rho (h)$ $:=$ $\mathbb{E}\epsilon _{t}\epsilon
_{t-h}/\mathbb{E}\epsilon _{t}^{2}$, \ 
\begin{eqnarray*}
&&\mathcal{D}_{t}(h):=\mathbb{E}x_{t-1}\epsilon _{t}\epsilon _{t-h}+\mathbb{E%
}\epsilon _{t}x_{t-1-h}\epsilon _{t-h}\text{ and }\mathfrak{D}_{n}(h):=\frac{%
1}{n}\sum_{t=1+h}^{n}\{\mathcal{D}_{t}(h)+\mathcal{D}_{t}(-h)-2\rho (h)%
\mathcal{D}_{t}(0)\} \\
&&z_{n,t}(h):=\frac{\left\{ \epsilon _{t}\epsilon _{t-h}-\mathbb{E}\epsilon
_{t}\epsilon _{t-h}\right\} -\rho (h)\left\{ \epsilon _{t}^{2}-\mathbb{E}%
\epsilon _{t}^{2}\right\} -\mathcal{H}^{-1}x_{t-1}\epsilon _{t}\mathfrak{D}%
_{n}(h)}{\mathbb{E}\epsilon _{t}^{2}} \\
&&\mathcal{Z}_{n}(h):=\frac{1}{\sqrt{n}}\sum_{t=1+h}^{n}z_{n,t}(h).
\end{eqnarray*}%
Assume $\upsilon _{t}$ $=$ $g_{t}(\varepsilon _{t},\varepsilon _{t-1},...)$
for measurable $g_{t}(\cdot )$ and iid $\{\varepsilon _{t}\}$.

\begin{assumption}[data generating process]
\label{assum:maxcor}$\medskip $\newline
$a$. $(\epsilon _{t},x_{t})$ are zero-mean, $\mathfrak{V}_{t}$-measurable,
second order stationary, $\mathbb{E}\epsilon _{t}^{2}$ $>$ $0$, and $\mathbb{%
E}(\epsilon _{t}x_{t-1})$ $=$ $0$ for unique $\phi _{0}$, an interior point
of compact $\Phi $ $\subset $ $\mathbb{R}^{k_{x}}$. Each $w_{t}$ $\in $ $%
\{\epsilon _{t},x_{t}\}$ is governed by a distribution satisfying $%
\max_{t\in \mathbb{N}}\mathbb{P}(|w_{t}|$ $>$ $z)$ $\leq $ $b_{1}\exp
\{-b_{2}z^{\gamma _{1}}\}$ $\forall z$ $\geq $ $0$ for some universal
constants $(b_{1},b_{2},\gamma _{1})$ $>$ $0$.$\medskip $\newline
$b.$ $\alpha (m)$ $\leq $ $a_{1}\exp \{-a_{2}m^{\gamma _{2}}\}$ for
constants $(a_{1},a_{2},\gamma _{2})$ $>$ $0$.$\medskip $\newline
$c.$ $\mathcal{H}$ is positive definite, and $\widehat{\mathcal{H}}_{n}$ is $%
a.s.$ positive definite $\forall n$ $\geq $ $n_{0}$, and some $n_{0}$ $\in $ 
$\mathbb{N}$.$\medskip $\newline
$d.$ $\lim \inf_{n\rightarrow \infty }\inf_{\lambda ^{\prime }\lambda =1}%
\mathbb{E}[(\sum_{h=1}^{\mathcal{L}}\lambda _{h}\mathcal{Z}_{n}(h))^{2}]$ $>$
$0$ for each $\mathcal{L}$ $\in $ $\mathbb{N}$.
\end{assumption}

\begin{remark}
\normalfont($a$)-($c$) allow us to use Theorem \ref{thm:max_LLN_t_mix} for
key summands by exploiting the fact that geometric $\alpha $-mixing implies
geometric $\tau ^{(p)}$-mixing. ($c$) is standard for least squares
identification. (d) is a conventional non-degeneracy property, required here
for a HD central limit theorem. It holds when $\boldsymbol{Z}_{n}(\mathcal{L}%
)$ $:=$ $[\mathcal{Z}_{n}(1),...,\mathcal{Z}_{n}(\mathcal{L})]^{\prime }$
satisfy a standard positive definiteness property: $\inf_{\lambda ^{\prime
}\lambda =1}\mathbb{E}[(\lambda ^{\prime }\boldsymbol{Z}_{n}(\mathcal{L}%
))^{2}]$ $>$ $0$ $\forall \mathcal{L}$, $\forall n$ $\geq $ $n_{1}$ and some 
$n_{1}$ $\in $ $\mathbb{N}$.
\end{remark}

\subsubsection{Main results}

In our first result, we may use $\sqrt{n}\max_{1\leq h\leq k_{n}}|\mathcal{%
\hat{Z}}_{n}(h)|$ to yield a HD Gaussian approximation for $\sqrt{n}%
\max_{1\leq h\leq k_{n}}\{|\hat{\rho}_{n}(h)$ $-$ $\rho (h)|\}$. The proof
exploits $\alpha $- and $\tau ^{(p)}$-mixing linkage Lemma C.3 in \cite%
{sm_max_LLN}, and $\tau ^{(p)}$-mixing max-LLN Theorem \ref%
{thm:max_LLN_t_mix}. See \citet[Appendix
D.1]{sm_max_LLN} for proofs.

\begin{lemma}
\label{lm:mc_nong}Under Assumption \ref{assum:maxcor}, $|\sqrt{n}\max_{1\leq
h\leq k_{n}}\left\vert \hat{\rho}_{n}(h)-\rho (h)\right\vert $ $-$ $\sqrt{n}%
\max_{1\leq h\leq k_{n}}\left\vert \mathcal{Z}_{n}(h)\right\vert |$ $\overset%
{p}{\rightarrow }$ $0$ for any $\{k_{n}\}$, $k_{n}$ $=$ $o(n).$
\end{lemma}

By Lemma \ref{lm:mc_nong} if suffices to yield a Gaussian approximation for $%
\max_{1\leq h\leq k_{n}}|\mathcal{Z}_{n}(h)|$. Define $\sigma _{n}^{2}(h)$ $%
:=$ $\mathbb{E}[\mathcal{Z}_{n}^{2}(h)]$. Let $\{\boldsymbol{Z}_{n}(h)$ $:$ $%
h$ $\in $ $\mathbb{N}\}_{h\geq 1}$ be an array of normally distributed
random variables, $\boldsymbol{Z}_{n}(h)$ $\sim $ $N(0,\sigma _{n}^{2}(h))$.
Assumption \ref{assum:maxcor}.a,b,e ensure $0$ $<$ $\sigma _{n}^{2}(h)$ $<$ $%
\infty $ uniformly in $h$ $\geq $ $1$ and $n$ $\geq $ $n_{1}$ for some $%
n_{1} $ $\in $ $\mathbb{N}$. The lower bound is Assumption \ref{assum:maxcor}%
.e. Consider the upper bound. Since $(\epsilon _{t},x_{t})$ are $\sigma
(\varepsilon _{\tau }$ $:$ $\tau $ $\leq $ $t)$-measurable and geometric $%
\alpha $-mixing, by uniform sub-exponentiality Lemma D.1 in \cite{sm_max_LLN}
and Lemma 2.1 in \cite{McLeish1975}$\ z_{n,t}(h)$ is an adapted geometric $%
\mathcal{L}_{2}$-mixingale with uniformly bounded constants. Thus $%
z_{n,t}(h) $ is geometrically physical dependent by Corollary 2.2 in \cite%
{Hill2025_mixg}. Hence $\max_{1\leq h\leq k_{n}}\sigma _{n}^{2}(h)$ $=$ $%
O(1) $ by Lemma \ref{lm:phys_dep_Lp}.

Define the Kolmogorov distance 
\begin{equation}
\delta _{n}:=\sup_{z\geq 0}\left\vert \mathbb{P}\left( \max_{1\leq h\leq
k_{n}}\left\vert \mathcal{Z}_{n}(h)\right\vert \leq z\right) -\mathbb{P}%
\left( \max_{1\leq h\leq k_{n}}\left\vert \boldsymbol{Z}_{n}(h)\right\vert
\leq z\right) \right\vert .  \label{Gauss_apprx}
\end{equation}

\begin{lemma}
\label{lm:Gauss_approx}Under Assumption \ref{assum:maxcor}, $\delta _{n}$ $%
\rightarrow $ $0$ for any $\{k_{n}\}$ satisfying $k_{n}$ $=$ $o(n^{1/9}(\ln
(n))^{1/3})$.
\end{lemma}

\begin{remark}
\label{rm:pd_mix}\normalfont The proof requires that describe the \emph{joint%
} mixing property of high dimensional $[z_{n,t}(h)]_{h=1}^{k_{n}}$: the
latter is $\alpha $-mixing with coefficients at displacement $m$ bounded
from above by $\alpha (\{m$ $-$ $k_{n}\}$ $\vee $ $0)$. This slows down
dependence decay, strongly impacting the allowed rate of divergence $k_{n}$ $%
\rightarrow $ $\infty $ \citep[cf.][Proposition
3]{ChangJiangShao2023}. We will see below that use of physical dependence is
a boon when the dimension involves such lags, since it need only hold for $%
z_{n,t}(h)$ \emph{uniformly} over $h$ rather than \emph{jointly} for $%
[z_{n,t}(h)]_{h=1}^{k_{n}}$, and allows for slower-than-geometric memory
decay.
\end{remark}

Lemmas \ref{lm:mc_nong} and \ref{lm:Gauss_approx} yield the desired
sharpening of Hill and Motegi's (\citeyear{HillMotegi2020}) HD Gaussian
approximation theory.

\begin{theorem}
\label{thm:rho_p_gaus}Under Assumption \ref{assum:maxcor}, for any $%
\{k_{n}\} $, $k_{n}$ $=$ $o(n^{1/9}(\ln (n))^{1/3})$, 
\begin{equation*}
\sup_{z\geq 0}\left\vert \mathbb{P}\left( \sqrt{n}\max_{1\leq h\leq
k_{n}}\left\vert \hat{\rho}_{n}(h)-\rho (h)\right\vert \leq z\right) -%
\mathbb{P}\left( \max_{1\leq h\leq k_{n}}\left\vert \boldsymbol{Z}%
_{n}(h)\right\vert \leq z\right) \right\vert \rightarrow 0.
\end{equation*}
\end{theorem}

\subsection{Physical dependence}

Now let $\{u_{t},v_{t}\}_{t\in \mathbb{Z}}$ be iid sequences, and assume
there exist measurable $\mathbb{R}^{k_{x}}$-valued and $\mathbb{R}$-valued
functions $g_{t}(\cdot )$ and $h_{t}(\cdot )$ satisfying $x_{t}$ $=$ $%
g_{t}(u_{t},u_{t-1},...)$ and $\epsilon _{t}$ $=$ $h_{t}(v_{t},v_{t-1},...).$
Let $\{u_{t}^{\prime },v_{t}^{\prime }\}_{t\in \mathbb{Z}}$ be independent
copies of $\{u_{t},v_{t}\}_{t\in \mathbb{Z}}$, and $\epsilon _{t}^{\prime
}(m)$ and $x_{t}^{\prime }(m)$ be the coupled versions based on $%
\{u_{t}^{\prime },v_{t}^{\prime }\}_{t\in \mathbb{Z}}$. Define $\mathcal{L}%
_{p}$-dependence coefficients $\theta _{t}^{(p)}(m)$ $:=$ $\left\Vert
x_{t}-x_{t}^{\prime }(m)\right\Vert _{p}$ and $\tilde{\theta}_{t}^{(p)}(m)$ $%
:=$ $\left\Vert \epsilon _{t}-\epsilon _{t}^{\prime }(m)\right\Vert _{p}$%
.\medskip \newline
\textbf{Assumption \ref{assum:maxcor}}.\textit{$b^{\ast }$ (physical
dependence). $(x_{i,t},\epsilon _{t})$ are $\mathcal{L}_{p}$-physical
dependent: there exist constants $d_{t}^{(p)}(h)$ and coefficients $\psi
_{m} $ satisfying for some $p$ $\geq $ $4$, and some size $\lambda $ $\geq $ 
$2$, 
\begin{equation*}
\left\{ \theta _{t}^{(p)}(m)\vee \tilde{\theta}_{t}^{(p)}(m)\vee \tilde{%
\theta}_{t-h}^{(p)}(m)\right\} \leq d_{t}^{(p)}(h)\psi _{m}\text{ and }\psi
_{m}=O(m^{-\lambda -\iota })\text{ and }\psi _{0}=1.
\end{equation*}%
}

The following result produces a significant improvement on $k_{n}$ for
reasons given above (see Remark \ref{rm:pd_mix}). See 
\citet[Appendix
D.2]{sm_max_LLN} for a proof.

\begin{theorem}
\label{thm:rho_p_phys_gauss}Under Assumption \ref{assum:maxcor}.a,b$^{\ast }$%
,c,d, for any $\{k_{n}\}$, $k_{n}$ $=$ $o(n)$,%
\begin{equation*}
\sup_{z\geq 0}\left\vert \mathbb{P}\left( \sqrt{n}\max_{1\leq h\leq
k_{n}}\left\vert \hat{\rho}_{n}(h)-\rho (h)\right\vert \leq z\right) -%
\mathbb{P}\left( \max_{1\leq h\leq k_{n}}\left\vert \boldsymbol{Z}%
_{n}(h)\right\vert \leq z\right) \right\vert \rightarrow 0.
\end{equation*}
\end{theorem}

\section{Application \#2: marginal screening\label{sec:applic:mr}}

\subsection{Test statistic}

Consider a scalar outcome $y$ and set of covariates $x$ $=$ $%
[x_{i}]_{i=1}^{k_{n}}$ with variances $v(y),v(x_{i})$ $\in (0,\infty )$. We
want to test the hypothesis that no covariate is linearly related to $y$: 
\begin{eqnarray}
&&H_{0}:cov\left( y,x_{i}\right) =0\text{ }\forall i=1,...,k_{n}\text{ and
each }n\text{ }  \label{H:marg_reg} \\
&&H_{1}:cov\left( y,x_{i}\right) \neq 0\text{ for some }\forall i=1,...,k_{n}%
\text{ as }n\rightarrow \infty ,  \notag
\end{eqnarray}%
where the number of covariates $k_{n}$ $\rightarrow $ $\infty $, and $%
k_{n}/n $ $\rightarrow $ $\infty $ is allowed. It is a simple generalization
to permit $k_{n}$ $\rightarrow $ $k^{\ast }$ for some $k^{\ast }$ $\in $ $%
\mathbb{N}\cup \infty $.

Now consider a sample of a covariance stationary process $%
\{y_{t},x_{t}\}_{t=1}^{n}$, and marginal regression models%
\begin{equation*}
y_{t}=\delta _{i}+\phi _{i}x_{i,t}+v_{i,t}=\beta _{i}^{\prime }\tilde{x}%
_{i,t}+v_{i,t},
\end{equation*}%
where $\tilde{x}_{i,t}$ $=$ $[1,x_{i,t}]^{\prime }$, $\mathbb{E}v_{i,t}$ $=$ 
$0$ for each model $i$ $=$ $1,...,k_{n}$, and the errors and covariates are
orthogonal $\mathbb{E}x_{i,t}v_{i,t}$ $=$ $0$. The sub-script
\textquotedblleft $i$" shows $\beta _{i}$ may be different for different
regressors $x_{i,t}$. Classically of course the (pseudo) true values are%
\begin{equation*}
\phi _{i}=\frac{cov\left( y,x_{i}\right) }{v(x_{i})}\text{ and }\delta _{i}=%
\mathbb{E}y_{t}-\phi _{i}\mathbb{E}x_{i,t}.
\end{equation*}%
Notice $y_{t}$ $=$ $\delta _{i}$ $+$ $v_{i,t}$ under $H_{0}$ for all $i$,
thus $\delta _{i}$ $:=$ $\mathbb{E}y_{t}$ $\forall i$, and tautologically $%
v_{i,t}$ $=$ $v_{t}$ $:=$ $y_{t}$ $-$ $\mathbb{E}y_{t}$ $\forall i.$

Define least squares estimators $\hat{\beta}_{i}$ $=$ $[\hat{\delta}_{i},%
\hat{\phi}_{i}]^{\prime }$ $:=$ $(\sum_{t=1}^{n}\tilde{x}_{i,t}\tilde{x}%
_{i,t}^{\prime })^{-1}$ $\times $ $\sum_{t=1}^{n}\tilde{x}_{i,t}y_{t}$, and,
e.g., $\bar{x}_{i,n}$ $:=$ $1/n\sum_{t=1}^{n}x_{i,t}$, hence%
\begin{equation*}
\hat{\phi}_{i}=\frac{1/n\sum_{t=1}^{n}\left( x_{i,t}-\bar{x}_{i,n}\right)
\left( y_{t}-\bar{y}_{n}\right) }{1/n\sum_{t=1}^{n}\left( x_{i,t}-\bar{x}%
_{i,n}\right) ^{2}}.
\end{equation*}%
\cite{McKeague_Qian_2015} study $\hat{\phi}_{(\hat{\imath}_{n}),n}$ as a
mechanism for testing (\ref{H:marg_reg}) based on an adaptively selected $%
\hat{\imath}_{n}$ $:=$ $\arg \max_{1\leq i\leq k_{n}}|\hat{\phi}_{i}|$ with
iid $\{y_{t},x_{t}\}$. They allow $k$ $>$ $n$ but for \textit{fixed} $k$.
They present an \textit{adaptive resampling test} in order to resolve
non-uniform and therefore non-standard asymptotics implicit in $\sqrt{n}|%
\hat{\phi}_{(\hat{\imath}_{n}),n}|$. See their introduction for historical
references.

A test of $H_{0}$, however, can also be based simply on $\sqrt{n}\max_{1\leq
i\leq k_{n}}|\hat{\phi}_{i}|$ without an endogenously selected $\hat{\imath}%
_{n}$. This alleviates the need for adaptive re-sampling, while inference is
easily gained by multiplier (block) bootstrap in high dimension, cf. \cite%
{Hill_2025_maxtest} and \cite{HillLi2025}. Historically, of course, there is
interest in an endogenously selected \textquotedblleft most informative"
regressor, and generally post-model-selection inference. See, e.g., \cite%
{LeebPotscher2006}, and consult \cite{McKeague_Qian_2015} for further
reading.

The present theory is related to the HD parameter test in \cite%
{Hill_2025_maxtest}. In that setting an iid linear regression model is
explored, with fixed (low) dimension nuisance parameter and a HD parameter
to be tested. The present setting allows for weak dependence and
non-stationarity with two tail settings: sub-exponentiality and $\mathcal{L}%
_{p}$-boundedness, yielding respectively exponential and polynomial bounds
on $k_{n}$. We impose second order stationarity to focus on the marginal
regression setting itself.

\subsection{Assumptions and main results\label{sec:mr_assum}}

Define compact parameter spaces $\Phi _{i},\mathcal{D}$ $\subset $ $\mathbb{R%
}$, and assume $0$ and $\phi _{i}$ are interior points of $\Phi _{i}$.
Define $\mathcal{B}_{i}$ $:=$ $\Phi _{i}\times\mathcal{D}$. As long as there
is no confusion we say \textit{uniformly} to denote $\limsup_{n\rightarrow
\infty }\max_{1\leq i\leq k_{n},1\leq t\leq n}$ or $\lim \inf_{n\rightarrow
\infty }\min_{1\leq i\leq k_{n}}$, depending on the case.

Let $\{\epsilon _{t}\}_{t\in \mathbb{Z}}$ be an iid sequence, and assume
there exist measurable $\mathbb{R}^{k_{n}+1}$-valued functions $g_{t}(\cdot
) $ satisfying $[x_{1,t},...,x_{k_{n},t},y_{t}]^{\prime }$ $=$ $%
g_{t}(\epsilon _{t},\epsilon _{t-1},\ldots ).$ Define $\mathcal{L}_{p}$%
-dependence coefficients $\theta _{i,t}^{(p)}(m)$ $:=$ $\left\Vert
x_{i,t}-x_{i,t}^{\prime }(m)\right\Vert _{p}$ and $\tilde{\theta}%
_{t}^{(p)}(m)$ $:=$ $\left\Vert y_{t}-y_{t}^{\prime }(m)\right\Vert _{p}$,
with accumulations $\Theta _{i,t}^{(p)}:=\sum_{m=0}^{\infty }\theta
_{i,t}^{(p)}(m)$ and $\tilde{\Theta}_{t}^{(p)}:=\sum_{m=0}^{\infty }\tilde{%
\theta}_{t}^{(p)}(m)$.

Write compactly $\left\Vert \check{z}\right\Vert _{n,p}$ $:=$ $\max_{1\leq
i\leq k_{n},1\leq t\leq n}\{\left\Vert x_{i,t}\right\Vert _{p}$ $\vee $ $%
\left\Vert y_{t}\right\Vert _{p}\}.$ In order to yield a clear upper bound
on $k_{n}$ that can be easily understood in terms of heterogeneity and
dependence decay, assume as before $\theta _{i,t}^{(p)}(m)$ $\vee $ $\tilde{%
\theta}_{t}^{(p)}(m)$ $\leq $ $d_{i,t}^{(p)}\psi _{i,m}$, where $\max_{i\in 
\mathbb{N}}\psi _{i,m}$ $=$ $O(m^{-\lambda -\iota })$ for some size $\lambda 
$ $\geq $ $1$, $\psi _{i,0}$ $=$ $1$. As with mixingales and near epoch
dependent sequences, by Minkowski's inequality logically $d_{i,t}^{(p)}$ $%
\leq $ $K\left\Vert \check{z}\right\Vert _{n,p}$ \citep[cf.][]{Davidson1994}%
. Then $\lambda $ $\geq $ $1$ yields $\max_{1\leq i\leq k_{n},1\leq t\leq
n}\{\Theta _{i,t}^{(p)}\vee \tilde{\Theta}_{t}^{(p)}\}$ $\leq $ $K\left\Vert 
\check{z}\right\Vert _{n,p}\sum_{m=1}^{\infty }m^{-\lambda -\iota }$ $=$ $%
K\left\Vert \check{z}\right\Vert _{n,p}$.

Define $\mathcal{H}_{i}$ $:=$ $\mathbb{E}\tilde{x}_{i,t}\tilde{x}%
_{i,t}^{\prime }$.

\begin{assumption}
\label{assum:marg-reg} \ \ \medskip \newline
$a$. $(x_{i,t},y_{t})$ are covariance stationary, and $\mathcal{L}_{p}$%
-physical dependent for some $p$ $\geq $ $4$, size $\lambda $ $\geq $ $1$,
and $\lim \sup_{n\rightarrow \infty }\left\Vert \check{z}\right\Vert _{n,p}$ 
$\leq $ $ap^{b}$ for some finite $a>0$ and $b$ $\in $ $[0,\infty )$.\medskip 
$\newline
b$. For each $i$ there exists a unique $\beta _{i}$ in the interior of $%
\mathcal{B}_{i}$ such that $\mathbb{E}(y_{t}$ $-$ $\beta _{i}^{\prime
}x_{i,t})x_{i,t}$ $=$ $0$.$\medskip $\newline
$c.$ $\lim \inf_{n\rightarrow \infty }\inf_{\lambda ^{\prime }\lambda =1}%
\mathbb{E}(1/\sqrt{n}\sum_{t=1}^{n}\lambda ^{\prime }\mathcal{H}_{i}^{-1}%
\tilde{x}_{i,t}v_{i,t}\lambda )^{2}$ $>$ $0$ for each $i$; $%
1/n\sum_{t=1}^{n}(x_{i,t}$ $-$ $\bar{x}_{i,n})^{2}$ $>$ $0$ $a.s.$
uniformly; $\mathbb{E}(x_{i,t}$ $-$ $\mathbb{E}x_{i,t})^{2}$ $>$ $0$
uniformly.
\end{assumption}

\begin{remark}
\normalfont($a$) restricts tail thickness prompting different exponential
bounds on $k_{n}$ based on $b$ $>$ $0$. Tails are sub-exponential when $b$ $%
\leq $ $1$, while moments grow too fast to be sub-exponential when $b$ $>$ $%
1 $ \citep[cf.][Proposition
2.7.1]{Vershynin2018}. ($b$) is a standard identification condition. ($c$)
implies $(\mathbb{E}\tilde{x}_{i,t}\tilde{x}_{i,t}^{\prime })^{-1}$ and $%
(1/n\sum_{t=1}^{n}\tilde{x}_{i,t}\tilde{x}_{i,t}^{\prime })^{-1}$ exist
uniformly in $i$ ($\lim \inf_{n\rightarrow \infty }\min_{1\leq i\leq
k_{n}}\inf_{\lambda ^{\prime }\lambda =1}\lambda \mathcal{H}_{i}\lambda $ $>$
$0$). The first item in ($c$), non-degeneracy,\ holds when $\boldsymbol{Z}%
_{n,i}$ $=$ $[x_{i,1}v_{i,1},$ $...,$ $x_{i,n}v_{i,n}]^{\prime }$ satisfies $%
\lim \inf_{n\rightarrow \infty }\inf_{\lambda ^{\prime }\lambda =1}\mathbb{E}%
(\lambda ^{\prime }\boldsymbol{Z}_{n,i})^{2}$ $>$ $0$, a conventional
positive definiteness property.
\end{remark}

We now present the main results. First, a HD first order approximation that
exploits Lemma \ref{lm:phys_dep_Lp} and max-WLLN Theorem \ref%
{thm:max_LLN_phys_dep}. See \citet[Appendix
E]{sm_max_LLN} for proofs. Write compactly%
\begin{equation}
z_{i,t}:=\left[ 0,1\right] \mathcal{H}_{i}^{-1}\tilde{x}_{i,t}v_{t}.
\label{zit}
\end{equation}

\begin{lemma}
\label{lm:mr_ng}Let Assumption \ref{assum:marg-reg} and $H_{0}$ hold, and
assume $\{x_{i,t},y_{t}\}$ are $\mathcal{L}_{p}$-physical dependent, $p$ $%
\geq $ $4$, with size $\lambda $ $\geq $ $1$. Then for any $\ln (k_{n})$ $=$ 
$o(n^{1/4})$%
\begin{equation*}
\left\vert \max_{1\leq i\leq k_{n}}\left\vert \sqrt{n}\hat{\phi}%
_{i}\right\vert -\max_{1\leq i\leq k_{n}}\left\vert \frac{1}{n^{1/2}}%
\sum_{t=1}^{n}z_{i,t}\right\vert \right\vert \overset{p}{\rightarrow }0
\end{equation*}
\end{lemma}

For a Gaussian approximation define $\sigma _{n,i}^{2}$ $:=$ $\mathbb{E}(1/%
\sqrt{n}\sum_{t=1}^{n}z_{i,t})^{2}$. Conditions imposed in Assumption \ref%
{assum:marg-reg} yield uniformly $\sigma _{n}^{2}$ $\in $ $(0,\infty )$. The
lower bound is ($c$). The upper bound is due to ($a$): By the proof of Lemma
E.1 in \cite{sm_max_LLN}, $\{x_{i,t}v_{t}\}$ is $\mathcal{L}_{2}$-physical
dependent when $\{x_{i,t},y_{t}\}$ are $\mathcal{L}_{4}$-physical dependent.
Hence by Lemma \ref{lm:phys_dep_Lp}.a $\max_{1\leq i\leq k_{n}}\left\Vert 1/%
\sqrt{n}\sum_{t=1}^{n}x_{i,t}v_{t}\right\Vert _{2}$ $\leq $ $2\max_{1\leq
i\leq k_{n},1\leq t\leq n}\{\Theta _{i,t}^{(4)}\vee \tilde{\Theta}%
_{t}^{(4)}\}$ $=$ $O(1)$ under uniform $\mathcal{L}_{4}$-boundedness $%
\left\Vert \check{z}\right\Vert _{n,4}$ $=$ $O(1)$, ruling out unbounded
fourth moment heterogeneity.

Now let $\{\mathcal{Z}_{n,i}$ $:$ $1$ $\leq $ $i$ $\leq $ $k_{n}\}_{n\geq 1}$
be a Gaussian array, $\mathcal{Z}_{n,i}$ $\sim $ $\mathcal{N}(0,\sigma
_{n,i}^{2})$, define%
\begin{equation*}
\rho _{n}:=\sup_{c>0}\mathbb{P}\left( \left\vert \max_{1\leq i\leq
k_{n}}\left\vert \frac{1}{\sqrt{n}}\sum_{t=1}^{n}z_{i,t}\right\vert
-\max_{1\leq i\leq k_{n}}\left\vert \mathcal{Z}_{n,i}\right\vert \right\vert
>c\right) ,
\end{equation*}%
and recall $b$ $>$ $0$ in Assumption \ref{assum:marg-reg}.a.

\begin{lemma}
\label{lm:mr_gauss}Let Assumption \ref{assum:marg-reg} and $H_{0}$ hold.
Assume $\{x_{i,t},y_{t}\}$ are $\mathcal{L}_{p}$-physical dependent, $p$ $%
\geq $ $4$, with size $\lambda $ $>$ $2$. Then for any $\{k_{n}\}$
satisfying $k_{n}$ $\rightarrow $ $\infty $ and $\ln (k_{n})$ $=$ $%
o(n^{g(b,\lambda )})$ where $g(b,\lambda )$ $:=$ $\frac{\lambda }{8+2\lambda 
}\frac{1}{(7/6)\vee (1+b)}$, we have $\rho _{n}$ $\rightarrow $ $0$.
Moreover $\max_{1\leq i\leq k_{n}}|1/\sqrt{n}\sum_{t=1}^{n}z_{i,t}|$ $%
\overset{d}{\rightarrow }$ $\max_{i\in \mathbb{N}}\left\vert \mathcal{Z}%
_{i}\right\vert $ for some Gaussian process $\{\mathcal{Z}_{i}\}$, $\mathcal{%
Z}_{i}$ $\sim $ $N(0,\sigma _{i}^{2})$ with $\sigma _{i}^{2}$ $=$ $%
\lim_{n\rightarrow \infty }\sigma _{n,i}^{2}$ $\in $ $(0,\infty )$.
\end{lemma}

\begin{remark}
\normalfont$k_{n}$ depends on tail conditions, memory decay, and
heterogeneity. As $\lambda $ $\searrow $ $2$ (far from independent) and $b$ $%
=$ $4$ (non-sub-exponential tails) then $\ln (k_{n})$ $=$ $o(n^{1/30})$.
Conversely, as $\lambda $ $\rightarrow $ $\infty $ (approaching geometric
memory/independence) with sub-exponential tails $b$ $=$ $1/6$ we have $%
g(b,\lambda )$ $\rightarrow $ $\frac{1}{2}\frac{1}{7/6}$ hence $\ln (k_{n})$ 
$=$ $o(n^{3/7})$.
\end{remark}

\begin{remark}
\normalfont The proof exploits HD Gaussian approximation Theorem 3($ii$) in 
\cite{ChangChenWu2024}. They propose two results: the first Theorem 3($i$)
supposedly imposing only their Condition 3 nondegeneracy, and the second
Theorem 3($ii$) imposing also their Condition 1 sub-exponential tails.
However, the \emph{dependence adjusted norms} that they exploit to bound $%
k_{n}$, based on ideas in \cite{WuWu2016}, only make sense when all moments
exist, e.g. $\lim \sup_{n\rightarrow \infty }\left\Vert \check{z}\right\Vert
_{n,p}$ $\leq $ $ap^{b}$ for some $b$ $>$ $0$. See especially 
\citet[Section 2.3, cf. eq.
(2.21)]{WuWu2016}.
\end{remark}

Lemmas \ref{lm:mr_ng} and \ref{lm:mr_gauss} imply the main result for the
max-test statistic $\max_{1\leq i\leq k_{n}}|\sqrt{n}\hat{\theta}_{i,n}|$.

\begin{theorem}
\label{thm:mr_gauss}Let Assumption \ref{assum:marg-reg} and $H_{0}$ hold.
Assume $\{x_{i,t},y_{t}\}$ are $\mathcal{L}_{p}$-physical dependent, $p$ $%
\geq $ $4$, with size $\lambda $ $>$ $2$. Then $\max_{1\leq i\leq k_{n}}|%
\sqrt{n}\hat{\theta}_{i,n}|$ $\overset{d}{\rightarrow }$ $\max_{i\in \mathbb{%
N}}|\mathcal{Z}_{i}|$ for any $\{k_{n}\}$ satisfying $\ln (k_{n})$ $=$ $%
o(n^{s(b,\lambda )})$ where by tail decay case%
\begin{eqnarray}
&&\text{if }b\in (0,1/6]\text{ then }s(b,\lambda )=\left\{ 
\begin{array}{ll}
\frac{1}{4} & \text{if }\lambda \geq \frac{28}{5} \\ 
\frac{\lambda }{8+2\lambda }\frac{1}{(7/6)\vee (1+b)} & \text{if }2<\lambda <%
\frac{28}{5}%
\end{array}%
\right.  \label{s(b,l)} \\
&&\text{if }b\in (1/6,1)\text{ then }s(b,\lambda )=\left\{ 
\begin{array}{ll}
\frac{1}{4} & \text{if }\lambda \geq \frac{4\left( 1+b\right) }{1-b} \\ 
\frac{\lambda }{8+2\lambda }\frac{1}{(7/6)\vee (1+b)} & \text{if }2<\lambda <%
\frac{4\left( 1+b\right) }{1-b}%
\end{array}%
\right.  \notag \\
&&\text{if }b\geq 1\text{ then }s(b,\lambda )=\frac{\lambda }{8+2\lambda }%
\frac{1}{1+b}\text{ }\forall \lambda >2.  \notag
\end{eqnarray}
\end{theorem}

\begin{remark}
\normalfont If $b$ $<$ $1/6$ then tails are sub-exponential %
\citep[cf.][Proposition 2.7.1]{Vershynin2018} and $\ln (k_{n})$ $=$ $%
o(n^{1/4})$ if memory decay is fast enough $\lambda $ $\geq $ $28/5$. If,
e.g., $b$ $<$ $1/6$ with slower decay $\lambda $ $=$ $4$ then $\ln (k_{n})$ $%
=$ $o(n^{3/14})$, a slower rate. In the intermediate range $b$ $\in $ $%
(1/6,1)$ tails are still sub-exponential, but $\ln (k_{n})$ $=$ $o(n^{1/4})$
when $\lambda $ $\geq $ $4/(\frac{2}{1+b}$ $-$ $1)$ $\searrow $ $28/5$ as $b$
$\searrow $ $1/6$: thinner tails allow for slower memory decay. Finally, $b$ 
$\geq $ $1$ allows for non-sub-exponential tails, yielding only $\ln (k_{n})$
$=$ $o(n^{\frac{\lambda }{8+2\lambda }\frac{1}{1+b}})$. If, e.g., $b$ $=$ $2$
then $\ln (k_{n})$ $=$ $o(n^{\frac{\lambda }{4+\lambda }\frac{1}{6}})$ where 
$\frac{\lambda }{4+\lambda }\frac{1}{6}$ $\searrow $ $\frac{1}{18}$ as $%
\lambda \searrow 2$ (far from independence) and $\frac{\lambda }{4+\lambda }%
\frac{1}{6}$ $\nearrow $ $1/6$ as $\lambda $ $\rightarrow $ $\infty $
(independence/geometric decay). In the hairline case $b$ $=$ $1$ tails are
sub-exponential, and $\ln (k_{n})$ $=$ $o(n^{\frac{\lambda }{4+\lambda }%
\frac{1}{4}})$ where $\frac{\lambda }{4+\lambda }\frac{1}{4}$ $\searrow $ $%
\frac{1}{12}$ as $\lambda \searrow 2$ and $\frac{\lambda }{4+\lambda }\frac{1%
}{4}$ $\nearrow $ $1/4$ as $\lambda $ $\rightarrow $ $\infty $.
\end{remark}

\section{Application \#3: testing parametric restrictions\label%
{sec:applic:param_test}}

Our final application combines methods in \cite{CattaneoJanssonNewey2018}
and \cite{Hill_2025_maxtest}. Consider a triangular array of observations $%
\{w_{n,t},x_{n,t},y_{n,t}$ $:$ $1$ $\leq $ $t$ $\leq $ $n\}_{n\geq 1}$ with
dependent variable $y_{n,t}$, and covariates $(w_{n,t},x_{n,t})$ of
dimensions $(k_{n},k_{\theta })$. The model is%
\begin{equation}
y_{n,t}=\delta _{n,0}^{\prime }w_{n,t}+\theta _{n,0}^{\prime }x_{n,t}+u_{n,t}
\label{model}
\end{equation}%
with error term $u_{n,t}$. Let $\mathbb{E}(y_{n,t}$ $-$ $\delta
_{n,0}^{\prime }w_{n,t}$ $-$ $\theta _{n,0}^{\prime
}x_{n,t})[w_{n,t}^{\prime },x_{n,t}^{\prime }]^{\prime }$ $=$ $0$ for unique 
$[\delta _{n,0}^{\prime },\theta _{n,0}^{\prime }]^{\prime }$. The model may
be pseudo-true in the sense $\mathbb{P}(\mathbb{E}_{\boldsymbol{w}_{nt},%
\boldsymbol{x}_{nt}}(y_{n,t}$ $-$ $\delta _{n}^{\prime }w_{n,t}$ $-$ $\theta
_{n}^{\prime }x_{n,t})$ $=$ $0)$ $<$ $1$ $\forall (\delta _{n},\theta _{n})$%
, where, e.g., $\boldsymbol{w}_{nt}$ $:=$ $\{w_{n,1},...,w_{n,t}\}$. The
array representation covers many cases in social sciences and statistics,
including $(i)$ linear models with increasing dimension via $k_{n}$; $(ii)$\
models with basis expansions of flexible functional forms, like partially
linear models $y_{t}$ $=$ $g(z_{t})$ $+$ $\theta _{0}^{\prime }x_{t}$ $+$ $%
u_{t}$ for some unknown measurable function $g$, and regressor set $z_{t}$;
and ($iii$) models with many dummy variables, e.g. panel models with
multi-way fixed effects. Cf. \citet[Section
3.3]{CattaneoJanssonNewey2018}.

\cite{CattaneoJanssonNewey2018} partial out the HD $\delta _{n,0}$ in order
to estimate the fixed low dimensional $\theta _{n,0}$, and propose HAC
methods for robust inference with arbitrary in-group dependence with finite
fixed group size. We consider the \textit{converse} problem in a far broader
setting.

We test the HD parameter $H_{0}$ $:$ $\delta _{n,0}$ $=$ $0$ vs. $H_{1}$ $:$ 
$\delta _{n,0}$ $\neq $ $0$ by partialling out $\theta _{n,0}$, but exploit
many low dimensional or \textit{parsimonious} models under $H_{0}$ as in 
\cite{Hill_2025_maxtest} to yield $\hat{\delta}_{i,n}$. We then use a
max-statistic $\max_{1\leq i\leq k_{n}}\sqrt{n}|\hat{\delta}_{i,n}|$ for
testing $H_{0}$. Partialling out is useful when $k_{\theta }$ is large
relative to $n$, or consistency of $\hat{\theta}_{n}$ is not guaranteed
(e.g. in panel settings with many fixed effects). Although we do not allow
for $x_{n,t}$ to be high dimensional, we anticipate the following will
extend to that case. The parsimonious approach alleviates the need for
regularization and therefore sparsity, as in de-biased Lasso, and is
significantly (potentially massively) faster to compute than de-biased Lasso %
\citep[see][]{Hill_2025_maxtest}. Moreover, a max-statistic sidesteps HAC
estimation and therefore inversion of a large dimension matrix, both of
which may lead to poor inference. See \cite{HillMotegi2020}, \cite%
{HillGhyselsMotegi2020} and \cite{Hill_2025_maxtest} for demonstrations of
asymptotic max-test superiority in models with potentially very many
parameters.

The partialled-out $\hat{\delta}_{n}$ is derived as follows. First, estimate
parsimonious models%
\begin{equation}
y_{n,t}=\delta _{i,n}^{\ast }w_{i,n,t}+\theta _{i,n}^{\ast \prime
}x_{n,t}+e_{i,n,t},\text{ }i=1,...,k_{n}.  \label{parsim}
\end{equation}%
Define $\delta _{n}^{\ast }$ $:=$ $[\delta _{i,n}^{\ast }]_{i=1}^{k_{n}}$.
By Theorem 2.1 in \cite{Hill_2025_maxtest} $\delta _{n}^{\ast }$ $=$ $0$ 
\textit{if and only if} $\delta _{n,0}$ $=$ $0$, hence $\theta _{n}^{\ast
}=\theta _{n,0}$ and $e_{i,n,t}$ $=$ $u_{n,t}$\ $\forall i$ under $H_{0}$.
Thus, we need only estimate each model in (\ref{parsim}) to yield some $\hat{%
\delta}_{n}$ $=$ $[\hat{\delta}_{i,n}]_{i=1}^{k_{n}}$ and thereby test $%
H_{0} $.

Define an $l_{2}$ orthogonal projection matrix $\mathcal{M}_{n}$ $:=$ $I_{n}$
$-$ $\boldsymbol{x}_{n}(\boldsymbol{x}_{n}^{\prime }\boldsymbol{x}_{n})^{-1}%
\boldsymbol{x}_{n}^{\prime }$ $\in $ $\mathbb{R}^{n\times n}$ with identity
matrix $I_{n}$, where $\boldsymbol{x}_{n}$ $:=$ $[x_{n,1},...,x_{n,n}]^{%
\prime }$. After partialling out based on a projection onto the linear space
spanned by $x_{n,t}$, yielding $\mathcal{M}_{n}\boldsymbol{y}_{n}$ $=$ $%
\delta _{i,n}^{\ast }\boldsymbol{\hat{v}}_{i,n}$ $+$ $\mathcal{M}_{n}%
\boldsymbol{e}_{i,n}$, where $\boldsymbol{\hat{v}}_{i,n}$ $:=$ $\mathcal{M}%
_{n}\boldsymbol{w}_{i,n}$ $\in $ $\mathbb{R}^{n\times 1}$, the estimator of $%
\delta _{i,n}^{\ast }$ reduces to%
\begin{equation*}
\hat{\delta}_{i,n}=\left( \boldsymbol{\hat{v}}_{i,n}^{\prime }\boldsymbol{%
\hat{v}}_{i,n}\right) ^{-1}\boldsymbol{\hat{v}}_{i,n}^{\prime }\boldsymbol{y}%
_{n}.
\end{equation*}%
The test statistic is $\mathcal{T}_{n}$ $:=$ $\max_{1\leq i\leq k_{n}}\sqrt{n%
}|\hat{\delta}_{i,n}|$. We assume below $\mathbb{E}(\gamma ^{\prime }%
\boldsymbol{w}_{n,t}+\delta ^{\prime }\boldsymbol{x}_{n,t})^{2}$ $>$ $0$
uniformly in ($n,t,\gamma ^{\prime }\gamma $ $=$ $\delta ^{\prime }\delta $ $%
=$ $1$), hence $\inf_{1\leq i\leq k_{n}}\{\boldsymbol{\hat{v}}_{i,n}^{\prime
}\boldsymbol{\hat{v}}_{i,n}/n\}$ $>$ $0$ $awp1$ \citep[Lemma F.3]{sm_max_LLN}%
. Thus logically $w_{n,t}$ and $x_{n,t}$\ cannot be perfectly linearly
related.

We assume stochastic components $\{w_{n,t},x_{n,t},u_{n,t}\}$ are $\rho _{i}$%
-Lipschitz Markov processes in order to focus ideas, implying both $\tau
^{(p)}$-mixing and $\mathcal{L}_{p}$-physical dependence. Define 
\begin{equation*}
\mathcal{\hat{Z}}_{n,i}:=\left( \frac{1}{n}\sum_{t=1}^{n}\mathbb{E}%
v_{i,n,t}v_{i,n,t}^{\prime }\right) ^{-1}\frac{1}{\sqrt{n}}%
\sum_{t=1}^{n}v_{i,n,t}u_{n,t}\text{ and }\sigma _{n,i}^{2}:=\mathbb{E}%
\mathcal{\hat{Z}}_{n,i}^{2}.
\end{equation*}

\begin{assumption}
\label{assume:LD_test}Let $z_{i,n,t}$ $\in $ $\{w_{i,n,t},x_{i,n,t},u_{n,t}%
\} $.$\medskip $\newline
$a.$ Each $z_{i,n,t}$ $=$ $f_{z_{i}}(z_{i,n,t-1})$ $+$ $\epsilon _{i,t}$,
for $\rho _{z_{i}}$-Lipschitz $f_{z_{i}}(\cdot )$, $\rho _{z_{i}}$ $\in $ $%
(0,e^{-\alpha _{z_{i}}}]$ for some $\alpha _{z_{i}}$ $\in $ $(0,\infty )$,
serially iid $\epsilon _{i,t}$, and $\mathcal{L}_{p}$-bounded $\{\epsilon
_{i,t},z_{i,n,t}\}$ for some $p$ $\geq $ $4$.\medskip \newline
$b$. $\max_{1\leq i\leq k_{n},1\leq t\leq n}\mathbb{P}(|z_{i,n,t}|$ $>$ $z)$ 
$\leq $ $a_{z}\exp \{b_{z}z^{-\gamma _{z}}\}$ $\forall n$ for some $%
(a_{z},b_{z},\gamma _{z})$ $\in $ $(0,\infty )$.\medskip \newline
$c$. $\lim \inf_{n\rightarrow \infty }\inf_{\lambda ^{\prime }\lambda
=1}\{\lambda ^{\prime }\boldsymbol{x}_{n}^{\prime }\boldsymbol{x}_{n}\lambda
/n$ $\wedge $ $\lambda ^{\prime }\boldsymbol{w}_{i,n}^{\prime }\boldsymbol{w}%
_{i,n}\lambda /n\}$ $>$ $0$ $a.s.$; and $\mathbb{E}(\gamma ^{\prime }w_{n,t}$
$+$ $\delta ^{\prime }x_{n,t})^{2}$ $>$ $0$ uniformly over $(n,t,\delta
^{\prime }\delta $ $=$ $1,\gamma ^{\prime }\gamma $ $=$ $1)$.\medskip 
\newline
$d$. $\sigma _{n,i}^{2}$ $\in $ $[K,\infty )$ for some $K$ $>$ $0$ uniformly
in $(i,n)$.
\end{assumption}

\begin{remark}
\normalfont($a$) implies $z_{i,n,t}$ $=$ $g_{z_{i}}(\epsilon _{i,t},\epsilon
_{i,t-1},...)$ for measurable $g_{z_{i}}$ %
\citep[e.g.][]{DiaconisFreedman1999}, and is geometrically $\tau ^{(p)}$%
-mixing by Example \ref{ex:p-Lip Mark}, and (therefore) geometrically
uniformly $\mathcal{L}_{p}$-physical dependent by Lemma C.4 in \cite%
{sm_max_LLN} and linkages in \cite{Hill2025_mixg}. See also 
\citet[p.
14152]{Wu2005}. Thus intertemporal dependence decays geometrically fast. We
can easily allow for arbitrary group-wise dependence for finite,
heterogeneously sized groups by assuming $z_{i,n,t}$ $=$ $r_{i,n,t}$ $+$ $%
s_{i,n,t}$ where $\rho _{z_{i}}$-Lipschitz $r_{i,n,t}$ $=$ $%
f_{r_{i}}(r_{i,n,t-1})$ $+$ $\epsilon _{i,t}$, and $s_{i,n,t}$ is $M_{s_{i}}$%
-dependent for finite heterogeneous $M_{s_{i}}$: $z_{i,n,t}$ is still
geometrically $\mathcal{L}_{p}$-physical dependent. Indeed, $M_{s_{i}}$%
-dependence can be replaced with arbitrary dependence in arbitrary groups
(e.g. $t_{1}^{i},...,t_{M_{s_{i}}}^{i}$), nesting the Assumption 1
independence setting in \cite{CattaneoJanssonNewey2018}. We work under ($a$)
instead to save notation.
\end{remark}

\begin{remark}
\normalfont($b$) ensures both a max-LLN and HD central limit theorem apply,
and implies $z_{i,n,t}$ are uniformly $\mathcal{L}_{p}$-bounded $\forall p$ $%
\geq $ $1$. In ($c$), $\mathbb{E}(\gamma ^{\prime }w_{n,t}+\delta ^{\prime
}x_{n,t})^{2}$ $>$ $0$ uniformly ensures positive definiteness $%
\inf_{\lambda ^{\prime }\lambda =1}\lambda ^{\prime }(1/n\sum_{t=1}^{n}$ $%
\mathbb{E}w_{n,t}w_{n,t}^{\prime })\lambda $ $>$ $0$ $\forall n$ and $%
\inf_{\lambda ^{\prime }\lambda =1}\lambda ^{\prime }(1/n\sum_{t=1}^{n}%
\mathbb{E}x_{n,t}x_{n,t}^{\prime })\lambda $ $>$ $0$ $\forall n$, and rules
out deviant cross-correlations, ensuring $\boldsymbol{\hat{v}}_{i,n}^{\prime
}\boldsymbol{\hat{v}}_{i,n}/n$ is positive definite $awp1$. Non-degeneracy ($%
d$) is standard: see remarks following Assumptions \ref{assum:maxcor} and %
\ref{assum:marg-reg}.
\end{remark}

\begin{remark}
\normalfont Our assumptions differ from 
\citet[Assumptions
1-3]{CattaneoJanssonNewey2018}. They impose cross-group independence with
finite heterogeneous group sizes, and allow for heteroscedasticity. They
need (\ref{model}) to be very close to the true model by several measures
(see their Assumption 3; e.g. $\mathbb{E}(\mathbb{E}_{\boldsymbol{w}_{nt},%
\boldsymbol{x}_{nt}}u_{n,t})^{2}$ $=$ $o(1/n)$). We allow for
nonstationarity, (\ref{model}) need not be the true model, and within-group
dependence can be arbitrary as discussed above. Nonstationarity allows for
heteroscedasticity and other forms of heterogeneity, and a max-test allows
us to by-pass covariance matrix estimation \emph{entirely} (it is \emph{ipso
facto} heteroscedasticity robust). Of course, they partial-out the high
dimensional term and estimate one model, while we $(i)$ partial out the
fixed (low) dimensional term, $(ii)$ estimate many low dimension models, and
therefore $(iii)$ use an entirely different asymptotic theory to test a HD
hypothesis.
\end{remark}

Let $\{\mathcal{Z}_{n,i}$ $:$ $1$ $\leq $ $i$ $\leq $ $k_{n}\}_{n\geq 1}$ be
Gaussian, $\mathcal{Z}_{n,i}\sim $ $\mathcal{N}(0,\sigma _{n,i}^{2})$ with $%
\sigma _{n,i}^{2}:=\mathbb{E}\mathcal{\hat{Z}}_{n,i}^{2}$, and define%
\begin{equation*}
\rho _{n}:=\sup_{c>0}\mathbb{P}\left( \left\vert \max_{1\leq i\leq
k_{n}}\left\vert \mathcal{\hat{Z}}_{n,i}\right\vert -\max_{1\leq i\leq
k_{n}}\left\vert \mathcal{Z}_{n,i}\right\vert \right\vert >c\right) .
\end{equation*}%
We require a moment growth parameter $b$ developed in 
\citet[Appendix
F]{sm_max_LLN}, similar to Assumption \ref{assum:marg-reg}.a. By Lemma F.4
each $z_{n,t}(i,j)$ $\in $ $\{w_{i,n,t}x_{j,n,t}$ $-$ $\mathbb{E}%
w_{i,n,t}x_{j,n,t}$, $x_{i,n,t}x_{j,n,t}$ $-$ $\mathbb{E}x_{i,n,t}x_{j,n,t}$
and $x_{i,n,t}u_{n,t}$ satisfies%
\begin{equation*}
\sup_{p\geq 2}\frac{1}{p^{b}}\left\{ \sup_{m\geq 0}\left( m+1\right)
\max_{1\leq i,j\leq k_{n},1\leq t\leq n}\left\Vert
z_{n,t}(i,j)-z_{n,t}^{\prime }(i,j)\right\Vert _{p/2}\right\} \leq K
\end{equation*}%
for some $b$ $>$ $0$ that depends only on the Assumption \ref{assume:LD_test}%
.b tail parameters. If $b$ $\leq $ $1$ then $z_{n,t}(i,j)$ have
sub-exponential tails. The following omnibus result characterizes first
order and Gaussian approximations, and the max-statistic limit. Max-WLLN
Theorem \ref{thm:max_LLN_phys_dep} is utilized in the proof. The proof is in %
\citet[Appendix F]{sm_max_LLN}.

\begin{theorem}
\label{thm:LD_test} Let Assumption \ref{assume:LD_test} and $H_{0}$
hold.\medskip \newline
$a.$ \emph{(Non-Gaussian Approximation)}. $\left\vert \sqrt{n}\max_{1\leq
i\leq k_{n}}\left\vert \hat{\delta}_{i,n}\right\vert -\max_{1\leq i\leq
k_{n}}\left\vert \mathcal{\hat{Z}}_{n,i}\right\vert \right\vert $ $=$ $%
o_{p}(1)$ for any $\{k_{n}\}$, $\ln (k_{n})$ $=$ $o(n^{1/4})$.\medskip 
\newline
$b.$ \emph{(Gaussian Approximation). }$\rho _{n}$ $\rightarrow $ $0$ for any 
$\{k_{n}\}$, $\ln (k_{n})$ $=$ $o(n^{1/[2(1+\varphi )]})$.\medskip \newline
$c$. $\mathcal{T}_{n}$ $\overset{d}{\rightarrow }\max_{i\in \mathbb{N}}%
\mathcal{Z}_{i}$ where $\mathcal{Z}_{i}$ $\sim $ $N(0,\lim_{n\rightarrow
\infty }\sigma _{n,i}^{2})$ for any $\{k_{n}\}$ satisfying $\ln (k_{n})$ $=$ 
$o(n^{s(b)})$ where $s(b)$ $=$ $\lim_{\lambda \rightarrow \infty
}s(b,\lambda )$, $s(b,\lambda )$\ is depicted in (\ref{s(b,l)}), and $b$ is
defined above. Thus $s(b)$ $=$ $1/4$ if $b$ $\in $ $(0,1)$ and $s(b)$ $=$ $%
1/[2(1+b)]$ if $b$ $\geq $ $1$.
\end{theorem}

\section{Conclusion\label{sec:conclud}}

We present weak laws of large numbers for the maximum sample average
\linebreak $\max_{1\leq i\leq k_{n}}|1/n\sum_{t=1}^{n}x_{i,n,t}|$ of a high
dimensional array $\{x_{i,n,t}$ $:$ $1$ $\leq $ $i$ $\leq $ $%
k_{n}\}_{t=1}^{n}$. We work under updated $\tau $-mixing and physical
dependence properties, while deriving new relational results. Certain
max-LLN's reveal a memory and dimension growth trade-off, depending on
nuances of the underlying dependence property. We work with and without
cross-coordinate dependence restrictions, where generally cross-coordinate
dependence can be wielded to achieve an improvement on $k_{n}$. The methods
are applied to independent sequences in the technical appendix in which we
also derive a max-SLLN. The results are applied to three settings: a
max-correlation white noise test; correlation screening under dependence and 
$k_{n}/n$ $\rightarrow $ $\infty $; and a high dimensional regression
parameter test under dependence.

As next steps, it would be interesting to $(i)$ extend the results to near
epoch dependent [NED] arrays which are nested under mixingales, or a spatial
setting; $(ii)$ study cross-coordinate dependence further in an attempt to
yield general results with applications; $(iii)$ extend the results to high
dimensional laws of iterated logarithm under dependence; $(iv)$ extend
results to uniform laws in high dimension; $(v)$ resolve difficulties
surrounding max-SLLN's for dependent sequences discussed in the technical
appendix. All such ideas are left for future consideration.

\setcounter{equation}{0} \renewcommand{\theequation}{{\thesection}.%
\arabic{equation}} \appendix

\section{Appendix}

\subsection{Technical proofs\label{app:proofs}}

\noindent \textbf{Proof of Lemma \ref{lm:tp_Bern}.} Under $\tau ^{(1)}$ and
mixing and tail decay (\ref{tao})-(\ref{ggg}) we have uniformly over $(i,n,t)
$ \citep[Theorem 1]{Merlevede_et_al_2011},%
\begin{eqnarray}
&&\mathbb{P}\left( \max_{1\leq l\leq n}\left\vert \frac{1}{n}%
\sum_{t=1}^{l}x_{i,n,t}\right\vert \geq \epsilon \right)   \label{t-ineq} \\
&&\text{ \ \ \ \ \ \ \ \ \ }\leq n\exp \left\{ -\mathcal{K}_{1}\epsilon
^{\gamma }n^{\gamma }\right\} +\exp \left\{ -\mathcal{K}_{2}\frac{\epsilon
^{2}n^{2}}{1+\mathcal{K}_{3}n}\right\} +\exp \left\{ -\mathcal{K}_{4}\frac{%
\epsilon ^{2}n^{2}}{n}e^{\frac{\mathcal{K}_{5}\left( \epsilon n\right)
^{\gamma (1-\gamma )}}{[\ln (\epsilon n)]^{\gamma }}}\right\} ,  \notag
\end{eqnarray}%
for some $\gamma $ $\in $ $(0,1)$. Since $n$ $\geq $ $4$ and $\epsilon $ $>$ 
$1/4$ by supposition, $\ln (\epsilon n)$ $>$ $0$. \cite{Merlevede_et_al_2011}
assume $d$ $=$ $\exp \{1\}$ in (\ref{sube}), but this can be generalized to
any $d$ $>$ $0$. Their proof, with coupling result Lemma C.2 in \cite%
{sm_max_LLN}, and arguments in \citet[Lemma
5]{Dedecker_Prieur_2004} and \citet[p. 460]{Merlevede_et_al_2011}, directly
imply (\ref{t-ineq}) holds under $\tau ^{(p)}$. Indeed $\max_{1\leq i\leq
k_{n}}\tau _{i,n}^{(1)}(m)\leq $ $\{\max_{1\leq i\leq k_{n}}\tau
_{i,n}^{(p)}(m)\}^{1/p}$ $\leq $ $a^{1/p}e^{-(b/p)m^{\gamma _{1}}}$ by
Lyapunov's inequality and (\ref{tao}). Hence arguments in 
\citet[proof of Theorem
1]{Merlevede_et_al_2011} go through with $(a,b)$ replaced with $(a^{1/p},b/p)
$. The upper bound in (\ref{t-ineq}) is not a function of $i$, thus proving (%
\ref{BE_t}). $\mathcal{QED}$.\bigskip \newline
\textbf{Proof of Theorem \ref{thm:max_LLN_t_mix}}. Jensen's inequality gives
a log-exp bound $\forall \lambda $ $>$ $0$,%
\begin{eqnarray}
\mathbb{E}\max_{1\leq i\leq k_{n}}\left\vert \frac{1}{n}%
\sum_{t=1}^{n}x_{i,n,t}\right\vert  &\leq &\frac{1}{\lambda }\ln \left( 
\mathbb{E}\exp \left\{ \lambda \max_{1\leq i\leq k_{n}}\left\vert \frac{1}{n}%
\sum_{t=1}^{n}x_{i,n,t}\right\vert \right\} \right)   \notag \\
&\leq &\frac{1}{\lambda }\ln \left( k_{n}\max_{1\leq i\leq k_{n}}\mathbb{E}%
\exp \left\{ \lambda \left\vert \frac{1}{n}\sum_{t=1}^{n}x_{i,n,t}\right%
\vert \right\} \right) .  \label{E1/nx}
\end{eqnarray}%
Furthermore, 
\begin{equation}
\mathbb{E}\exp \left\{ \lambda \left\vert \frac{1}{n}\sum_{t=1}^{n}x_{i,n,t}%
\right\vert \right\} =\int_{0}^{\infty }\mathbb{P}\left( \left\vert \frac{1}{%
n}\sum_{t=1}^{n}x_{i,n,t}\right\vert >\frac{1}{\lambda }\ln \left( u\right)
\right) du.  \label{E1nXIntP}
\end{equation}%
In (\ref{t-ineq}), cf. (\ref{BE_t}) in Lemma \ref{lm:tp_Bern}, because $%
\gamma $ $\in $ $(0,1)$ the first term trivially dominates the third, and
dominates the second for all $\epsilon $ $\geq $ $1$ and $n$ $\geq $ $n_{0}$%
, and finite $n_{0}$ $\geq $ $1$ depending on $(\mathcal{K}_{1},\mathcal{K}%
_{2},\mathcal{K}_{3},\gamma )$. Hence for some $\mathcal{K}$ depending on $(%
\mathcal{K}_{1},\mathcal{K}_{2},\mathcal{K}_{3},\gamma )$ that may be
different in different places, and for any $i$%
\begin{equation*}
\mathbb{P}\left( \max_{1\leq l\leq n}\left\vert \frac{1}{n}%
\sum_{t=1}^{l}x_{i,n,t}\right\vert \geq \epsilon \right) \leq 3n\exp \left\{
-\mathcal{K}\epsilon ^{\gamma }n^{\gamma }\right\} \text{ }\forall \epsilon
\geq 1\text{ and }n\geq n_{0}.
\end{equation*}%
Moreover, $3n\exp \{-\mathcal{K}\epsilon ^{\gamma }n^{\gamma }\}$ $\leq $ $%
\exp \{-(\mathcal{K}/2)\epsilon ^{\gamma }n^{\gamma }\}$ $\forall \epsilon $ 
$\geq $ $1$, $n$ $\geq $ $n_{1}$, and finite $n_{1}$ $\geq $ $1$. Therefore, 
$\forall n$ $\geq $ $n_{0}\vee n_{1}$, and any $\lambda $ $>$ $0$,%
\begin{eqnarray*}
\max_{1\leq i\leq k_{n}}\mathbb{E}\exp \left\{ \lambda \left\vert \frac{1}{n}%
\sum_{t=1}^{n}x_{i,n,t}\right\vert \right\}  &\leq &e+\max_{1\leq i\leq
k_{n}}\int_{e}^{\infty }\mathbb{P}\left( \left\vert \frac{1}{n}%
\sum_{t=1}^{n}x_{i,n,t}\right\vert >\frac{1}{\lambda }\ln \left( u\right)
\right) du \\
&\leq &e+\int_{e}^{\infty }\exp \left\{ -\frac{\mathcal{K}n^{\gamma }}{%
2\lambda ^{\gamma }}\left( \ln \left( u\right) \right) ^{\gamma }\right\} du
\\
&=&e+\frac{1}{\gamma }\int_{1}^{\infty }\frac{1}{v^{(\gamma -1)/\gamma }}%
\exp \left\{ v^{1/\gamma }-\frac{\mathcal{K}n^{\gamma }}{2\lambda ^{\gamma }}%
v\right\} dv \\
&\leq &e+\frac{1}{\gamma }\int_{1}^{\infty }\frac{1}{v^{(\gamma -1)/\gamma }}%
\exp \left\{ -\left( \frac{\mathcal{K}n^{\gamma }}{2\lambda ^{\gamma }}%
-1\right) v\right\} dv \\
&\leq &e+\frac{1}{\gamma }\int_{1}^{\infty }\exp \left\{ -\left( \frac{%
\mathcal{K}n^{\gamma }}{2\lambda ^{\gamma }}-1\right) v\right\} dv.
\end{eqnarray*}%
The second equality uses a change of variables $v$ $=$ $(\ln (u))^{\gamma }$ 
$\in $ $[0,\infty )$, the third inequality uses $\gamma $ $\geq $ $1$ from (%
\ref{ggg}), and the fourth uses $v$ $>$ $1$. Notice for all $v$ $>$ $1$, all 
$n$, some $\mathcal{\tilde{K}}$\ $\in $ $(0,\mathcal{K}/2)$ and any $\lambda 
$ $\leq $ $(\mathcal{K}/2-\mathcal{\tilde{K})}^{1/\gamma }n$, 
\begin{equation*}
\exp \left\{ -\left( \frac{\mathcal{K}n^{\gamma }}{2\lambda ^{\gamma }}%
-1\right) v\right\} \leq \exp \left\{ -\mathcal{\tilde{K}}\frac{n^{\gamma }}{%
\lambda ^{\gamma }}v\right\} .
\end{equation*}%
Therefore, $\forall n$ $\geq $ $n_{0}\vee n_{1}$ and any $\lambda $ $\leq $ $%
(\mathcal{K}/2-\mathcal{\tilde{K})}^{1/\gamma }n,$%
\begin{eqnarray*}
\max_{1\leq i\leq k_{n}}\mathbb{E}\exp \left\{ \lambda \left\vert \frac{1}{n}%
\sum_{t=1}^{n}x_{i,n,t}\right\vert \right\}  &\leq &e+\frac{1}{\gamma }%
\int_{1}^{\infty }\exp \left\{ -\mathcal{\tilde{K}}\frac{n^{\gamma }}{%
\lambda ^{\gamma }}v\right\} dv \\
&\leq &e+\frac{\lambda ^{\gamma }}{\gamma \mathcal{\tilde{K}}n^{\gamma }}%
\exp \left\{ -\mathcal{\tilde{K}}\frac{n^{\gamma }}{\lambda ^{\gamma }}%
\right\} \leq e+\frac{\lambda ^{\gamma }}{\gamma \mathcal{\tilde{K}}%
n^{\gamma }}.
\end{eqnarray*}%
Now use (\ref{E1/nx}) and take $\lambda $ $=$ $an$ for any $0$ $<$ $a$ $\leq 
$ $(\mathcal{K}/2-\mathcal{\tilde{K})}^{1/\gamma }$, and $\ln (k_{n})$ $=$ $%
o(n)$, to yield%
\begin{eqnarray*}
\mathbb{E}\max_{1\leq i\leq k_{n}}\left\vert \frac{1}{n}%
\sum_{t=1}^{n}x_{i,n,t}\right\vert  &\leq &\frac{1}{\lambda }\left\{ \ln
(k_{n})+\ln \left( e+\frac{\lambda ^{\gamma }}{\gamma \mathcal{\tilde{K}}%
n^{\gamma }}\right) \right\}  \\
&=&\frac{1}{an}\ln (k_{n})+\frac{1}{n^{a}}\ln \left( e+\frac{a^{\gamma
}n^{\gamma }}{\gamma \mathcal{\tilde{K}}n^{\gamma }}\right) =\frac{\ln
(k_{n})}{an}+\frac{1}{an}\ln \left( e+\frac{a^{\gamma }}{\gamma \mathcal{%
\tilde{K}}}\right) \rightarrow 0.
\end{eqnarray*}%
Hence $\mathcal{M}_{n}$ $\overset{\mathcal{L}_{1}}{\rightarrow }$ $0$
whenever $\ln (k_{n})$ $=$ $o(n)$.

Finally, the above arguments with $\lambda $ $=$ $\sqrt{n\ln \left(
k_{n}\right) }$\ and $\ln (k_{n})$ $=$ $O(n)$\ imply identically 
\begin{eqnarray*}
\mathbb{P}\left( \max_{1\leq i\leq k_{n}}\left\vert \frac{1}{\sqrt{n\ln
\left( k_{n}\right) }}\sum_{t=1}^{n}x_{i,n,t}\right\vert >c\right) &\leq &%
\frac{1}{c}\sqrt{\frac{n}{\ln \left( k_{n}\right) }}\frac{1}{\lambda }%
\left\{ \ln (k_{n})+\ln \left( e+\frac{1}{\gamma \mathcal{\tilde{K}}}\left( 
\frac{\lambda }{n}\right) ^{\gamma }\right) \right\} \\
&=&\frac{1}{c}\frac{1}{\ln \left( k_{n}\right) }\left\{ \ln (k_{n})+\ln
\left( e+\frac{1}{\gamma \mathcal{\tilde{K}}}\left( \sqrt{\frac{\ln \left(
k_{n}\right) }{n}}\right) ^{\gamma }\right) \right\} \\
&=&\frac{1}{c}\left\{ 1+\frac{1}{\ln \left( k_{n}\right) }\ln \left(
e+O\left( 1\right) \right) \right\} =O(1)\text{ \ }\forall c>0,
\end{eqnarray*}%
completing the proof. $\mathcal{QED}$.\medskip \newline
\textbf{Proof of Lemma \ref{lm:phys_dep_Lp}.} Write $\mathcal{Z}_{i,l}$ $:=$ 
$1/\sqrt{n}\sum_{t=1}^{l}x_{i,n,t}$.\medskip \newline
\textbf{Claim (a).} For similar arguments see 
\citet[Lemma
21]{JirakKostenberger2024} when $p$ $>$ $1$ and 
\citet[Theorem
2(\textit{i})]{Wu2005} when $p$ $\geq $ $2$. Recall $\xi _{i,t}$ $:=$ $%
\{\epsilon _{i,t},\epsilon _{i,t-1},..\}$.

Define $\mathcal{M}_{r,m}$ $:=$ $\sum_{l=1}^{m}y_{i,n,l}^{(r)}$ where $%
y_{i,n,l}^{(r)}$ $:=$ $\mathbb{E}(x_{i,n,l}|\xi _{i,l-r})$ $-$ $\mathbb{E}%
(x_{i,n,l}|\xi _{i,l-r-1})$. Then\ $\sum_{t=1}^{n}x_{i,n,t}$ $=$ $%
\sum_{r=0}^{\infty }\mathcal{M}_{r,n}$, hence by triangle and Minkowski
inequalities, and Doob's martingale inequality when $p$ $>$ $1$ 
\citep[e.g.][Theorem
2.2]{HallHeyde1980},

\begin{equation}
\left\Vert \max_{1\leq l\leq n}\left\vert \sum_{t=1}^{l}x_{i,n,t}\right\vert
\right\Vert _{p}\leq \sum_{r=0}^{\infty }\left\Vert \max_{1\leq l\leq
n}\left\vert \sum_{t=1}^{l}y_{i,n,t}^{(r)}\right\vert \right\Vert _{p}\leq 
\frac{p}{p-1}\sum_{r=0}^{\infty }\left\Vert
\sum_{t=1}^{n}y_{i,n,t}^{(r)}\right\Vert _{p}.  \label{phys_mart}
\end{equation}%
Define $\mathcal{A}_{i,n,j}^{(r)}$ $:=$ $\sigma
(y_{i,n,1}^{(r)},...,y_{i,n,j}^{(r)})$, hence $\mathcal{A}_{i,n,j}^{(r)}$ $=$
$\sigma (\xi _{i,j-r})$. Define \cite{Burkholder1973}'s constant $\mathcal{C}%
_{p}$ $:=$ $18p^{3/2}/(p$ $-$ $1)^{1/2}$, and $\mathcal{C}_{p}^{\prime }$ $%
:= $ $p\mathcal{C}_{p}/(p$ $-$ $1)$.\medskip

\textbf{Case 1 ($p$ $\in $ $(1,2)$)}. Apply Lemma 2.2 in \cite{Li2003} to$%
\left\Vert \sum_{l=1}^{n}y_{i,n,l}^{(r)}\right\Vert _{p}$, cf. 
\citet[Lemma
1]{WuShao2007}, to yield 
\begin{equation*}
\left\Vert \sum_{t=1}^{n}y_{i,n,t}^{(r)}\right\Vert _{p}\leq \mathcal{C}%
_{p}\left( \sum_{t=1}^{n}\left\Vert y_{i,n,t}^{(r)}\right\Vert
_{p}^{p}\right) ^{1/p}\leq \mathcal{C}_{p}n^{1/p}\max_{1\leq t\leq
n}\left\Vert y_{i,n,t}^{(r)}\right\Vert _{p}.
\end{equation*}%
Hence $\left\Vert \max_{1\leq t\leq n}\left\vert \mathcal{Z}%
_{i,t}\right\vert \right\Vert _{p}$ $\leq $ $\mathcal{C}_{p}^{\prime
}n^{1/p-1/2}\max_{1\leq t\leq n}\sum_{r=0}^{\infty }\left\Vert
y_{i,n,t}^{(r)}\right\Vert _{p}$. By definition $\left\Vert
y_{i,n,t}^{(r)}\right\Vert _{p}$ $=$ $\left\Vert \mathbb{E}(x_{i,n,t}|\xi
_{i,t-r})-\mathbb{E}(x_{i,n,t}|\xi _{i,t-r-1})\right\Vert _{p}$ $=:$ $\rho
_{i,n,t}^{(p)}(r)$, thus 
\begin{equation*}
\left\Vert \max_{1\leq t\leq n}\left\vert \mathcal{Z}_{i,t}\right\vert
\right\Vert _{p}\leq \mathcal{C}_{p}^{\prime }n^{1/p-1/2}\sum_{m=0}^{\infty
}\max_{1\leq t\leq n}\rho _{i,n,t}^{(p)}(m).
\end{equation*}%
Hence $\left\Vert \max_{1\leq t\leq n}\left\vert \mathcal{Z}%
_{i,t}\right\vert \right\Vert _{p}$ $\leq $ $\mathcal{C}_{p}^{\prime
}n^{1/p-1/2}\max_{1\leq t\leq n}\Theta _{i,n,t}^{(p)}$ by Theorem 2.1 in 
\cite{Hill2025_mixg}.\medskip

\textbf{Case 2 ($p$ }$\mathbf{\geq }$\textbf{\ $2$)}. The above argument
exploit's Burkholder's inequality and carries over to any $p$ $>$ $1$ %
\citep[see][Lemma 21]{JirakKostenberger2024}. We get a better a constant,
however, when $p$ $\geq $ $2$ based on arguments in \cite%
{DedeckerDoukhan2003}, cf. \citet[Chapt. 2.5]{Rio2017}. Apply Proposition 4
in \cite{DedeckerDoukhan2003} to $\left\Vert
\sum_{l=1}^{n}y_{i,n,l}^{(r)}\right\Vert _{p}$ in (\ref{phys_mart})\ to
yield 
\begin{eqnarray}
\left\Vert \sum_{l=1}^{n}y_{i,n,l}^{(r)}\right\Vert _{p} &\leq &\sqrt{2p}%
\left( \sum_{j=1}^{n}\max_{j\leq l\leq n}\left\Vert
y_{i,n,j}^{(r)}\sum_{m=j}^{l}\mathbb{E}\left( y_{i,n,m}^{(r)}|\mathcal{A}%
_{i,n,j}^{(r)}\right) \right\Vert _{p/2}\right) ^{1/2}  \label{Lpynt} \\
&=&\sqrt{2p}\left( \sum_{j=1}^{n}\left\Vert y_{i,n,j}^{(r)}\mathbb{E}\left(
y_{i,n,j}^{(r)}|\mathcal{A}_{i,n,j}^{(r)}\right) \right\Vert _{p/2}\right)
^{1/2}\leq \sqrt{2p}\sqrt{n}\max_{1\leq t\leq n}\left\Vert
y_{i,n,t}^{(r)}\right\Vert _{p}.  \notag
\end{eqnarray}%
The equality follows from the martingale difference property of $%
y_{i,n,m}^{(r)}$, measurability, and iterated expectations since%
\begin{eqnarray*}
\mathbb{E}\left( y_{i,n,m}^{(r)}|\mathcal{A}_{i,n,j}^{(r)}\right) &=&\mathbb{%
E}\left[ \mathbb{E}\left( y_{i,n,m}^{(r)}|\sigma (\xi _{i,j-r})\right) |%
\mathcal{A}_{i,n,j}^{(r)}\right] \\
&=&\mathbb{E}\left[ \mathbb{E}\left\{ \mathbb{E}\left( x_{i,n,m}|\xi
_{i,m-r}\right) -\mathbb{E}\left( x_{i,n,m}|\xi _{i,m-r-1}\right) |\sigma
(\xi _{i,j-r})\right\} |\mathcal{A}_{i,n,j}^{(r)}\right] \\
&=&0\text{ }\forall m\geq j+1
\end{eqnarray*}%
The second inequality uses Cauchy-Schwartz and Lyapunov inequalities. Now
use (\ref{phys_mart}) and repeat the argument in Case 1 to complete the
proof.\medskip \newline
\textbf{Claim (b).} Recall $\Theta _{i}^{(p)}$ $:=$ $\lim \sup_{n\rightarrow
\infty }\max_{1\leq t\leq n}\Theta _{i,n,t}^{(p)}$ and $\gamma _{i}(\alpha )$
$:=$ $\lim \sup_{p\rightarrow \infty }p^{1/2-1/\alpha }\Theta _{i}^{(p)}$,
and by assumption $(\Theta _{i}^{(p)},\gamma _{i}(\alpha ))$ $\in $ $%
(0,\infty )$ uniformly in $i$\ for some $1$ $<$ $\alpha $ $\leq $ $2$.
Define $\bar{\gamma}$ $:=$ $\max_{i\in \mathbb{N}}\gamma _{i}(\alpha )$ $>$ $%
0$ and $\bar{\lambda}_{0}$ $:=$ $(e\alpha \bar{\gamma}^{\alpha
})^{-1}2^{-\alpha /2}$.

By Stirling's formula and $\bar{\gamma}$ $<$ $\infty $, for any $0$ $<$ $%
\lambda $ $\leq $ $\bar{\lambda}_{0}$ \citep[proof of Theorem 2.(ii)]{Wu2005}%
\begin{eqnarray*}
\limsup_{p\rightarrow \infty }\frac{\lambda \left\{ \sqrt{2\alpha p}%
\max_{i\in \mathbb{N}}\Theta _{i}^{(\alpha p)}\right\} ^{\alpha }}{\left(
p!\right) ^{1/p}} &=&\limsup_{p\rightarrow \infty }\frac{\lambda \left\{ 
\sqrt{2\alpha p}\max_{i\in \mathbb{N}}\Theta _{i}^{(\alpha p)}\right\}
^{\alpha }}{p/e} \\
&=&\limsup_{p\rightarrow \infty }\lambda e\alpha \left\{ \sqrt{2}\left(
\alpha p\right) ^{1/2-1/\alpha }\max_{i\in \mathbb{N}}\Theta _{i}^{(\alpha
p)}\right\} ^{\alpha } \\
&=&\lambda e\alpha 2^{\alpha /2}\bar{\gamma}^{\alpha }<1.
\end{eqnarray*}%
Thus from ($a$) and uniform boundedness%
\begin{equation*}
\max_{i\in \mathbb{N}}\sum_{p=[2/\alpha ]+1}^{\infty }\frac{\mathbb{E}\left(
\lambda \max_{1\leq l\leq n}\left\vert \mathcal{Z}_{i,l}\right\vert ^{\alpha
}\right) ^{p}}{p!}\leq \sum_{p=[2/\alpha ]+1}^{\infty }\frac{\lambda
^{p}\left( \sqrt{2\alpha p}\frac{\alpha p}{\alpha p-1}\max_{i\in \mathbb{N}%
}\Theta _{i}^{(\alpha p)}\right) ^{\alpha p}}{p!}=O(1).
\end{equation*}%
Hence by the Maclaurin series $\max_{i\in \mathbb{N}}\mathbb{E}\exp \left\{
\lambda \max_{1\leq l\leq n}\left\vert \mathcal{Z}_{i,l}\right\vert ^{\alpha
}\right\} $ $<$ $\infty $. The proof now mimics 
\citet[proof of Theorem
2(ii)]{Wu2005} by choosing any $\lambda $ $\in $ $(0,\bar{\lambda}_{0}).$ $%
\mathcal{QED}$.\bigskip \newline
\textbf{Proof of Theorem \ref{thm:max_LLN_phys_dep}.}\medskip \newline
\textbf{Claim (a).} Lemma \ref{lm:phys_dep_Lp}.a and (\ref{Mn2}) yield for $p
$ $>$ $1$, and some $\mathcal{B}_{p}$ $\in $ $(0,\infty )$,%
\begin{equation*}
\mathbb{E}\mathcal{M}_{n}\leq k_{n}^{1/p}\max_{1\leq i\leq k_{n}}\left\Vert 
\frac{1}{n}\sum_{t=1}^{n}x_{i,n,t}\right\Vert _{p}\leq \mathcal{B}%
_{p}k_{n}^{1/p}\frac{1}{n^{1-1/p^{\prime }}}\max_{1\leq i\leq k_{n},1\leq
t\leq n}\Theta _{i,n,t}^{(p)}.
\end{equation*}%
Therefore $\sqrt{n}\mathcal{M}_{n}$ $=$ $O_{p}(k_{n}^{1/p}n^{1/p^{\prime
}-1/2}\max_{1\leq i\leq k_{n},1\leq t\leq n}\Theta _{i,n,t}^{(p)})$. Thus $%
\mathcal{M}_{n}$ $\overset{p}{\rightarrow }$ $0$ as claimed if $k_{n}$ $=$ $%
o(n^{p(1-1/p^{\prime })}/\max_{1\leq i\leq k_{n},1\leq t\leq n}\{\Theta
_{i,n,t}^{(p)}\}^{p})$\textbf{.}\medskip \newline
\textbf{Claim (b).} Use Lemma \ref{lm:phys_dep_Lp}.b with $\mathcal{C}$ $=$ $%
1$ (to reduce notation) together with (\ref{E1/nx}) and (\ref{E1nXIntP}).
First, for some $1$ $<$ $\alpha $ $\leq $ $2$ and any $\lambda $ $>$ $0$,
and by a change of variables $v$ $=$ $(\ln \left( u\right) )^{\alpha }$,%
\begin{eqnarray*}
\mathbb{E}\exp \left\{ \lambda \left\vert \frac{1}{n}\sum_{t=1}^{n}x_{i,n,t}%
\right\vert \right\}  &\leq &e+\int_{e}^{\infty }\mathbb{P}\left( \left\vert 
\frac{1}{n}\sum_{t=1}^{n}x_{i,n,t}\right\vert >\frac{1}{\lambda }\ln \left(
u\right) \right) du \\
&\leq &e+\int_{e}^{\infty }\exp \left\{ -\mathcal{K}\left( n^{1/2}\frac{1}{%
\lambda }\ln \left( u\right) \right) ^{\alpha }\right\} du \\
&=&e+\frac{1}{\alpha }\int_{1}^{\infty }\frac{1}{v^{1-1/\alpha }}\exp
\left\{ v^{1/\alpha }-\mathcal{K}\frac{n^{\alpha /2}}{\lambda ^{\alpha }}%
v\right\} dv \\
&\leq &e+\int_{1}^{\infty }\exp \left\{ v-\mathcal{K}\frac{n^{\alpha /2}}{%
\lambda ^{\alpha }}v\right\} dv.
\end{eqnarray*}%
Notice the last inequality uses $(a,v)$ $\geq $ $1$. Hence for any $\lambda $
$<$ $\mathcal{K}^{1/\alpha }\sqrt{n}$,%
\begin{eqnarray}
\max_{1\leq i\leq k_{n}}\mathbb{E}\exp \left\{ \lambda \left\vert \frac{1}{n}%
\sum_{t=1}^{n}x_{i,n,t}\right\vert \right\}  &\leq &e+\int_{1}^{\infty }\exp
\left\{ -\left( \mathcal{K}n^{\alpha /2}/\lambda ^{\alpha }-1\right)
v\right\} dv  \notag \\
&\leq &e+\frac{\exp \left\{ -\left( \mathcal{K}n^{\alpha /2}/\lambda
^{\alpha }-1\right) \right\} }{\mathcal{K}n^{\alpha /2}/\lambda ^{\alpha }-1}
\notag \\
&\leq &e+\frac{1}{\mathcal{K}n^{\alpha /2}/\lambda ^{\alpha }-1}.
\label{Ee_phys1}
\end{eqnarray}%
Now use (\ref{E1/nx}) and (\ref{Ee_phys1}) to deduce for $\lambda $ $=$ $\xi 
\sqrt{n}$ and any $\xi $ $\in $ $(0,\mathcal{K}^{1/\alpha })$, and $\ln
(k_{n})$ $=$ $o(\sqrt{n})$,%
\begin{eqnarray*}
\mathbb{E}\max_{1\leq i\leq k_{n}}\left\vert \frac{1}{n}%
\sum_{t=1}^{n}x_{i,n,t}\right\vert  &\leq &\frac{1}{\lambda }\ln \left( k_{n}%
\left[ e+\frac{1}{\mathcal{K}n^{\alpha /2}/\lambda ^{\alpha }-1}\right]
\right)  \\
&=&\frac{\ln (k_{n})}{\xi \sqrt{n}}+\frac{1}{\xi \sqrt{n}}\ln \left( e+\frac{%
1}{\mathcal{K}\xi ^{-\alpha }-1}\right) \rightarrow 0.
\end{eqnarray*}

Finally, set $\lambda $ $=$ $\xi \sqrt{n}$ for any $\xi $ $\in $ $(0,%
\mathcal{K}^{1/\alpha })$ to yield 
\begin{eqnarray*}
\mathbb{E}\max_{1\leq i\leq k_{n}}\left\vert \frac{1}{n}%
\sum_{t=1}^{n}x_{i,n,t}\right\vert  &\leq &\frac{\ln (k_{n})}{\xi \sqrt{n}}+%
\frac{1}{\xi \sqrt{n}}\ln \left( e+\frac{1}{\mathcal{K}\left( \sqrt{n}/\left[
\xi \sqrt{n}\right] \right) ^{\alpha }-1}\right)  \\
&=&\frac{\ln (k_{n})}{\xi \sqrt{n}}+\frac{1}{\xi \sqrt{n}}\ln \left( e+\frac{%
1}{\mathcal{K}/\xi ^{\alpha }-1}\right) =\frac{\ln (k_{n})}{\xi \sqrt{n}}%
+O\left( \frac{1}{\sqrt{n}}\right) ,
\end{eqnarray*}%
hence $\max_{1\leq i\leq k_{n}}|1/\sqrt{n}\sum_{t=1}^{n}x_{i,n,t}|$ $=$ $%
O_{p}(\ln \left( k_{n}\right) )$ by Markov's inequality. $\mathcal{QED}$%
.\medskip \newline
\textbf{Proof of Theorem \ref{thm:max_WLLN_mix_units}.} Under $\alpha $%
-mixing $\lim \sup_{n\rightarrow \infty }\alpha _{n}(m)$ $=$ $O(m^{-\lambda
-\iota })$, $\lambda $ $>$ $qp/(q$ $-$ $p)$ and $q$ $>$ $p,$ it follows $%
x_{i,t}$ is an $\mathcal{L}_{q}$-bounded $\mathcal{L}_{p}$-mixingale for
each $i$, $1$ $\leq $ $p$ $\leq $ $q$, with size $\lambda (1/p$ $-$ $1/q)$ %
\citep[Lemma 1.6]{McLeish1975}. Thus $x_{i,t}$ is $\mathcal{L}_{p}$-physical
dependent given $\lambda $ $>$ $qp/(q$ $-$ $p)$ for each $i$\ 
\citep[Theorem
2.1]{Hill2025_mixg}. Moreover, by measurability $\sqrt{n}\mathcal{S}_{i,n}$
is mixing with coefficients $\lim \sup_{n\rightarrow \infty }\tilde{\alpha}%
_{n}(m)$ $=$ $O(m^{-2-\iota })$. Hence $\sqrt{n}\mathcal{S}_{i,n}$ satisfies %
\citet{Leadbetter1974,Leadbetter1983}'s $\mathcal{D}(u_{n})$ property with 
\begin{equation*}
u_{n}:=u/a_{n}+b_{n},
\end{equation*}%
for all $u$ $\in $ $\mathbb{R}$, and some $a_{n}$ $>$ $0$ and $b_{n}$ $\in $ 
$\mathbb{R}$. Furthermore, \citet{Leadbetter1974,Leadbetter1983}'s $\mathcal{%
D}^{\prime }(u_{n})$ property also holds since for any $l$ $>$ $0$%
\begin{eqnarray*}
&&k_{n}\sum_{m=2}^{k_{n}}\mathbb{P}\left( \sqrt{n}\mathcal{S}%
_{i,n}>u_{k_{n}l},\sqrt{n}\mathcal{S}_{i+m,n}>u_{k_{n}l}\right)  \\
&&\text{ \ }=k_{n}\sum_{m=2}^{k_{n}}\left\{ \mathbb{P}\left( \sqrt{n}%
\mathcal{S}_{i,n}>u_{k_{n}l},\sqrt{n}\mathcal{S}_{i+m,n}>u_{k_{n}l}\right) -%
\mathbb{P}\left( \sqrt{n}\mathcal{S}_{i,n}>u_{k_{n}l}\right) \mathbb{P}%
\left( \sqrt{n}\mathcal{S}_{i+m,n}>u_{k_{n}l}\right) \right\}  \\
&&\text{ \ \ \ \ \ \ \ \ \ \ \ \ \ \ \ }+k_{n}\sum_{m=2}^{k_{n}}\mathbb{P}%
\left( \sqrt{n}\mathcal{S}_{i,n}>u_{k_{n}l}\right) \mathbb{P}\left( \sqrt{n}%
\mathcal{S}_{i+m,n}>u_{k_{n}l}\right)  \\
&&\text{ \ }\leq k_{n}\sum_{m=1}^{k_{n}-1}\tilde{\alpha}_{n}(m)+\frac{1}{%
l^{2}}\times lk_{n}\mathbb{P}\left( \sqrt{n}\mathcal{S}_{i,n}>u_{k_{n}l}%
\right) \times \frac{1}{k_{n}}\sum_{m=2}^{k_{n}}lk_{n}\mathbb{P}\left( \sqrt{%
n}\mathcal{S}_{i+m,n}>u_{k_{n}l}\right)  \\
&&\text{ \ }\leq Kk_{n}\sum_{m=1}^{k_{n}-1}m^{-2-\iota }+\tau ^{2}\frac{1}{%
l^{2}}\left( 1+o\left( 1\right) \right)  \\
&&\text{ \ }\simeq Kk_{n}\frac{1}{k_{n}^{1+\iota }}+\tau ^{2}\frac{1}{l^{2}}%
\left( 1+o\left( 1\right) \right) =o\left( 1/l\right) \text{, cf. \citet[eq.
(3.2)]{Leadbetter1974}}.
\end{eqnarray*}%
The second and third inequalities use $\max_{1\leq i\leq k_{n}}k_{n}\mathbb{P%
}(\sqrt{n}\mathcal{S}_{i,n}$ $>$ $u_{k_{n}})$ $\rightarrow $ $\tau $. The
first uses the $\alpha $-mixing coefficient construction implication 
\begin{equation*}
\left\vert \mathbb{P}\left( \sqrt{n}\mathcal{S}_{i,n}>u_{k_{n}l},\sqrt{n}%
\mathcal{S}_{i+m,n}>u_{k_{n}l}\right) -\mathbb{P}\left( \sqrt{n}\mathcal{S}%
_{m,n}>u_{k_{n}l}\right) \times \mathbb{P}\left( \sqrt{n}\mathcal{S}%
_{i+m,n}>u_{k_{n}l}\right) \right\vert \leq \tilde{\alpha}_{n}(m).
\end{equation*}%
The conditions of Theorem 1.2 in \cite{Leadbetter1983} therefore hold for $%
\max_{1\leq i\leq k_{n}}|\sqrt{n}\mathcal{S}_{i,n}|$: 
\begin{equation}
\mathbb{P}\left( a_{k_{n}}\left\{ \max_{1\leq i\leq k_{n}}\left\vert \sqrt{n}%
\mathcal{S}_{i,n}\right\vert -b_{k_{n}}\right\} >u\right) =\mathbb{P}\left(
\max_{1\leq i\leq k_{n}}\left\vert \sqrt{n}\mathcal{S}_{i,n}\right\vert \leq
u_{k_{n}}\right) \rightarrow \exp \{-\tau \}\text{ }\forall u\in \mathbb{R}.
\label{Lead1983}
\end{equation}%
Clearly the centering sequence $\{b_{n}\}$ must satisfy $b_{n}$ $>$ $0$
given $|\sqrt{n}\mathcal{S}_{i,n}|$ $\geq $ $0$. Therefore $\forall u$ $>$ $0
$ and $\sqrt{n}a_{k_{n}}$ $\rightarrow $ $\infty $ and $b_{k_{n}}/\sqrt{n}$ $%
\rightarrow $ $0$,%
\begin{equation*}
\mathbb{P}\left( \mathcal{M}_{n}>u\right) =\mathbb{P}\left( a_{k_{n}}\left\{
\max_{1\leq i\leq k_{n}}\left\vert \sqrt{n}\mathcal{S}_{i,n}\right\vert
-b_{k_{n}}\right\} >\sqrt{n}a_{k_{n}}\left\{ u-\frac{b_{k_{n}}}{\sqrt{n}}%
\right\} \right) \rightarrow 0.
\end{equation*}%
Simply note that (\ref{Lead1983}) yields $a_{k_{n}}\{\max_{1\leq i\leq
k_{n}}|\sqrt{n}\mathcal{S}_{i,n}|$ $-$ $b_{k_{n}}\}$ $=$ $O_{p}(1)$, hence%
\begin{equation*}
\frac{a_{k_{n}}\left\{ \max_{1\leq i\leq k_{n}}\left\vert \sqrt{n}\mathcal{S}%
_{i,n}\right\vert -b_{k_{n}}\right\} }{\sqrt{n}a_{k_{n}}}+\frac{b_{k_{n}}}{%
\sqrt{n}}=O_{p}\left( \frac{1}{\sqrt{n}a_{k_{n}}}+\frac{b_{k_{n}}}{\sqrt{n}}%
\right) =o_{p}(1).
\end{equation*}%
This suffices to prove $\mathcal{M}_{n}$ $\overset{p}{\rightarrow }$ $0$ if $%
\sqrt{n}a_{k_{n}}$ $\rightarrow $ $\infty $ and $b_{k_{n}}/\sqrt{n}$ $%
\rightarrow $ $0$ as required. $\mathcal{QED}$.

\newpage \clearpage

\begin{sidewaystable}
\centering
\caption{Summary of Results in Order of Appearance} \label{tlb:sum} 
\begin{tabular}{llllll}
\textbf{Grouping} & \textbf{Dependence} & \textbf{Tails} & \textbf{Coordinate%
}$^{\text{1}}$ & \textbf{Rate $\boldsymbol{k_{n}}$}$^{\text{2}}$ & \textbf{%
Reference} \\ \hline\hline
$(i)$ & \multicolumn{1}{|l}{$\tau $-mixing} & sub-exp & unrestricted & $\ln
(k_{n})=o(n)$ & Theorem \ref{thm:max_LLN_t_mix} \\ \hline
$(ii)$ & \multicolumn{1}{|l}{$a.$ $\mathcal{L}_{p}$-phys. dep. ($p>1$)} & 
unrestricted & unrestricted & $k_{n}$ $=$ $o\left( \frac{n^{p(1-1/p^{\prime
})}}{\Theta _{n}^{(p)}}\right) $ & Theorem \ref{thm:max_LLN_phys_dep}.a \\ 
& \multicolumn{1}{|l}{$\ \ \ \,\mathcal{L}_{p}$-phys. dep. ($p>1$)} & sub-exp
& unrestricted & $\ln (k_{n})=O(\sqrt{n})$ & Theorem \ref%
{thm:max_LLN_phys_dep}.b \\ 
& \multicolumn{1}{|l}{$b.$ $\mathcal{L}_{p}$-phys. dep. ($p>1$)} & 
unrestricted & martingale & $k_{n}$ arbitrary & Theorem \ref%
{thm:max_WLLN_mart} \\ 
& \multicolumn{1}{|l}{$c.$ $\mathcal{L}_{p}$-phys. dep. ($p>1$)} & 
unrestricted & nearly mart. & $\ln (k_{n})=o(n)$ & Theorem \ref%
{thm:max_WLLN_nearmart} \\ 
& \multicolumn{1}{|l}{$d.$ $\mathcal{L}_{p}$-phys. dep. ($p>1$)} & 
sub-Gaussian & mixing & $\ln (k_{n})=o(n)$ & Theorem \ref%
{thm:max_WLLN_mix_units} \\ 
& \multicolumn{1}{|l}{$\ \ \ \,\mathcal{L}_{p}$-phys. dep. ($p>1$)} & stable
domain & mixing & $\ln (k_{n})=O(\sqrt{n})$ & Theorem \ref%
{thm:max_WLLN_mix_units} \\ 
& \multicolumn{1}{|l}{$\ \ \ \,\mathcal{L}_{p}$-phys. dep. ($p>1$)} & stable
domain & mixing & $k_{n}=O\left( \frac{n^{a}}{h(n)}\right) $ & Theorem \ref%
{thm:max_WLLN_mix_units} \\ \hline
$(iii)$ & \multicolumn{1}{|l}{independence} & sub-exp & unrestricted & $\ln
(k_{n})=O\left(\frac{\sqrt{n}}{h(n)}\right)$ & Theorem B.2$^{\text{3}}$ \\ 
& \multicolumn{1}{|l}{independence} & unrestricted & unrestricted & $\ln
(k_{n})=o\left( \frac{n}{g(n)}\right) $ & Theorem B.4%
\end{tabular}
\bigskip

\centering\caption{Summary of Results Ranked by $k_{n}$} \label{tlb:rank} 
\begin{tabular}{lllll}
\textbf{Dependence} & \textbf{Tails} & \textbf{Coordinate}$^{\text{1}}$ & 
\textbf{Rate $\boldsymbol{k_{n}}$}$^{\text{2}}$ & \textbf{Reference} \\ 
\hline\hline
$\mathcal{L}_{p}$-phys. dep. ($p>1$) & unrestricted & martingale & $k_{n}$
arbitrary & Theorem \ref{thm:max_WLLN_mart} \\ \hline
$\tau $-mixing & sub-exp & unrestricted & $\ln (k_{n})=o(n)$ & Theorem \ref%
{thm:max_LLN_t_mix} \\ \hline
$\mathcal{L}_{p}$-phys. dep. ($p>1$) & unrestricted & nearly mart. & $\ln
(k_{n})=o(n)$ & Theorem \ref{thm:max_WLLN_nearmart} \\ \hline
$\mathcal{L}_{p}$-phys. dep. ($p>1$) & sub-Gaussian & mixing & $\ln
(k_{n})=o(n)$ & Theorem \ref{thm:max_WLLN_mix_units} \\ \hline
independence & unrestricted & unrestricted & $\ln (k_{n})=o\left( \frac{n}{%
g(n)}\right) $ & Theorem B.4 \\ \hline
$\mathcal{L}_{p}$-phys. dep. ($p>1$) & sub-exp & unrestricted & $\ln
(k_{n})=O(\sqrt{n})$ & Theorem \ref{thm:max_LLN_phys_dep}.b \\ \hline
$\mathcal{L}_{p}$-phys. dep. ($p>1$) & stable domain & mixing & $\ln
(k_{n})=O(\sqrt{n})$ & Theorem \ref{thm:max_WLLN_mix_units} \\ \hline
Independence & sub-exp & unrestricted & $\ln (k_{n})=O\left( \frac{\sqrt{n}}{%
h(n)}\right) $ & Theorem B.2$^{\text{3}}$ \\ \hline
$\mathcal{L}_{p}$-phys. dep. ($p>1$) & stable domain & mixing & $%
k_{n}=O\left( \frac{n^{a}}{h(n)}\right) $ & Theorem \ref%
{thm:max_WLLN_mix_units} \\ \hline
$\mathcal{L}_{p}$-phys. dep. ($p>1$) & unrestricted & unrestricted & $k_{n}$ 
$=$ $o\left( \frac{n^{p(1-1/p^{\prime })}}{\Theta _{n}^{(p)}}\right) $ & 
Theorem \ref{thm:max_LLN_phys_dep}.a%
\end{tabular}
\medskip \newline
\parbox{0.9\linewidth}{
\raggedright
\footnotesize
1. Coordinate restrictions are over $i$; mart. = martingale.\newline
2. Rates are for WLLN's; $\Theta _{n}^{(p)}=\sum_{m=0}^{\infty }\theta
_{i,n,t}^{(p)}(m)$ with \noindent $\mathcal{L}_{p}$-physical dependence
measure $\theta _{i,n,t}^{(p)}(m)$; $a$ $>$ $0$ and $h(n)$ is slowly
varying; $g(n)$ $=$ $o(n)$ depends on tail decay.\newline
3. See \cite{sm_max_LLN} for Theorems B.2 and B.4.}
\end{sidewaystable}

\newpage \clearpage
{\small \setstretch{.75} 
\bibliographystyle{cas-model2-names.bst}
\bibliography{refs_max_LLN}

\begin{thebibliography}{90}
\expandafter\ifx\csname natexlab\endcsname\relax\def\natexlab#1{#1}\fi
\providecommand{\url}[1]{\texttt{#1}}
\providecommand{\href}[2]{#2}
\providecommand{\path}[1]{#1}
\providecommand{\DOIprefix}{doi:}
\providecommand{\ArXivprefix}{arXiv:}
\providecommand{\URLprefix}{URL: }
\providecommand{\Pubmedprefix}{pmid:}
\providecommand{\doi}[1]{\href{http://dx.doi.org/#1}{\path{#1}}}
\providecommand{\Pubmed}[1]{\href{pmid:#1}{\path{#1}}}
\providecommand{\bibinfo}[2]{#2}
\ifx\xfnm\relax \def\xfnm[#1]{\unskip,\space#1}\fi
\bibitem[{Adamek et~al.(2023)Adamek, Smeekes and Wilms}]{Adamek_et_al2023}
\bibinfo{author}{Adamek, R.}, \bibinfo{author}{Smeekes, S.},
  \bibinfo{author}{Wilms, I.}, \bibinfo{year}{2023}.
\newblock \bibinfo{title}{Lasso inference for high-dimensional time series}.
\newblock \bibinfo{journal}{J. Econometrics} \bibinfo{volume}{235},
  \bibinfo{pages}{1114--1143}.
\bibitem[{Andrews(1984)}]{Andrews1984}
\bibinfo{author}{Andrews, D.W.K.}, \bibinfo{year}{1984}.
\newblock \bibinfo{title}{Nonstrong mixing autoregressive processes}.
\newblock \bibinfo{journal}{J. Appl. Probab.} \bibinfo{volume}{21},
  \bibinfo{pages}{930--934}.
\bibitem[{Andrews(1987)}]{Andrews1987}
\bibinfo{author}{Andrews, D.W.K.}, \bibinfo{year}{1987}.
\newblock \bibinfo{title}{Consistency in nonlinear econometric models: A
  generic uniform law of large numbers}.
\newblock \bibinfo{journal}{Econometrica} \bibinfo{volume}{55},
  \bibinfo{pages}{1465--1471}.
\bibitem[{Andrews(1988)}]{Andrews1988}
\bibinfo{author}{Andrews, D.W.K.}, \bibinfo{year}{1988}.
\newblock \bibinfo{title}{Laws of large numbers for dependent non-identically
  distributed random variables}.
\newblock \bibinfo{journal}{Econometric Theory} \bibinfo{volume}{4},
  \bibinfo{pages}{458--467}.
\bibitem[{Babii et~al.(2019)Babii, Ghysels and Striaukas}]{Babii_etal2019}
\bibinfo{author}{Babii, A.}, \bibinfo{author}{Ghysels, E.},
  \bibinfo{author}{Striaukas, J.}, \bibinfo{year}{2019}.
\newblock \bibinfo{title}{Estimation and hac-based inference for machine
  learning time series regressions}.
\newblock \bibinfo{note}{Dept. of Economics, University of North Carolina -
  Chapel Hill}.
\bibitem[{Belloni et~al.(2014)Belloni, Chernozhukov and
  Hansen}]{BelloniChernozhukovHansen2014}
\bibinfo{author}{Belloni, A.}, \bibinfo{author}{Chernozhukov, V.},
  \bibinfo{author}{Hansen, C.}, \bibinfo{year}{2014}.
\newblock \bibinfo{title}{High-dimensional methods and inference on structural
  and treatment effects}.
\newblock \bibinfo{journal}{J. Econom. Perspect.} \bibinfo{volume}{28},
  \bibinfo{pages}{29--50}.
\bibitem[{Bentkus(2008)}]{Bentkus2008}
\bibinfo{author}{Bentkus, V.}, \bibinfo{year}{2008}.
\newblock \bibinfo{title}{An extension of the {H}oeffding inequality to
  unbounded random variables}.
\newblock \bibinfo{journal}{Lith. Math. J.} \bibinfo{volume}{48},
  \bibinfo{pages}{137--157}.
\bibitem[{Berman(1964)}]{Berman1964}
\bibinfo{author}{Berman, S.M.}, \bibinfo{year}{1964}.
\newblock \bibinfo{title}{Limit theorems for the maximum term in stationary
  sequences}.
\newblock \bibinfo{journal}{Ann. Math. Stat.} \bibinfo{volume}{35},
  \bibinfo{pages}{502--516}.
\bibitem[{Bosq(1993)}]{Bosq1993}
\bibinfo{author}{Bosq, D.}, \bibinfo{year}{1993}.
\newblock \bibinfo{title}{Bernstein-type large deviations inequalities for
  partial sums of strong mixing processes}.
\newblock \bibinfo{journal}{Statistics} \bibinfo{volume}{24},
  \bibinfo{pages}{59--70}.
\bibitem[{B{\"u}hlmann and van~de Geer(2011)}]{BuhlmannVanDeGeer2011}
\bibinfo{author}{B{\"u}hlmann, P.}, \bibinfo{author}{van~de Geer, S.},
  \bibinfo{year}{2011}.
\newblock \bibinfo{title}{Statistics for High-Dimensional Data}.
\newblock \bibinfo{publisher}{Springer}, \bibinfo{address}{Berlin}.
\bibitem[{Burkholder(1973)}]{Burkholder1973}
\bibinfo{author}{Burkholder, D.L.}, \bibinfo{year}{1973}.
\newblock \bibinfo{title}{Distribution function inequalities for martingales}.
\newblock \bibinfo{journal}{Ann. Probab.} \bibinfo{volume}{1},
  \bibinfo{pages}{19--42}.
\bibitem[{Cattaneo et~al.(2018)Cattaneo, Jansson and
  Newey}]{CattaneoJanssonNewey2018}
\bibinfo{author}{Cattaneo, M.D.}, \bibinfo{author}{Jansson, M.},
  \bibinfo{author}{Newey, W.}, \bibinfo{year}{2018}.
\newblock \bibinfo{title}{Inference in linear regression models with many
  covariates and heteroscedasticity}.
\newblock \bibinfo{journal}{J. Amer. Statist. Assoc.} \bibinfo{volume}{113},
  \bibinfo{pages}{1350--1361}.
\bibitem[{Chang et~al.(2024)Chang, Chen and Wu}]{ChangChenWu2024}
\bibinfo{author}{Chang, J.}, \bibinfo{author}{Chen, X.}, \bibinfo{author}{Wu,
  M.}, \bibinfo{year}{2024}.
\newblock \bibinfo{title}{Central limit theorems for high dimensional dependent
  data}.
\newblock \bibinfo{journal}{Bernoulli} \bibinfo{volume}{30},
  \bibinfo{pages}{712--742}.
\bibitem[{Chang et~al.(2023)Chang, Jiang and Shao}]{ChangJiangShao2023}
\bibinfo{author}{Chang, J.}, \bibinfo{author}{Jiang, Q.},
  \bibinfo{author}{Shao, X.}, \bibinfo{year}{2023}.
\newblock \bibinfo{title}{Testing the martingale difference hypothesis in high
  dimension}.
\newblock \bibinfo{journal}{J. Econometrics} \bibinfo{volume}{235},
  \bibinfo{pages}{972--1000}.
\bibitem[{Chazottes and Gouezel(2012)}]{ChazottesGouezel2012}
\bibinfo{author}{Chazottes, J.R.}, \bibinfo{author}{Gouezel, S.},
  \bibinfo{year}{2012}.
\newblock \bibinfo{title}{Optimal concentration inequalities for dynamical
  systems}.
\newblock \bibinfo{journal}{Comm. Math. Phys.} \bibinfo{volume}{316},
  \bibinfo{pages}{843--889}.
\bibitem[{Chernick(1981)}]{Chernick1981}
\bibinfo{author}{Chernick, M.R.}, \bibinfo{year}{1981}.
\newblock \bibinfo{title}{A limit theorem for the maximum of autoregessive
  processes with uniform marginal distribution}.
\newblock \bibinfo{journal}{Ann. Probab.} \bibinfo{volume}{9},
  \bibinfo{pages}{145--149}.
\bibitem[{Chernozhukov et~al.(2013)Chernozhukov, Chetverikov and
  Kato}]{Chernozhukov_etal2013}
\bibinfo{author}{Chernozhukov, V.}, \bibinfo{author}{Chetverikov, D.},
  \bibinfo{author}{Kato, K.}, \bibinfo{year}{2013}.
\newblock \bibinfo{title}{Gaussian approximations and multiplier bootstrap for
  maxima of sums of high-dimensional random vectors}.
\newblock \bibinfo{journal}{Ann. Statist.} \bibinfo{volume}{41},
  \bibinfo{pages}{2786--2819}.
\bibitem[{Cohen(1982)}]{Cohen1982}
\bibinfo{author}{Cohen, J.P.}, \bibinfo{year}{1982}.
\newblock \bibinfo{title}{The penultimate form of approximation to normal
  extremes}.
\newblock \bibinfo{journal}{Adv. Appl. Probab.} \bibinfo{volume}{14},
  \bibinfo{pages}{324--339}.
\bibitem[{Collet et~al.(2002)Collet, Martinez and Schmitt}]{Collet_etal2002}
\bibinfo{author}{Collet, P.}, \bibinfo{author}{Martinez, S.},
  \bibinfo{author}{Schmitt, B.}, \bibinfo{year}{2002}.
\newblock \bibinfo{title}{Exponential inequalities for dynamical measures of
  expanding maps of the interval}.
\newblock \bibinfo{journal}{Probab. Theory Rel.} \bibinfo{volume}{123},
  \bibinfo{pages}{301--322}.
\bibitem[{Davidson(1994)}]{Davidson1994}
\bibinfo{author}{Davidson, J.}, \bibinfo{year}{1994}.
\newblock \bibinfo{title}{Stochastic Limit Theory}.
\newblock \bibinfo{publisher}{Oxford University Press},
  \bibinfo{address}{Oxford, U. K.}
\bibitem[{Dedecker and Doukhan(2003)}]{DedeckerDoukhan2003}
\bibinfo{author}{Dedecker, J.}, \bibinfo{author}{Doukhan, P.},
  \bibinfo{year}{2003}.
\newblock \bibinfo{title}{A new covariance inequality and applications}.
\newblock \bibinfo{journal}{Stochastic Process. Appl.} \bibinfo{volume}{106},
  \bibinfo{pages}{63--80}.
\bibitem[{Dedecker et~al.(2007)Dedecker, Doukhan, Lang, Leon, Louhichi and
  Prieur}]{Dedecker_etal2007}
\bibinfo{author}{Dedecker, J.}, \bibinfo{author}{Doukhan, P.},
  \bibinfo{author}{Lang, G.}, \bibinfo{author}{Leon, J.R.},
  \bibinfo{author}{Louhichi, S.}, \bibinfo{author}{Prieur, C.},
  \bibinfo{year}{2007}.
\newblock \bibinfo{title}{Weak Dependence: With Examples and Applications}.
\newblock \bibinfo{publisher}{Springer}.
\bibitem[{Dedecker and Prieur(2004)}]{Dedecker_Prieur_2004}
\bibinfo{author}{Dedecker, J.}, \bibinfo{author}{Prieur, C.},
  \bibinfo{year}{2004}.
\newblock \bibinfo{title}{Coupling for $\tau$-dependent sequences and
  applications}.
\newblock \bibinfo{journal}{J. Theoret. Probab.} \bibinfo{volume}{17},
  \bibinfo{pages}{861--885}.
\bibitem[{Dedecker and Prieur(2005)}]{DedeckerPrieur2005}
\bibinfo{author}{Dedecker, J.}, \bibinfo{author}{Prieur, C.},
  \bibinfo{year}{2005}.
\newblock \bibinfo{title}{New dependence coefficients. examples and
  applications to statistics}.
\newblock \bibinfo{journal}{Probab. Theory Rel.} \bibinfo{volume}{132},
  \bibinfo{pages}{203--235}.
\bibitem[{Dezeure et~al.(2017)Dezeure, B{\"u}hlmann and
  Zhang}]{Dezeure_etal_2017}
\bibinfo{author}{Dezeure, R.}, \bibinfo{author}{B{\"u}hlmann, P.},
  \bibinfo{author}{Zhang, C.H.}, \bibinfo{year}{2017}.
\newblock \bibinfo{title}{High-dimensional simultaneous inference with the
  bootstrap}.
\newblock \bibinfo{journal}{Test} \bibinfo{volume}{26},
  \bibinfo{pages}{685--719}.
\bibitem[{Diaconis and Freedman(1999)}]{DiaconisFreedman1999}
\bibinfo{author}{Diaconis, P.}, \bibinfo{author}{Freedman, D.},
  \bibinfo{year}{1999}.
\newblock \bibinfo{title}{Iterated random functions}.
\newblock \bibinfo{journal}{SIAM Rev.} \bibinfo{volume}{41},
  \bibinfo{pages}{45--76}.
\bibitem[{Dumbgen et~al.(2010)Dumbgen, van~de Geer, Veraar and
  Wellner}]{Dumbgen_etal_2010}
\bibinfo{author}{Dumbgen, L.}, \bibinfo{author}{van~de Geer, S.},
  \bibinfo{author}{Veraar, M.C.}, \bibinfo{author}{Wellner, J.A.},
  \bibinfo{year}{2010}.
\newblock \bibinfo{title}{Nemirovski's inequalities revisited}.
\newblock \bibinfo{journal}{Amer. Math. Monthly} \bibinfo{volume}{117},
  \bibinfo{pages}{138--160}.
\bibitem[{Fan and Li(2006)}]{FanLi2006}
\bibinfo{author}{Fan, J.}, \bibinfo{author}{Li, R.}, \bibinfo{year}{2006}.
\newblock \bibinfo{title}{Statistical challenges with high dimensionality:
  Feature selection in knowledge discovery}, in: \bibinfo{editor}{Sanz-Sole,
  M.}, \bibinfo{editor}{Soria, J.}, \bibinfo{editor}{Varona, J.L.},
  \bibinfo{editor}{Verdera, J.} (Eds.), \bibinfo{booktitle}{Proceedings of the
  International Congress of Mathematicians}, \bibinfo{organization}{European
  Mathematical Society}, \bibinfo{address}{Zurich}. pp.
  \bibinfo{pages}{595--622}.
\bibitem[{Fan et~al.(2011)Fan, Lv and Qi}]{FanLvQi2011}
\bibinfo{author}{Fan, J.}, \bibinfo{author}{Lv, J.}, \bibinfo{author}{Qi, .L.},
  \bibinfo{year}{2011}.
\newblock \bibinfo{title}{Sparse high-dimensional models in economics}.
\newblock \bibinfo{journal}{Annu. Rev. Economics} \bibinfo{volume}{3},
  \bibinfo{pages}{291--317}.
\bibitem[{Genovese et~al.(2012)Genovese, Jin, Wasserman and
  Yao}]{Genovese_etal_2012}
\bibinfo{author}{Genovese, C.R.}, \bibinfo{author}{Jin, J.},
  \bibinfo{author}{Wasserman, L.}, \bibinfo{author}{Yao, Z.},
  \bibinfo{year}{2012}.
\newblock \bibinfo{title}{A comparison of the lasso and marginal regression}.
\newblock \bibinfo{journal}{J. Mach. Learn. Res.} \bibinfo{volume}{13},
  \bibinfo{pages}{2107--2143}.
\bibitem[{Gordin(1969)}]{Gordin1969}
\bibinfo{author}{Gordin, M.I.}, \bibinfo{year}{1969}.
\newblock \bibinfo{title}{The central limit theorem for stationary processes}.
\newblock \bibinfo{journal}{Soviet. Math. Dokl.} \bibinfo{volume}{10},
  \bibinfo{pages}{1174--1176}.
\bibitem[{Hall(1979)}]{Hall1979}
\bibinfo{author}{Hall, P.}, \bibinfo{year}{1979}.
\newblock \bibinfo{title}{On the rate of convergence of normal extremes}.
\newblock \bibinfo{journal}{J. Appl. Probab.} \bibinfo{volume}{16},
  \bibinfo{pages}{433--439}.
\bibitem[{Hall and Heyde(1980)}]{HallHeyde1980}
\bibinfo{author}{Hall, P.}, \bibinfo{author}{Heyde, C.C.},
  \bibinfo{year}{1980}.
\newblock \bibinfo{title}{Martingale Limit Theory and Its Application}.
\newblock \bibinfo{publisher}{AcademicPress}, \bibinfo{address}{New York}.
\bibitem[{Hang and Steinwart(2017)}]{HangSteinwart2017}
\bibinfo{author}{Hang, H.}, \bibinfo{author}{Steinwart, I.},
  \bibinfo{year}{2017}.
\newblock \bibinfo{title}{A {B}ernstein-type inequality for some mixing
  processes and dynamical systems with an application to learning}.
\newblock \bibinfo{journal}{Ann. Statist.} \bibinfo{volume}{45},
  \bibinfo{pages}{708--743}.
\bibitem[{Hannan(1973)}]{Hannan1973}
\bibinfo{author}{Hannan, E.J.}, \bibinfo{year}{1973}.
\newblock \bibinfo{title}{Central limit theorems for time series regression}.
\newblock \bibinfo{journal}{Z. Wahrscheinlichkeitstheorie verw. Geb.}
  \bibinfo{volume}{26}, \bibinfo{pages}{157--170}.
\bibitem[{Hannan and Deistler(1988)}]{HannanDiestler1988}
\bibinfo{author}{Hannan, E.J.}, \bibinfo{author}{Deistler, M.},
  \bibinfo{year}{1988}.
\newblock \bibinfo{title}{The Statistical Theory of Linear Systems}.
\newblock \bibinfo{publisher}{Wiley}, \bibinfo{address}{New York}.
\bibitem[{Hansen(1991)}]{Hansen1991}
\bibinfo{author}{Hansen, B.E.}, \bibinfo{year}{1991}.
\newblock \bibinfo{title}{Strong laws for dependent heterogeneous processes}.
\newblock \bibinfo{journal}{Econometric Theory} \bibinfo{volume}{7},
  \bibinfo{pages}{213--221}.
\bibitem[{Herrmann et~al.(2014)Herrmann, Hofert and
  Ne\v{s}lehov\'{a}}]{Herrmann_etal_2014}
\bibinfo{author}{Herrmann, K.}, \bibinfo{author}{Hofert, M.},
  \bibinfo{author}{Ne\v{s}lehov\'{a}, J.G.}, \bibinfo{year}{2014}.
\newblock \bibinfo{title}{Limiting behavior of maxima under dependence}.
\newblock \bibinfo{note}{{h}ttps://arxiv.org/abs/2405.02833}.
\bibitem[{Hill(2025a)}]{Hill2025_mixg}
\bibinfo{author}{Hill, J.B.}, \bibinfo{year}{2025}a.
\newblock \bibinfo{title}{Mixingale and physical dependence equality with
  applications}.
\newblock \bibinfo{journal}{Statist. Probab. Lett.} \bibinfo{volume}{221},
  \bibinfo{pages}{110380}.
\bibitem[{Hill(2025b)}]{Hill_2025_maxtest}
\bibinfo{author}{Hill, J.B.}, \bibinfo{year}{2025}b.
\newblock \bibinfo{title}{Testing many zero restrictions in a high dimensional
  linear regression setting}.
\newblock \bibinfo{journal}{J. Bus. Econom. Statist.} \bibinfo{volume}{43},
  \bibinfo{pages}{55--67}.
\bibitem[{Hill(2026)}]{sm_max_LLN}
\bibinfo{author}{Hill, J.B.}, \bibinfo{year}{2026}.
\newblock \bibinfo{title}{Supplemental material for `\textit{Max-Laws of Large
  Numbers for High Dimensional Arrays with Applications}'}.
\newblock \bibinfo{note}{Dept. of Economics, University of North Carolina -
  Chapel Hill}.
\bibitem[{Hill et~al.(2020)Hill, Ghysels and Motegi}]{HillGhyselsMotegi2020}
\bibinfo{author}{Hill, J.B.}, \bibinfo{author}{Ghysels, E.},
  \bibinfo{author}{Motegi, K.}, \bibinfo{year}{2020}.
\newblock \bibinfo{title}{Testing a large set of zero restrictions in
  regression models, with an application to mixed frequency granger causality}.
\newblock \bibinfo{journal}{J. Econometrics} \bibinfo{volume}{218},
  \bibinfo{pages}{633--654}.
\bibitem[{Hill and Li(2025)}]{HillLi2025}
\bibinfo{author}{Hill, J.B.}, \bibinfo{author}{Li, T.}, \bibinfo{year}{2025}.
\newblock \bibinfo{title}{A bootstrapped test of covariance stationarity based
  on orthonormal transformations}.
\newblock \bibinfo{journal}{Bernoulli} \bibinfo{volume}{31},
  \bibinfo{pages}{1527--1551}.
\bibitem[{Hill and Motegi(2020)}]{HillMotegi2020}
\bibinfo{author}{Hill, J.B.}, \bibinfo{author}{Motegi, K.},
  \bibinfo{year}{2020}.
\newblock \bibinfo{title}{A max-correlation white noise test for weakly
  dependent time series}.
\newblock \bibinfo{journal}{Econometric Theory} \bibinfo{volume}{36},
  \bibinfo{pages}{907--960}.
\bibitem[{Hsing et~al.(1996)Hsing, H{\"u}sler and Reiss}]{HsingHuslerReiss1996}
\bibinfo{author}{Hsing, T.}, \bibinfo{author}{H{\"u}sler},
  \bibinfo{author}{Reiss, R.D.}, \bibinfo{year}{1996}.
\newblock \bibinfo{title}{The extremes of a triangular array of normal random
  variables}.
\newblock \bibinfo{journal}{Ann. Appl. Probab.} \bibinfo{volume}{6},
  \bibinfo{pages}{671--686}.
\bibitem[{Ibragimov(1962)}]{Ibragimov1962}
\bibinfo{author}{Ibragimov, I.A.}, \bibinfo{year}{1962}.
\newblock \bibinfo{title}{Some limit theorems for stationary processes}.
\newblock \bibinfo{journal}{Theory Probab. Appl.} \bibinfo{volume}{7},
  \bibinfo{pages}{349--382}.
\bibitem[{Jin et~al.(2015)Jin, Wang and Wang}]{JinWangWang2015}
\bibinfo{author}{Jin, L.}, \bibinfo{author}{Wang, S.}, \bibinfo{author}{Wang,
  H.}, \bibinfo{year}{2015}.
\newblock \bibinfo{title}{A new non-parametric stationarity test of time series
  in the time domain}.
\newblock \bibinfo{journal}{J. Roy. Stat. Soc. Ser. B} \bibinfo{volume}{77},
  \bibinfo{pages}{893--922}.
\bibitem[{Jirak and K{\"o}stenberger(2024)}]{JirakKostenberger2024}
\bibinfo{author}{Jirak, M.J.}, \bibinfo{author}{K{\"o}stenberger, G.},
  \bibinfo{year}{2024}.
\newblock \bibinfo{title}{Sharp oracle inequalities and universality of the aic
  and fpe}.
\newblock \bibinfo{note}{ArXiv:2406.13513v1}.
\bibitem[{Kallenberg(2021)}]{Kallenberg2021}
\bibinfo{author}{Kallenberg, O.}, \bibinfo{year}{2021}.
\newblock \bibinfo{title}{Foundations of Modern Probability}.
\newblock \bibinfo{edition}{3rd} ed., \bibinfo{publisher}{Springer Nature},
  \bibinfo{address}{Switzerland}.
\bibitem[{Keenan(1997)}]{Keenan1997}
\bibinfo{author}{Keenan, D.M.}, \bibinfo{year}{1997}.
\newblock \bibinfo{title}{A central limit theorem for m(n) autocovariances}.
\newblock \bibinfo{journal}{J. Time Series Anal.} \bibinfo{volume}{18},
  \bibinfo{pages}{61--78}.
\bibitem[{Krebs(2018a)}]{Krebs2018a}
\bibinfo{author}{Krebs, J.T.N.}, \bibinfo{year}{2018}a.
\newblock \bibinfo{title}{A {B}ernstein inequality for exponentially growing
  graphs}.
\newblock \bibinfo{journal}{Comm. Statist. - Theory and Method}
  \bibinfo{volume}{47}, \bibinfo{pages}{5097--5106}.
\bibitem[{Krebs(2018b)}]{Krebs2018b}
\bibinfo{author}{Krebs, J.T.N.}, \bibinfo{year}{2018}b.
\newblock \bibinfo{title}{A large deviation inequality for $\beta$ mixing time
  series and its applications to the functional kernel regression mode}.
\newblock \bibinfo{journal}{Statist. Probab. Lett.} \bibinfo{volume}{113},
  \bibinfo{pages}{50--58}.
\bibitem[{K{\"u}hn and Schilling(2023)}]{KuhnSchilling2023}
\bibinfo{author}{K{\"u}hn, F.}, \bibinfo{author}{Schilling, R.L.},
  \bibinfo{year}{2023}.
\newblock \bibinfo{title}{Maximal inequalities and some applications}.
\newblock \bibinfo{journal}{Probab. Surveys} \bibinfo{volume}{20},
  \bibinfo{pages}{382--485}.
\bibitem[{Leadbetter(1974)}]{Leadbetter1974}
\bibinfo{author}{Leadbetter, M.R.}, \bibinfo{year}{1974}.
\newblock \bibinfo{title}{On extreme values in stationary sequences}.
\newblock \bibinfo{journal}{Z. Wahrscheinlichkeitstheorie verw. Geb.}
  \bibinfo{volume}{28}, \bibinfo{pages}{289--303}.
\bibitem[{Leadbetter(1983)}]{Leadbetter1983}
\bibinfo{author}{Leadbetter, M.R.}, \bibinfo{year}{1983}.
\newblock \bibinfo{title}{Extremes and local dependence in stationary
  sequences}.
\newblock \bibinfo{journal}{Z. Wahrscheinlichkeitstheorie verw. Geb.}
  \bibinfo{volume}{65}, \bibinfo{pages}{291--306}.
\bibitem[{Leeb and P\"{o}tscher(2006)}]{LeebPotscher2006}
\bibinfo{author}{Leeb, H.}, \bibinfo{author}{P\"{o}tscher, B.M.},
  \bibinfo{year}{2006}.
\newblock \bibinfo{title}{Can one estimate the conditional distribution of
  post-model-selection estimators}.
\newblock \bibinfo{journal}{Ann. Statist.} \bibinfo{volume}{34},
  \bibinfo{pages}{2554--2591}.
\bibitem[{Li(2003)}]{Li2003}
\bibinfo{author}{Li, Y.L.}, \bibinfo{year}{2003}.
\newblock \bibinfo{title}{A martingale inequality and large deviations}.
\newblock \bibinfo{journal}{Statist. Probab. Lett.} \bibinfo{volume}{62},
  \bibinfo{pages}{317--321}.
\bibitem[{Liu et~al.(2013)Liu, Xiao and Wu}]{LiuXiaoWu2013}
\bibinfo{author}{Liu, W.}, \bibinfo{author}{Xiao, H.}, \bibinfo{author}{Wu,
  W.B.}, \bibinfo{year}{2013}.
\newblock \bibinfo{title}{Probability and moment inequalities under
  dependence}.
\newblock \bibinfo{journal}{Statist. Sinica} \bibinfo{volume}{23},
  \bibinfo{pages}{257--1272}.
\bibitem[{Ljung and Box(1978)}]{LjungBox1978}
\bibinfo{author}{Ljung, G.M.}, \bibinfo{author}{Box, G.E.P.},
  \bibinfo{year}{1978}.
\newblock \bibinfo{title}{On a measure of lack of fit in time series models}.
\newblock \bibinfo{journal}{Biometrika} \bibinfo{volume}{65},
  \bibinfo{pages}{297--303}.
\bibitem[{Maume-Deschamps(2006)}]{Maume-Deschamps2006}
\bibinfo{author}{Maume-Deschamps, V.}, \bibinfo{year}{2006}.
\newblock \bibinfo{title}{Exponential inequalities and functional estimations
  for weak dependent data; applications to dynamical systems}.
\newblock \bibinfo{journal}{Stochastic Dyn.} \bibinfo{volume}{6},
  \bibinfo{pages}{535--560}.
\bibitem[{McKeague and Qian(2015)}]{McKeague_Qian_2015}
\bibinfo{author}{McKeague, I.}, \bibinfo{author}{Qian, M.},
  \bibinfo{year}{2015}.
\newblock \bibinfo{title}{An adaptive resampling test for detecting the
  presence of significant predictors}.
\newblock \bibinfo{journal}{J. Amer. Statist. Assoc.} \bibinfo{volume}{110},
  \bibinfo{pages}{1422--1433}.
\bibitem[{McLeish(1975)}]{McLeish1975}
\bibinfo{author}{McLeish, D.L.}, \bibinfo{year}{1975}.
\newblock \bibinfo{title}{A maximal inequality and dependent strong laws}.
\newblock \bibinfo{journal}{Ann. Probab.} \bibinfo{volume}{3},
  \bibinfo{pages}{829--839}.
\bibitem[{Merlev{\`e}de and Peligrad(2002)}]{MerlevedePeligrad2002}
\bibinfo{author}{Merlev{\`e}de, F.}, \bibinfo{author}{Peligrad, M.},
  \bibinfo{year}{2002}.
\newblock \bibinfo{title}{On the coupling of dependence random variables and
  applications}, in: \bibinfo{booktitle}{Empical Processes Techniques for
  Dependent Data}. \bibinfo{publisher}{Birkh{\"a}user}, pp.
  \bibinfo{pages}{171--193}.
\bibitem[{Merlev{\`e}de et~al.(2011)Merlev{\`e}de, Peligrad and
  Rio}]{Merlevede_et_al_2011}
\bibinfo{author}{Merlev{\`e}de, F.}, \bibinfo{author}{Peligrad, M.},
  \bibinfo{author}{Rio, E.}, \bibinfo{year}{2011}.
\newblock \bibinfo{title}{Bernstein inequality and moderate deviations for
  weakly dependent sequences}.
\newblock \bibinfo{journal}{Probab. Theory Rel.} \bibinfo{volume}{151},
  \bibinfo{pages}{435--474}.
\newblock \bibinfo{note}{Volume 5}.
\bibitem[{Mies and Steland(2023)}]{MiesSteland2023}
\bibinfo{author}{Mies, F.}, \bibinfo{author}{Steland, A.},
  \bibinfo{year}{2023}.
\newblock \bibinfo{title}{Sequential gaussian approximation for nonstationary
  time series in high dimensions}.
\newblock \bibinfo{journal}{Bernoulli} \bibinfo{volume}{29},
  \bibinfo{pages}{3114--3140}.
\bibitem[{Nagaev(1957)}]{Nagave1957}
\bibinfo{author}{Nagaev, S.V.}, \bibinfo{year}{1957}.
\newblock \bibinfo{title}{Some limit theorems for stationary markov chains}.
\newblock \bibinfo{journal}{Theory Probab. Appl.} \bibinfo{volume}{11},
  \bibinfo{pages}{378--406}.
\bibitem[{Nemirovski(2000)}]{Nemirovski2000}
\bibinfo{author}{Nemirovski, A.S.}, \bibinfo{year}{2000}.
\newblock \bibinfo{title}{Topics in nonparametric statistics}, in:
  \bibinfo{editor}{Emery, M.}, \bibinfo{editor}{Nemirovski, A.},
  \bibinfo{editor}{Voiculescu, D.}, \bibinfo{editor}{Bernard, P.} (Eds.),
  \bibinfo{booktitle}{Ecole d'Ete{\'e} de Probabilitis de Saint-Flour XXVII}.
  \bibinfo{publisher}{Springer}, \bibinfo{address}{Berlin}. volume
  \bibinfo{volume}{1738}, pp. \bibinfo{pages}{87--285}.
\newblock \bibinfo{note}{Lectures Notes on Mathematics}.
\bibitem[{Newey(1991)}]{Newey1991}
\bibinfo{author}{Newey, W.K.}, \bibinfo{year}{1991}.
\newblock \bibinfo{title}{Uniform convergence in probability and stochastic
  equicontinuity}.
\newblock \bibinfo{journal}{Econometrica} \bibinfo{volume}{59},
  \bibinfo{pages}{1161--1167}.
\bibitem[{Peligrad(2002)}]{Peligrad2002}
\bibinfo{author}{Peligrad, M.}, \bibinfo{year}{2002}.
\newblock \bibinfo{title}{Some remarks on coupling of dependent random
  variables}.
\newblock \bibinfo{journal}{Statist. Probab. Lett.} \bibinfo{volume}{60},
  \bibinfo{pages}{201--209}.
\bibitem[{Pinelis(1994)}]{Pinelis1994}
\bibinfo{author}{Pinelis, I.}, \bibinfo{year}{1994}.
\newblock \bibinfo{title}{Optimum bounds for the distributions of martingales
  in banach spaces}.
\newblock \bibinfo{journal}{Ann. Probab.} \bibinfo{volume}{22},
  \bibinfo{pages}{1679--1706}.
\bibitem[{Pollard(1984)}]{Pollard1984}
\bibinfo{author}{Pollard, D.}, \bibinfo{year}{1984}.
\newblock \bibinfo{title}{Convergence of Stochastic Processes}.
\newblock \bibinfo{publisher}{Springer Verlag}, \bibinfo{address}{New York}.
\bibitem[{P{\"o}tscher and Prucha(1989)}]{PotscherPrucha1991}
\bibinfo{author}{P{\"o}tscher, B.M.}, \bibinfo{author}{Prucha, I.R.},
  \bibinfo{year}{1989}.
\newblock \bibinfo{title}{A uniform law of large numbers for dependent and
  heterogeneous data processes}.
\newblock \bibinfo{journal}{Econometrica} \bibinfo{volume}{5},
  \bibinfo{pages}{675--683}.
\bibitem[{Rio(1995)}]{Rio1995}
\bibinfo{author}{Rio, E.}, \bibinfo{year}{1995}.
\newblock \bibinfo{title}{The functional law of the iterated logarithm for
  stationary strongly mixing sequences}.
\newblock \bibinfo{journal}{Ann. Probab.} \bibinfo{volume}{23},
  \bibinfo{pages}{1188--1203}.
\bibitem[{Rio(1996)}]{Rio1996}
\bibinfo{author}{Rio, E.}, \bibinfo{year}{1996}.
\newblock \bibinfo{title}{Sur le theoreme de {B}erry-{E}sseen pour les suites
  faiblement}.
\newblock \bibinfo{journal}{Probab. Theory Rel.} \bibinfo{volume}{104},
  \bibinfo{pages}{255--282}.
\bibitem[{Rio(2017)}]{Rio2017}
\bibinfo{author}{Rio, E.}, \bibinfo{year}{2017}.
\newblock \bibinfo{title}{Asymptotic Theory of Weakly Dependent Random
  Processes}.
\newblock \bibinfo{publisher}{Springer}.
\bibitem[{Samson(2000)}]{Samson2000}
\bibinfo{author}{Samson, P.M.}, \bibinfo{year}{2000}.
\newblock \bibinfo{title}{Concentration of measure inequalities for markov
  chains and $\phi$-mixing processes}.
\newblock \bibinfo{journal}{Ann. Probab.} \bibinfo{volume}{28},
  \bibinfo{pages}{416--461}.
\bibitem[{Talagrand(1995a)}]{Talagrand1995a}
\bibinfo{author}{Talagrand, M.}, \bibinfo{year}{1995}a.
\newblock \bibinfo{title}{Concentration of measure and isoperimetric
  inequalities in product spaces}.
\newblock \bibinfo{journal}{Pub. Math. de l'IHÉS} \bibinfo{volume}{81},
  \bibinfo{pages}{73--205}.
\bibitem[{Talagrand(1995b)}]{Talagrand1995b}
\bibinfo{author}{Talagrand, M.}, \bibinfo{year}{1995}b.
\newblock \bibinfo{title}{The missing factor in hoeffding’s inequalities,}.
\newblock \bibinfo{journal}{Ann. Henri Poincare} .
\bibitem[{Talagrand(2003)}]{Talagrand2003}
\bibinfo{author}{Talagrand, M.}, \bibinfo{year}{2003}.
\newblock \bibinfo{title}{Spin Glasses: A Challenge for Mathematicians: Cavity
  and Mean}.
\newblock \bibinfo{publisher}{Springer}, \bibinfo{address}{Berlin}.
\bibitem[{van~der Vaart and Wellner(1996)}]{vanderVaartWellner1996}
\bibinfo{author}{van~der Vaart, A.}, \bibinfo{author}{Wellner, J.},
  \bibinfo{year}{1996}.
\newblock \bibinfo{title}{Weak Convergence and Empirical Processes}.
\newblock \bibinfo{publisher}{Springer}, \bibinfo{address}{New York}.
\bibitem[{Valenzuela-Dominguez et~al.(2017)Valenzuela-Dominguez, Krebs and
  Franke}]{Valenzuela-Dominguez_et_al_2017}
\bibinfo{author}{Valenzuela-Dominguez, E.}, \bibinfo{author}{Krebs, J.T.N.},
  \bibinfo{author}{Franke, J.E.}, \bibinfo{year}{2017}.
\newblock \bibinfo{title}{A {B}ernstein inequality for spatial lattice
  processes}.
\newblock \bibinfo{note}{ArXiv preprint arXiv:1702.02023}.
\bibitem[{Vershynin(2018)}]{Vershynin2018}
\bibinfo{author}{Vershynin, R.}, \bibinfo{year}{2018}.
\newblock \bibinfo{title}{High-Dimensional Probability}.
\newblock \bibinfo{publisher}{Cambridge University Press},
  \bibinfo{address}{Cambridge, UK}.
\bibitem[{Viennet(1997)}]{Viennet1997}
\bibinfo{author}{Viennet, G.}, \bibinfo{year}{1997}.
\newblock \bibinfo{title}{Inequalities for absolutely regular sequences:
  Application to density estimation}.
\newblock \bibinfo{journal}{Probab. Theory Rel.} \bibinfo{volume}{107},
  \bibinfo{pages}{467--492}.
\bibitem[{Wintenberger(2010)}]{Wintenberger2010}
\bibinfo{author}{Wintenberger, O.}, \bibinfo{year}{2010}.
\newblock \bibinfo{title}{Deviation inequalities for sums of weakly dependent
  time series}.
\newblock \bibinfo{journal}{Electron. Commun. Probab.} \bibinfo{volume}{15},
  \bibinfo{pages}{489--503}.
\bibitem[{Wu(2005)}]{Wu2005}
\bibinfo{author}{Wu, W.B.}, \bibinfo{year}{2005}.
\newblock \bibinfo{title}{Nonlinear system theory: Another look at dependence}.
\newblock \bibinfo{journal}{Proc. Natl. Acad. Sci.} \bibinfo{volume}{102},
  \bibinfo{pages}{14150--14154}.
\bibitem[{Wu(2011)}]{Wu2011}
\bibinfo{author}{Wu, W.B.}, \bibinfo{year}{2011}.
\newblock \bibinfo{title}{Asymptotic theory for stationary processes}.
\newblock \bibinfo{journal}{Statist. Interface} \bibinfo{volume}{0},
  \bibinfo{pages}{1--20}.
\bibitem[{Wu and Min(2005)}]{WuMin2005}
\bibinfo{author}{Wu, W.B.}, \bibinfo{author}{Min, M.}, \bibinfo{year}{2005}.
\newblock \bibinfo{title}{On linear processes with dependent innovations}.
\newblock \bibinfo{journal}{Stochastic Process. Appl.} \bibinfo{volume}{115},
  \bibinfo{pages}{939--958}.
\bibitem[{Wu and Shao(2004)}]{WuShao2004}
\bibinfo{author}{Wu, W.B.}, \bibinfo{author}{Shao, X.}, \bibinfo{year}{2004}.
\newblock \bibinfo{title}{Limit theorems for iterated random functions}.
\newblock \bibinfo{journal}{J. Appl. Probab.} \bibinfo{volume}{41},
  \bibinfo{pages}{425--436}.
\bibitem[{Wu and Shao(2007)}]{WuShao2007}
\bibinfo{author}{Wu, W.B.}, \bibinfo{author}{Shao, X.}, \bibinfo{year}{2007}.
\newblock \bibinfo{title}{A limit theorem for quadratic forms and its
  applications}.
\newblock \bibinfo{journal}{Econometric Theory} \bibinfo{volume}{23},
  \bibinfo{pages}{930--951}.
\bibitem[{Wu and Wu(2016)}]{WuWu2016}
\bibinfo{author}{Wu, W.B.}, \bibinfo{author}{Wu, Y.N.}, \bibinfo{year}{2016}.
\newblock \bibinfo{title}{Performance bounds for parameter estimates of
  high-dimensional linear models with correlated errors}.
\newblock \bibinfo{journal}{Electron. J. Statist.} \bibinfo{volume}{10},
  \bibinfo{pages}{352--379}.

\end{thebibliography}
\singlespacing\setstretch{1} \clearpage
}

\end{document}